\documentclass[10pt]{article}
\usepackage{hyperref}
\usepackage{amssymb}
\usepackage{amsmath}
\input{liemacs10.sty} 
\numberwithin{equation}{section}

\usepackage{color}

\newcommand\ba{{\bf{a}}}
\newcommand\bb{{\bf{b}}}
\newcommand\bc{{\bf{c}}}
\newcommand\be{{\bf{e}}}

\newcommand\bw{{\bf{w}}}
\newcommand\bv{{\bf{v}}}

\newcommand{\jV}{{\rm V}} 

\newcommand{\fco}{\mathop{{\mathfrak{co}}}\nolimits}

\newcommand\bK{{\mathbb K}}

\renewcommand{\bO}{\mathbb O}

\newcommand{\rank}{\mathop{{\rm rk}}\nolimits}
\renewcommand{\rk}{\mathop{{\rm rk}}\nolimits}

\newcommand{\Co}{\mathop{{\rm Co}}\nolimits}

\renewcommand{\phi}{\varphi}

\newcommand{\AdS}{\mathop{{\rm AdS}}\nolimits}
\newcommand{\CFT}{\mathop{{\rm CFT}}\nolimits}

\newcommand{\bHy}[1]{\mathrm{Hyp}^{#1}}

\renewcommand\mlabel{\label}

\begin{document}

\title{Open orbits in causal flag manifolds, \\
modular flows and wedge regions} 
\author{Karl-Hermann Neeb}

\maketitle

\begin{center}
  {\large {\it To the memory of Joseph Anthony Wolf}}
\end{center}

\begin{abstract}
  {We study open orbits of symmetric subgroups of a simple connected
    Lie group $G$ on a causal flag manifold. 
    First we show that a flag manifold $M$ of 
$G$ carries an invariant causal structure  if and only if $G$ is hermitian
  of tube type and   $M$ is the conformal completion  
  of the corresponding simple euclidean Jordan algebra,
  resp., the Shilov boundary of the associated symmetric tube domain.
  We then study open orbits in $M$ under symmetric subgroups,
  also  called causal
  Makarevi\v c spaces, from the perspective of applications in 
Algebraic Quantum Field Theory (AQFT). 
A key motivation is the geometry of corresponding modular flows.
  
  The open orbits are reductive causal symmetric
  spaces, which arise in two flavors:
  compactly causal and non-compactly causal ones.
  In the non-compactly causal case we determine the corresponding Euler elements 
  and their positivity regions. For compactly causal spaces,
  modular flows do not always exist and we determine
  when this is the case. Then the   positivity regions of the modular flows 
are not globally hyperbolic, 
  but these spaces contain other interesting
  globally hyperbolic subsets that can be described
  in terms of the conformally flat Jordan coordinates via Cayley charts. 
  We discuss the Lorentzian case, including de Sitter
  and anti-de Sitter space, in some detail.   

\nin  {\bf Keywords:} causal manifolds, causal symmetric space,
flag manifold, Jordan algebra, Jordan spacetime, $\AdS/\CFT$-correspondence }\\
\end{abstract}
  
\tableofcontents

\section{Introduction}

We call a homogeneous space $M = G/H$ of a Lie group $G$
{\it causal} if there exists a $G$-invariant field $(C_m)_{m \in M}$ of pointed generating
closed convex cones $C_m \subeq T_m(M)$. This is equivalent to the existence
of such a cone for the adjoint action of $H$ on $T_{eH}(M) \cong \g/\fh$.
Recent interest in causal homogeneous spaces in relation with 
representation theory arose from their role as analogs of
spacetime manifolds in the context of
Algebraic Quantum Field Theory (AQFT) 
in the sense of Haag--Kastler, where one considers 
families (also  called nets) of von Neumann algebras $\cM(\cO)$
on a fixed Hilbert space $\cH$,
associated to open subsets $\cO$ in some space-time manifold~$M$ 
(\cite{Ha96}). The hermitian elements of the algebra $\cM(\cO)$ represent
observables  that can be measured in the ``laboratory'' $\cO$.
We consider causal homogeneous spaces as abstract variants of spacetimes
that provide a geometric context in which various features of
nets of operator algebras can be studied in representation theoretic terms.
For some of them, the cone field comes from a
time-oriented Lorentzian metric, but our setting is more general,
focusing on the causal structure, rather than on the Lorentzian metric. 
This allows us to study causality aspects of AQFT
in a highly symmetric context without the need of an invariant Lorentzian form.

An important feature of AQFT that ties the local algebras $\cM(\cO)$
to unitary group representations is that one expects modular groups
of some of these algebras (associated to cyclic separating vectors by
the Tomita--Takesaki Theorem \cite[Thm.~2.5.14]{BR87}) to correspond to
flows on $M$, generated by Lie algebra elements $h \in \g$
(\cite{BDFS00, BMS01, BS04, BGL93, BGL02}).
In \cite{MN24} we have shown that $h$ is a so-called
{\it Euler element}, i.e., $\ad h$ is diagonalizable with eigenvalues in $\{-1,0,1\}$,
so that $\g$ inherits a $3$-grading
\[ \g = \g_1(h) \oplus \g_0(h) \oplus \g_{-1}(h).\]
The positivity regions $W \subeq M$, on which
the vector field $X_M^h(m) = \frac{d}{dt}\big|_{t = 0} \exp(th).m$
is ``timelike'' in the sense that it takes values in the open cone $C_m^\circ$, 
generalize the wedge regions in AQFT and thus deserve particular attention.
We therefore call the flows generated by Euler elements {\it modular}. 

Causal symmetric spaces $M = G/H$ form an important class of causal
homogeneous spaces: Here $G$ carries a non-trivial involution $\tau$,
$H$ is an open subgroup of the group $G^\tau$ of $\tau$-fixed points,
and $\fq := \ker(\1 + \tau) \subeq \g$ can be identified with the tangent space
$T_{eH}(G/H)$ and the causal structure is determined by an
$\Ad(H)$-invariant cone in $C$. We then say that $M$ is
compactly/non-compactly causal if $C^\circ$ consists of elliptic/hyperbolic
elements of $\g$. In \cite{MNO23} we describe a classification
of irreducible non-compactly causal symmetric spaces
in terms of Euler elements.
Irreducible compactly causal spaces are Cartan dual to non-compactly causal
ones, i.e., the Lie algebra of the corresponding group
is $\g^c := \fh + i \fq$ (cf.\ \cite{HO97}).
For the Lorentz group $G = \SO_{1,d}(\R)_e$,
we thus obtain the non-compactly causal de Sitter space $\dS^d$
and its Cartan dual is the compactly causal anti-de Sitter space
$\AdS^d \subeq \R^{2,d-1}$ on which $\SO_{2,d-1}(\R)$ acts transitively.
Wedge regions in irreducible causal symmetric spaces
have been studied in \cite{NO23a, NO23b}.
For representation theoretic constructions of corresponding
(free) nets of local algebras, we refer to
\cite{NO23a, FNO23, MN24}. 
Local nets on causal flag manifolds are studied in detail in the forthcoming
paper \cite{MN25}.

The main goal of the present paper is to provide
more geometric background for this theory, such as a classification of
causal flag manifolds and of the open orbits of symmetric subgroups therein. 
A~key feature of these orbits is that they possess natural flat
coordinates for the causal structure. 
These open orbits are known as causal Makarevi\v c spaces
and have been studied in the 1990s in the context of the Gelfand--Gindikin
program (cf.\ \cite{OO99, Gi92, KOe96, KOe97}).

We start by showing that flag manifolds of simple Lie groups
are causal if and only if $G$ is hermitian and $M$ is the conformal completion
of a simple euclidean Jordan algebra $\jV$ (Section~\ref{sec:2}). 
Then $M$ can be identified with the Shilov boundary of the corresponding
symmetric tube domain (cf.\ \cite{FK94}). 
The open orbits under question
are thus parametrized by involutive Jordan algebra automorphisms~$\alpha$ 
and we review the main facts concerning this classification, following
  work of W.~Bertram and J.~Hilgert (\cite{Be96, BH98, Be00}).
These open orbits are reductive causal symmetric
spaces, thus arising in two flavors: compactly causal
spaces $M^{(\alpha)}$ and non-compactly causal spaces~$M^{(-\alpha)}$.
Here $M^{(\pm \alpha)}$ is the orbit of a connected subgroup $G^{(\pm \alpha)}$,
whose Lie algebra is the fixed point algebra of an involutive
automorphism $\theta_{\pm \alpha}$ of~$\g$. 
In the non-compactly causal case we determine 
the corresponding Euler elements 
 and study their positivity regions~$W \subeq M^{(-\alpha)}$.
    For compactly causal spaces,
    modular flows 
    do not always exist and we determine
  when they do. In this case the
  positivity regions are not 
  globally hyperbolic in the sense that
  their causal intervals are relatively compact subsets
  (cf.\ Section~\ref{subsec:7.6}), 
  but the spaces $M^{(\alpha)}$ 
  contain other interesting globally hyperbolic subsets that can be described
  in terms of the conformally flat Jordan coordinates via Cayley charts.

  The following table lists the simple hermitian Lie algebras of tube type,
  the only non-simple Lie algebra listed is 
  $\so_{2,2}(\R) \cong \so_{1,2}(\R)^{\oplus 2}$, 
  corresponding to the non-simple Jordan algebra 
  $\jV = \R^{1,1} \cong \R \oplus \R$
  (the Minkowski plane, decomposing in lightray coordinates). 

\begin{center} 
  \begin{tabular}{|l|l|l|l|l|l|l|}\hline
\text{Hermitian Lie algebra} & $\g$ &  $\fsp_{2r}(\R)$ & $\su_{r,r}(\C)$  &  $\so^*(4r)$ 
& $\fe_{7(-25)}$  &   $\so_{2,d}(\R)$\\ \hline
\text{Euclidean Jordan algebra} &  $\jV$ &$\Sym_r(\R)$& $\Herm_r(\C)$  & $\Herm_r(\H)$
  & $\Herm_3(\bO)$ &$\R^{1,d-1}\phantom{\Big)}$  \\ \hline
{\rm rank of} $\jV$  & $\rk\jV$ &$r$  &$r$ &$r$  & $3$ & $2$   \\
\hline 
\end{tabular}\\[2mm]
{Table 1: Hermitian Lie algebras 
  of tube type and euclidean Jordan algebras}  
\end{center}

The corresponding flag manifolds $M$ have interesting
geometric interpretations.
For $\g = \so_{2,d}(\R)$, the manifold $M$  is the isotropic
quadric $Q = Q(\R^{2,d})$ in the real projective space $\bP(\R^{2,d})$
(cf.\ Section~\ref{sec:lorentz}), and for 
\[ \Omega := \Omega_{2r} := \pmat{ 0 & \1_r \\ -\1_r & 0} \in M_{2r}(\K), \quad \K = \R, \C, \H,\]
we obtain a uniform realization of the Lie algebras 
$\sp_{2r}(\R), \fu_{r,r}(\C)$ and $\so^*(4r)$ as
\begin{equation}
  \label{eq:fuomega}
  \fu(\Omega,\K^{2r}) := \{ x \in \gl_{2r}(\K) \: x^* \Omega + \Omega x =0\}.
\end{equation}
Then $M$ is the space of
maximal isotropic subspaces $L \subeq \K^{2r}$ with respect to the
skew-hermitian form $\beta(z,w) := z^* \Omega w$ on $\K^{2r}$
(cf.\ Section~\ref{sec:8}).

In \cite{MdR07} Mack and de Riese classify a class of
causal homogeneous spaces they call {\it simple space-time manifolds}
(cf.\ also Remark~\ref{rem:2.6}).
These manifolds are the simply connected coverings of the
causal flag manifolds.
The same class of manifolds has been studied by G\"unaydin in 
\cite{Gu93, Gu75}, where they are called {\it Jordan spacetimes}.

Various aspects of the conformal geometry
  of $Q(\R^{2,d})$ and the submanifolds $M^{(\pm \alpha)}$ 
appear frequently in the Mathematical Physics literature. 
  For instance \cite{BS04} deals with quantum systems on $\AdS^d$,
  where concepts appear that resemble  ``wedge dualities'' in AQFT.
  The unitary group $\U_2(\C)$, with its action of the conformal group
  $\SU_{2,2}(\C)$, as a compactification of Minkowski space $\R^{1,3}$, 
  also appears frequently in conformal field theory (CFT)
  (cf.\ \cite{To19}, \cite{Mo21})
    and  the compactification of anti-de Sitter space $\AdS^4$
  with a $3$-dimensional quadric $Q(\R^{2,3})$ at infinity 
  is the geometric framework of the
  celebrated $\AdS/\CFT$-correspondence
  and its holographic aspects (see \cite{dB01}, \cite{St01},
 \cite{BCW25}, and also 
  \cite{Wo24} for a recent popular
  article relating to quantum cosmology).
  On the more geometric side, \cite{HQR23} studies
  embeddings of Lorentzian manifolds
  into $Q(\R^{2,d})$ to implement $\SO_{2,d}(\R)$ as a group
  of conformal transformations.
  
In the 1990s a principal motivation to study
those spaces  $M^{(\pm \alpha)}$ which are dense in $M$, was that
the conformal group of $M$ can also
be considered as a conformal group of $M^{(\pm \alpha)}$,
acting with singularities, such as fractional linear maps on~$\C$. 
Causal manifolds with open embeddings
into some causal flag manifold $M$ are conformally flat.
In this context Bertram's Liouville Theorem \cite{Be96b} asserts
that local causal diffeomorphisms on open subsets
of $\jV$ always extend to elements of the conformal
group $\Co(\jV)$, whose identity component is~$G$. 

We list a few of these applications. 
Open orbits of Cayley type were used by \'Olafsson and \O{}rsted in
\cite{OO99} to study the Hardy space of the Cayley type spaces~$G/G^h$. 
For spaces $M^{(\alpha)}$ which are groups (considered as causal
symmetric spaces), this technique has also been
used by Gindikin in \cite{Gi92} and by Koufany/\O{}rsted
in \cite{KOe96}. Specifically,
\cite{KOe96} discusses the group $G = \U_{1,1}(\C)$ as a
  $4$-dimensional Lorentzian manifold, embedded into the conformal
  compactification $\U_2(\C)$ of $\R^{1,3}$, which enlarges the
  $7$-dimensional isometry group $(G \times G)/\Delta_{Z(G)}$, acting on $G$, 
  to the $15$-dimensional conformal group
  $G^\flat \cong \SU_{2,2}(\C)$, acting on $\U_{1,1}(\C)$. This is used to study
  the wave representation of $G^\flat$.
  It continues ideas from I.~Segal who
  used the compactification $\U_2(\C)$ of 
  $\Herm_2(\C) \cong \R^{1,3}$   in \cite{Se71} and \cite{Se76}.
  In \cite{KOe97} embeddings
  for the groups $G = \U_{p,q}(\C), \Sp_{2r}(\R)$ and $\SO^*(2\ell)$
  and the corresponding conformal groups
  $G^\flat = \SU_{p+q,p+q}(\C), \Sp_{4r}(\R)$ and $\SO^*(4\ell)$
(see Section~\ref{sec:8})  are used to analyze the Hardy space of $G$ in terms of
  the classical Hardy space of the bounded symmetric domain
  of $G^\flat$. \\

The {\bf structure of this paper} is as follows. In Section~\ref{sec:2}
we show that the only causal flag manifolds $M = G/Q$ are conformal completions
of euclidean Jordan algebras~$\jV$, i.e., $G$ is hermitian of tube type
and $M$ is the minimal flag manifold specified by an Euler element.

In Section~\ref{sec:3} we turn to causal Makarevi\v c spaces
and their  parametrization in terms of involutive automorphisms
$\alpha$ of euclidean Jordan algebras $\jV$.
This section is partly expository, presenting the classification
from the perspective developed in \cite{MNO23} and \cite{NO23a, NO23b},
where causal symmetric spaces are described in terms of Euler elements.
It also contains various refinements concerning the
existence of modular flows on reductive compactly causal symmetric spaces.
In particular, we describe several aspects of 
Jordan involutions, including their close connection
to the classification of modular compactly causal symmetric spaces,
i.e., those with a modular flow having a non-trivial
positivity region (cf.\ Subsection~\ref{subsubsec:3.3.2} and \cite{NO23a}).
Subsection~\ref{subsec:classif-be96} contains a table
with the complete classification of causal Makarevi\v c spaces
and the relevant data. We divide irreducible involutions
of euclidean Jordan algebras into several types:
the flip involution on $\jV \oplus \jV$ (corresponding to Cayley type (C)),
Pierce involutions fixing a Jordan frame (type (P)), and the remaining
involutions of are divided into 
split type (S) ($\rk \jV = \rk \jV^\alpha$) and non-split type (NS)
($\rk \jV = 2\rk \jV^\alpha$).
Section~\ref{sec:inv-jordan} describes
the involutions on the Jordan algebras $\Herm_r(\K)$.

The orbits corresponding to flip involutions
(Cayley type spaces) are discussed in
the short Section~\ref{subsec:flip-ct}. 
In Section~\ref{sec:pierce}, we take a closer look at the properties
of the Pierce involutions and their associated 
non-compactly causal Makarevi\v c spaces~$M^{(-\alpha)}$, which can be realized
as connected components of the set $\jV^\times$ of invertible elements
of~$\jV$. As they surround the cones $\pm\jV_+\subeq \jV$, they are also called
satellite spaces. Here $\theta_{-\alpha}$ is a Cayley type involution on $\g$,
hence   conjugate to $\tau_h$. 
  
  For applications in AQFT it is important to know which
  spaces $M^{(\pm\alpha)} \cong G^{(\pm \alpha)}/H^{(\pm\alpha)}$
  support a modular flow, i.e., $\alpha_t(m) =\exp(th).m$,
  where $h \in \g^{(\alpha)}$ is an Euler element.
This is always the case for the non-compactly
causal spaces $M^{(-\alpha)}$ (Proposition~\ref{prop:cc-alpha},
\cite{NO23b}, \cite{MNO23}),
  but not for the compactly causal spaces $M^{(\alpha)}$.
In Section~\ref{sec:7},
we describe those cases where
  Euler elements $h \in \g^{(\alpha)}$ exist. In general they are not unique,
  but if they exist, there is always one which is also
  Euler in $\g$, and these are all conjugate under $G^{(\alpha)}$.
  These Euler elements define a
  modular structure on $M^{(\alpha)}$, so that the general
  results from \cite{NO23a} on modular compactly causal spaces
  imply that their positivity regions
  (wedge regions) are non-empty (cf.~Theorem~\ref{thm:6.5}).

  Motivated by the importance of Lorentzian manifolds
  as spacetimes in physics, 
we take a detailed look at the Lorentzian case in 
Section~\ref{sec:lorentz}.
Here we connect the conformal compactification of the Jordan algebra
$\jV = \R^{1,d-1}$ (Minkowski space) to concepts from classical
conformal geometry, such as stereographic projections.
The conformal completion $Q = Q(\R^{2,d})$ of Minkowski space $\R^{1,d-1}$
is the only Lorentzian causal flag manifold, but the exceptional isomorphisms
\[ \Sym_2(\R) \cong \R^{1,2}, \quad 
 \Herm_2(\C) \cong \R^{1,3}, \quad 
 \Herm_2(\H) \cong \R^{1,5} \]
also lead to interesting matrix models, in particular for
spacetime dimension $4$, where $Q \cong \U_2(\C)$.
The Makarevi\v c spaces in $Q$ are conformally embedded copies of 
the compactly causal symmetric spaces $(\bS^{d - k} \times \AdS^k)/\{\pm \1\}$ 
and the non-compactly causal spaces 
$\bHy{d - k} \times \dS^k$ 
    for $k = 1, \ldots, d$, where
    $\AdS^k$ is $k$-dimensional Anti-de Sitter space, 
    $\bHy{k}$ is $k$-dimensional hyperbolic space,
    and $\dS^k$ is $k$-dimensional de Sitter space.
    This provides in particular natural flat conformal coordinates on these spaces.

There are three classes of Makarevi\v c spaces
$M^{(\alpha)}$ which are groups, namely
$\Sp_{2n}(\R)$, $\U_{p,q}(\C)$ and $\SO^*(2r)$.
These spaces are discussed in a uniform fashion in Section~\ref{sec:8},
where we explain in particular how the Cayley charts of
these groups appear naturally, 
which implies in particular that it is compatible with the causal structure.
We refer to \cite{LY24} for a recent investigation of Cayley charts of Lie groups
in a much larger context. 

We conclude this paper with a short section listing 
open problems (Section~\ref{sec:perspec}).

  \nin {\bf Acknowledgment:} 
  We thank Wolfang Bertram and Gestur \'Olafsson
  for pointing out references and for discussions
  that eventually led to some results that were formulated
  as open problems in an earlier version of this manuscript.

We also
thank Olaf M\"uller, Klaus Fredenhagen, Karl-Henning Rehren
  and Michael Keyl for illuminating discussions on
  global hyperbolicity of domains in Minkowski space.  \\
  
  \nin {\bf Personal note:} This paper deals with a
  subject matter that was also dear to Joe Wolf, to whom this this volume 
  is dedicated. Open orbits of subgroups in flag manifolds appear
  in many places of his work, see for instance   \cite{Wo69}.
  Here we deal with the specific case of causal flag manifolds. 
  Joe's work had a substantial influence on my own. First of all,
  his beautiful book  on spaces of constant curvature \cite{Wo67}, 
  which I studied as a student,  the importance of his work on flag domains
  in the theory of Lie semigroups and compression semigroups
  \cite{HN93}, his work on direct limit Lie groups \cite{NRW91} and their
  representations \cite{NRW01},  \cite{Wo05}
  and, finally, his work on cycle spaces   \cite{FHW06},
  which relates closely with crown domains of Riemannian symmetric spaces
  and their importance in constructing nets of local algebras in AQFT 
  (cf.~\cite{FNO23, FNO25}).

  The closest contacts I had with Joe were during extended visits
  of Alan Huckleberry's complex analysis group in Bochum in the 1990s,
  where we often had dinner together.
  I frequently asked him for references for ``well-known'' results.
  Then he would typically need only
  a few minutes of digging in his computer account
  to come up with a very precise reference to
  one of his own papers from 20 or more
  years ago. It was amazing to see how efficient he organized his database.
  In 2004 we organized an Oberwolfach meeting on complex analysis and
  representation theory,   together with Alan Huckleberry.

  I only have the fondest memories of Joe as mathematician, a colleague
  and a friend.   With him we lost a mathematician with a very broad view
  of mathematics, firmly grounded in Analysis, Geometry and Algebra. \\

\nin {\bf Notation:}
\begin{itemize}
\item We write $e \in G$ for the identity element in the Lie group~$G$ 
and $G_e$ for its identity component. 
\item $\Inn(\g) = \la e^{\ad \g} \ra$ is the group of inner automorphisms
  of the Lie algebra $\g$ and,
  for a Lie subalgebra $\fs \subeq \g$, we write 
$\Inn_\g(\fs)= \la e^{\ad \fs} \ra \subeq \Inn(\g)$. 
\item $\jV$ denotes a euclidean Jordan algebra,  $\jV_+ \subeq \jV$
  the open positive cone (\cite{FK94}), and for $a, b \in \jV$ with
  $a-b \in \jV_+$, we write
  \begin{equation}
    \label{eq:dab}
    \cD_{a,b} := (a - \jV_+) \cap (b+ \jV_+)
  \end{equation}
  for the corresponding {\it double cone} (causal diamond)
  (\cite{Be18}). We call an open subset $\cO \subeq \jV$ {\it causally convex}
  if $a,b \in \cO$ implies $\cD_{a,b} \subeq \cO$. 
\item $G :=\Co(\jV)_e$ is the identity component of the corresponding
  conformal group $\Co(\jV)$ (\cite{Be96}).
\item $M = \hat\jV \cong G/Q_h$ 
  is the conformal completion of $\jV$ (cf.\ Section~\ref{sec:2}).
\item $G^x$ is the stabilizer of a point $x$ in a $G$-space.
  For the adjoint action, we have in particular
  \[ G^x = \{ g \in G \: \Ad(g)x= x\}.\]
\item $\theta$ is a Cartan involution of the semisimple Lie algebra~$\g$,
  the corresponding involution on~$G$ is denoted $\theta^G$,
  and  $\g = \fk \oplus \fp$ is 
  the associated eigenspace decomposition. 
\item For $x \in \g$ and $\lambda \in \R$, we write 
$\g_\lambda(x) := \ker(\ad x - \lambda \1)$ for the corresponding eigenspace 
in the adjoint representation.
\item $h \in \fp$ is an Euler element defining the $3$-grading
  with $\g_1(h) = \jV$, 
  $\tau_h = e^{\pi i \ad h}$ is the corresponding involution, 
  $\tau_h^G$ the corresponding involution on~$G$,
  and $\kappa_h := e^{\frac{\pi i}{2} \ad h} \in \Aut(\g_\C)$.
\item We write $\cE(\g)$ for the set of (non-central) Euler elements in~$\g$.   
\item For an involutive automorphism $\alpha \in \Aut(\jV)$ of
  the unital euclidean Jordan algebra $\jV$, we write
  $\sigma_\alpha \in \Aut(\g)$ for its canonical extension to $\g$
  with $\sigma_\alpha\res_{\g_1(h)} = \alpha$
  and $\sigma_\alpha^G \in \Aut(G)$ for the corresponding
  automorphism of $G$, given by $\sigma_\alpha^G(g) = \alpha g \alpha^{-1}$
  on~$M$. 
\item $\theta_\alpha = \theta^G \sigma_\alpha^G = \sigma_\alpha^G \theta^G$ and
  $\theta_{-\alpha} := \tau_h^G \theta_\alpha = \theta_\alpha \tau_h^G$
  (involutions on $G$; see \eqref{eq:theta-alpha}).
\item $I_{p,q} = \pmat{\1_p & 0 \\ 0 & -\1_q} = \diag(\1_p, -\1_q)$.  
\item For a causal $G$-manifold $M$ and   $h \in \g$, the positivity region
  in $M$ is
  \[ W_M^+(h) =
  \Big\{ m \in M \: \frac{d}{dt}\Big|_{t = 0} \exp(th).m \in C_m^\circ \Big\}.\]
\end{itemize}

\section{Causal flag manifolds}
\mlabel{sec:2}

In this section we classify causal flag manifolds
(Theorem~\ref{thm:causal-flagman}) and characterize
those which are orientable
(Theorem~\ref{thm:orient}). 
To define flag manifolds, let $G$ be a connected semisimple Lie group
and $x \in \g$ be an element for which $\ad x$ is diagonalizable over
$\R$. Then
\[ \fq_x := \sum_{\lambda \leq 0} \g_\lambda(x), \qquad
  \g_\lambda(x) =\ker(\ad x - \lambda \1), \]
is called the corresponding {\it parabolic subalgebra}, and
\[ Q_x := \{  g \in G \: \Ad(g) \fq_x = \fq_x \} \]
the corresponding {\it parabolic subgroup}. The homogeneous spaces
$G/Q_x$ are called {\it flag manifolds}.

We determine those flag manifolds $G/Q$ of
semisimple connected
Lie groups $G$ that carry a $G$-invariant causal structure.
The main result is Theorem~\ref{thm:causal-flagman}
which implies, in the irreducible case, that $G$ is hermitian of tube 
type and that $Q = Q_h$ is a maximal parabolic subgroup 
corresponding to an Euler element $h$ 
(this is equivalent to the nilradical of $\fq$ being abelian).

We also show in Subsection~\ref{sec:2.2} that, for causal flag manifolds
$M = G/Q_h$, orientability is equivalent to the connectedness of the
stabilizer $\Ad(G)^h$. This is a special feature of
causal flag manifolds because there are also
Euler elements $h$ for which the flag manifold $G/Q_h$ is orientable,
but the centralizer $\Ad(G)^h$ is not connected. 

\subsection{A classification of causal flag manifolds}

The following lemma can be derived from
\cite[Prop.~VII.3.4]{Ne99}. For
the sake of completeness, we include the nice short argument.

\begin{lem} \mlabel{lem:nilpo-invcon} If $\fn$ is  a nilpotent Lie algebra and
  $C \subeq \fn$ is a pointed generating $\Inn(\fn)$-invariant closed
  convex cone, then $\fn$ is abelian. 
\end{lem}

\begin{prf}   
  We start with an abelian ideal $\fa \subeq \fn$
  and a pointed generating invariant cone $C \subeq \fn$.
  For $x \in C$, we then obtain an affine subspace
  $e^{\ad \fa} x = x + [x,\fa] \subeq C$, so that the pointedness of
  $C$ implies that $[x,\fa]= \{0\}$. As $C$ spans $\fn$,
  the ideal $\fa$ is central.

  Now pick $N \geq 2$ maximal such that the descending central series
  defined by $C^1(\fn) = \fn$ and  $C^{k+1}(\fn) := [\fn, C^k(\fn)]$,
  $k \in \N$,   satisfies $C^N(\fn) \not=\{0\}$. We have to show that
  $N = 1$. 
  For $N \geq 3$, we have
  \[ [C^{N-1}(\fn), C^{N-1}(\fn)] \subeq C^{2N-2}(\fn)
    \subeq C^{N+1}(\fn) = \{0\},\]
  so that $C^{N-1}(\fn)$ is an abelian ideal, hence central.
  This leads to the contradiction $C^N(\g) = \{0\}$,
  so that $N \leq 2$. Then the commutator algebra 
  $C^2(\fn) = [\fn,\fn]$ is central in $\fn$. 
  For any $x \in \fn$, the subspace $\fa := \R x + [\fn,\fn]$ now
  is an abelian ideal, hence central, and therefore $[x,\fn] = \{0\}$.
  This shows that $\fn$ is abelian.   
\end{prf}

In the following $\g$ is a non-compact semisimple real Lie algebra,
$\g = \fk \oplus \fp$ is a Cartan decomposition,
$\theta$ the corresponding involution, and  
$\fa \subeq \fp$ a maximal abelian linear subspace.
We write
\[ \g_\alpha := \{  y \in \g \: (\forall x \in \fa) \ [x,y]
  = \alpha(x) y\}, \quad \alpha \in \fa^*, \]
for the corresponding root spaces and
\[ \Sigma := \Sigma(\g,\fa) := \{ \alpha \in \fa^* \:
  0 \not= \alpha, \g_\alpha \not=\{0\}\}\]
for the corresponding root system. 
Note that $\theta(\g_\alpha) = \g_{-\alpha}$ for all $\alpha \in
\Sigma$.

\begin{lem} \mlabel{lem:irred}
  The root spaces $(\g_\alpha)_{\alpha \in \Sigma}$,
  are irreducible $\g_0$-modules.
\end{lem}

\begin{prf} (\cite[Cor.~5.8]{Koo82}, \cite[Thm.~1.9]{Kos10} \begin{footnote}
    {We thank Olivier Guichard for the reference to
      Kostant's paper.}  \end{footnote}  and \cite[Thm.~8.13.3]{Wo67})
  Each element of $\g_\alpha$ is a positive multiple
  of some $x_\alpha \in \g_\alpha$, for which 
  $[x_\alpha, \theta(x_{\alpha})]\in \fa$ 
  is the unique element  $\alpha^\vee\in [\g_\alpha, \g_{-\alpha}] \cap \fa$ with
  $\alpha(\alpha^\vee) = 2$ (\cite[Lemma~5.1]{Koo82}). 
  Then $\fs := \Spann \{ x_\alpha, \theta(x_\alpha), \alpha^\vee\}
  \cong \fsl_2(\R)$ and $\sum_{k \in \Z} \g_{k\alpha}$ is an
  $\fs$-module. Now the representation theory of $\fsl_2(\R)$
  implies that $[x_\alpha, \g_0] = \g_\alpha$. This implies
  the lemma.   
\end{prf}

\begin{lem} \mlabel{lem:root-sum}
  $[\g_\alpha, \g_\beta] = \g_{\alpha + \beta}$ for
  $\alpha, \beta \in \Sigma$ with $\alpha + \beta \not=0$. 
\end{lem}

\begin{prf} (\cite[Thm.~2.3]{Kos10})
  If $\alpha + \beta\not\in \Sigma$, then $\g_{\alpha + \beta} = \{0\}$
  and the assertion is trivial. We may therefore assume that
  $\alpha + \beta \in \Sigma$. 
  Since all root spaces $(\g_\alpha)_{\alpha \in \Sigma}$, are irreducible
  $\g_0$-modules (Lemma~\ref{lem:irred}), it suffices to show that
    $[\g_\alpha, \g_\beta] \not=\{0\}$. Suppose that this is not the
    case. Then
    $[\g_{-\alpha}, \g_{-\beta}] = \theta([\g_\alpha, \g_\beta]) = \{0\}$.
    For the Cartan--Killing form $\kappa$ of $\g$, this leads to 
    \[ \{0\} = \kappa([\g_{-\alpha}, \g_{-\beta}], \g_{\alpha+\beta})
      = \kappa(\g_{-\beta}, [\g_{-\alpha}, \g_{\alpha+\beta}]),\]
    and hence to 
    \begin{equation}
      \label{eq:brak20}
      [\g_{-\alpha}, \g_{\alpha+\beta}] =\{0\}.
    \end{equation}

As in the proof of Lemma~\ref{lem:irred}, we pick $x_\alpha \in \g_\alpha$ with
$[x_\alpha, \theta(x_{\alpha})] =  \alpha^\vee \in \fa$. 
From $[x_\alpha, \g_\beta] = \{0\}$ it follows that
    $\beta(\alpha^\vee) \geq 0$, so that
    $(\beta +\alpha)(\alpha^\vee) \geq 2 > 0$.
    The representation theory of
    $\fs :=\Spann_\R \{x_\alpha, x_{-\alpha}, \alpha^\vee\} \cong \fsl_2(\R)$
    now implies that $[x_{-\alpha}, \g_{\alpha+ \beta}] \not=\{0\}$,
    contradicting \eqref{eq:brak20}.   
\end{prf}

\begin{lem} Let $\g$ be a semisimple real Lie algebra
  and $\fq = \fn \rtimes \fl\subeq \g$ be a parabolic subalgebra
  with nilradical $\fn$ and $\fl$ reductive. 
  Then $\fn$ is abelian if and only if there exists an Euler element
  $h \in \g$ with
  \[ \fn = \g_1(h) \quad \mbox{ with } \quad \fl = \g_0(h).\] 
\end{lem}

\begin{prf} This can be derived from \cite[Lemma~2.2]{RRS92}. 
For the sake of completeness and simplicity,
  we give a direct Lie algebraic proof. 

  If $h\in \g$ is an Euler element, then
  $\g = \g_1(h) + \g_0(h) + \g_{-1}(h)$ is a $3$-grading of $\g$,
  and $\fq :=\g_1(h) + \g_0(h)$ is a parabolic subalgebra with
  abelian nilradical $\g_1(h)$.

  Suppose, conversely, that $\fn$ is abelian. Writing
  $\g$ as a direct sum of simple ideals, $\fq$ is adapted to this
  decomposition, so that we may w.l.o.g.\ assume that $\g$ is simple.
  We choose a maximal $\ad$-diagonalizable subalgebra $\fa \subeq \g$ 
such that 
\[ \fq = \g_0 + \sum_{\alpha \in \Sigma_\fq} \g_\alpha
  \quad \mbox{ for } \quad \Sigma_\fq \subeq \Sigma,\] 
where $\Sigma = \Sigma(\g,\fa)$ is the corresponding system of restricted
roots. It is irreducible because $\g$ is simple.
  With $\Sigma_\fq^0 := \Sigma_\fq \cap - \Sigma_\fq$ and
  $\Sigma_\fq^+ := \Sigma_\fq \setminus \Sigma_\fq^0$, we then have
  $\g = \fn \rtimes \fl$, $\fl = \g_0 + \sum_{\alpha \in \Sigma_\fq^0} \g_\alpha$,
  and 
  \begin{equation}
    \label{eq:n-root}
 \fn = \sum_{\alpha \in \Sigma_\fq^+} \g_\alpha  
    \quad \mbox{ with  } \quad
    (\Sigma_\fq^+ + \Sigma_\fq^+) \cap \Sigma= \eset, 
  \end{equation}
  by Lemma~\ref{lem:root-sum}. 
  Let $\Sigma^+ \subeq \Sigma_\fq$ be a positive system
of $\Sigma$  and $\Pi = \{ \alpha_1, \ldots, \alpha_r\}
  \subeq \Sigma^+$ the subset of simple roots.
  We enumerate them in such  a way that
  $\alpha_1, \ldots, \alpha_k \in \Sigma_\fq^0$ and
  $\alpha_j \in \Sigma_\fq^+$ for $j > k$. 
  Let $\beta = \sum_{j = 1}^r m_j \alpha_j \in \Sigma_\fq^+$
  be the highest root and $\fm \subeq \fn$ be the $\fq$-submodule
  generated by $\g_\beta$. As $\fn$ is abelian, only $\fl$ acts on this space,
  so that all $\fa$-weights of $\fm$ are of the form
  $\gamma = \sum_{j = 1}^r m_j' \alpha_j$ with $m_j' = m_j$ for
  $j > k$. For the lowest $\fa$-weight $\gamma$ of $\fm$,
  there exists a simple root $\alpha_j, j > k,$ with
  $\gamma - \alpha_j \in \Sigma$.   Then
  $\gamma = (\gamma - \alpha_j) + \alpha_j$ 
  implies with \eqref{eq:n-root} that
  $\gamma - \alpha_j \not\in \Sigma_\fq^+$.
  From $\alpha_j = \gamma + (\alpha_j - \gamma)$ we infer that
  also $\gamma - \alpha_j \not\in - \Sigma_\fq^+$, hence that
$\gamma - \alpha_j  \in \Sigma_\fq^0$. The fact that 
$\Sigma_\fq^0$ is generated by $\alpha_1, \ldots, \alpha_k$,
now shows that 
\[ m_j =  1 \quad \mbox{ and } \quad m_i = 0 \quad \mbox{ for } \quad
  k < i \not= j \leq r.\] 
As all $m_j$ are positive, $k = r-1$ and $j = r$.
Hence $\fq$ is a maximal
parabolic subalgebra and the coefficient $m_r$
of $\alpha_r$ in $\beta$ is~$1$.
Then the dual element $h := \alpha_r^* \in \fa$,
specified by $\alpha_j(h) = 0$ for $j < r$ and
$\alpha_r(h) = 1$ is an Euler element
with the asserted properties. 
\end{prf}

\begin{thm} \mlabel{thm:causal-flagman} {\rm(Classification of causal flag manifolds)}
  Let $G$ be a connected semisimple Lie group
  and $Q \subeq G$ be a parabolic subgroup such that $\fq$ contains no
  non-zero ideals of $\g$.   Suppose that the corresponding flag manifold  $G/Q$
  carries a causal structure, i.e., a $G$-invariant field of
  pointed generating closed convex cones
  $C_m \subeq T_m(G/Q)$. Then $\g$ is direct sum of hermitian simple ideals
  and   there exists an Euler element $h \in \g$ such that
  \[ \fq = \fq_h = \g_0(h) + \g_{-1}(h).\]
  If, conversely, this is the case, then $G/Q_h$ is a causal flag manifold. 
\end{thm}

\begin{prf} Let $\fq = \fn \rtimes \fl$ be the Levi decomposition
  of the Lie algebra of $Q$ and $\oline\fn = \theta(\fn)$ the opposite
  nilpotent Lie algebra of $\fn$, so that
$\g = \fn \oplus \fl \oplus \oline\fn$ 
  is a vector space direct sum. For the Cartan--Killing form
  $\kappa$ of $\g$, we have $\kappa(\fn,\fn) = \{0\}$, 
  $\kappa$ is non-degenerate on $\fl \times \fl$
  and $\kappa(\fl, \fn + \oline\fn) = \{0\}$.
  We thus  obtain a natural $\Ad(Q)$-equivariant isomorphism 
  \begin{equation}
    \label{eq:dual-n}
  \fn \to (\g/\fq)^* \cong \fq^\bot, \quad x \mapsto \kappa(x,\cdot). 
  \end{equation}

  Causal structures on $G/Q$ correspond to
  $\Ad(Q)$-invariant pointed generating closed convex cones
  $C \subeq \g/\fq \cong T_{eQ}(G/Q)$. 
  Assume that such a cone exists. Then its dual
  cone
  \[  C^\star := \{ x \in \fn \: \kappa(x,C) \subeq [0, \infty[\} \]
  is a pointed generating $\Ad(Q)$-invariant cone in $\fn$.
  In particular this cone is invariant under $\Inn(\fn)$,
  so that Lemma~\ref{lem:nilpo-invcon} implies that $\fn$ is abelian.

  We claim that the cone $C \subeq \g/\fn \cong \g_1(h)$
  is adapted to the decomposition
  $\g = \g_1 + \cdots + \g_\ell$ of $\g$ into simple ideals.
  By assumption, $\g_0(h)$ contains no non-zero ideals of $\g$, so that
  the Euler element $h$ can be written as $h_1 + \cdots + h_\ell$ with
  (non-zero) Euler elements $h_j \in \g_j \cap \g_0(h)$.
  Then the cone $C$ is invariant under the $1$-parameter group
  induced by $e^{\R \ad h_j}$ in $\g_1(h)$. Further, $\ad h_j$
  has on $\g_1(h)$ only the eigenvalues
  $0$ and $1$ ($1$ on $\g_{j,1}(h)$ and $0$ on $\sum_{k \not= j}
  \g_{k,1}(h)$). This implies that $C$ is adapted to the
  eigenspace decomposition of $\g_1(h)$ under $\ad h_j$
  (Lemma~\ref{lem:ext-eigenvalue}), and thus 
  \[ C = C_1 + \cdots + C_\ell \quad \mbox{ with } \quad
  C_j = C \cap \g_j = p_{\g_j}(C) \]
with pointed generating cones $C_j \subeq \g_{j,1}(h_j)$.
To show that each $\g_j$ is hermitian, we may therefore
assume that $\g$ is simple.

  Let $P_{\rm min} = Z_K(\fa) AN_1 \subeq Q$ be a minimal parabolic subgroup.
  Then $\exp(\fn) \trile N_1$ and the cone $C^\star \subeq \fn$ is
  $\Ad(P_{\rm min})$-invariant.
  Let $x \in \fa$ be such that $\alpha(x) > 0$ for all
  $\alpha \in \Sigma^+$, so that $\beta(x)$, where $\beta$ is
  the highest root, is the maximal eigenvalue of $\ad x$ on $\fn$. 
  Hence the invariance of
  $C^\star \subeq \fn$ under $e^{\R \ad x}$ implies that
  \[ D := C^\star \cap \g_\beta \not=\{0\}\]
  (Lemma~\ref{lem:ext-eigenvalue}).
  The closed convex cone $D$ is fixed pointwise by $\fn$ and invariant
  under the compact group $Z_K(\fa)$. Hence it contains a
  non-zero $Z_K(\fa)$-fixed point~$d$. Then $\R_+ d$ is a $P_{\rm min}$-invariant
  ray in $\g$.  Therefore the Kostant--Vinberg Theorem implies that
  $\g$ contains pointed generating invariant cones 
  (\cite[Thm.~1]{Vi80}), hence that $\g$ is hermitian
  (\cite[Thm.~4]{Vi80}).

  For the converse, we assume that $M \cong M_1 \times \cdots \times M_\ell$
  is a product of flag manifolds corresponding to the simple
  hermitian ideals $\g_j$ and Euler elements $h_j \in \g_j$.
  Let $C_{\g_j} \subeq \g_j$ be  a pointed
  generating invariant closed convex cone. Its invariance under
  $e^{\R \ad h_j}$ implies that
  \[ C_j := p_{\g_1}(C_{\g_j}) = C_{\g_j} \cap \g_1(h) \]
(cf.\ Lemma~\ref{lem:ext-eigenvalue}) is also pointed and generating, and invariant under the
  action of the parabolic subgroup $Q_{h_j}$. We thus  obtain on each
  flag manifold $M_j = G_j/Q_{h_j}$ a causal structure,
  hence also on~$M$.   
\end{prf}

\begin{rem} \mlabel{rem:2.6} In their classification
  of simple space-time manifolds,
  Mack and de Riese obtain the simply connected coverings
  of the causal flag manifolds $G/Q_h$ from Theorem~\ref{thm:causal-flagman} (\cite{MdR07}).
  They start with  the assumptions that
  \begin{itemize}
  \item[(A1)] $M = G/Q_e$ for a parabolic subgroup $Q \subeq G$
  \item[(A2)] For the Levi decomposition $\fq = \fn \rtimes \fl$,
    the nilpotent ideal $\fn$ acts trivially on the tangent space
    $\g/\fq$ of the base point in $M$. This is actually equivalent to
    $\fn$ being abelian.
  \item[(A3)] There exists a global causal order on $M$. This
    implies that $M$ is simply connected and that $\g$ is hermitian.  
  \end{itemize}
 As our arguments show, (A2) follows from
 Lemma~\ref{lem:nilpo-invcon} and
 the existence of the causal structure on $G/Q_e$, hence is redundant.
\end{rem}

\subsection{The corresponding Jordan algebras}

If $\g$ is simple hermitian and
$\Sigma = \Sigma(\g,\fa)$ its restricted root system, then
either $\Sigma$ is of type $C_r$ or~$BC_r$. Moreover,
$\g$ contains an Euler element if and only if
$\Sigma$ is of type $C_r$, i.e., $\g$ is of tube type,
and in this case all Euler elements are $\Inn(\g)$-conjugate 
(cf.\ \cite[Prop.~3.11]{MN21}). In particular $h$ is {\it symmetric} 
in the sense that $-h \in \cO_h := \Inn(\g)h$. 
 Then there exist  
  $e \in \g_1(h)$ and $f =- \theta(e)/2\in \g_{-1}(h)$ such that
  \begin{equation}
    \label{eq:hef}
    [e, f] = h, \quad \mbox{ and also } \quad
    [h,e] = e, \quad [h,f] = -f
  \end{equation}
  (cf.\ Appendix~\ref{app:a.3} for notation related to $\fsl_2(\R)$).
  We then consider on $\jV := \g_1(h)$ the bilinear product
  \begin{equation}
    \label{eq:x*y}
    x * y := [[x,f],y], 
  \end{equation}
  which defines a unital euclidean Jordan algebra
  $(\jV,*,e)$ (cf.~\cite{FK94}). In particular, the set  of squares 
  is a closed pointed generating convex cone
  (Koecher--Vinberg Theorem; see \cite[Thm.~III.2.1]{FK94}).
  Its interior $\jV_+$ is the cone of invertible squares. 
  We refer to Table 1 in the introduction for
  a list of these Lie algebras and the corresponding
  Jordan algebras

  \begin{rem} Although we will not use it explicitly in the following,
    we recall the Jordan products of some Jordan algebras
    that appear below. 
The Jordan algebra structure on the spaces $\Herm_r(\bK)$ is given by
the symmetrized matrix product
\begin{equation}
  \label{eq:jordanmatrix}
  A * B := \frac{1}{2}(AB + BA),
\end{equation}
and on Minkowski space $\R^{1,d-1}$ by
\begin{equation}
  \label{eq:jordan-mink}
  (x_0, \bx) * (y_0, \by)  = (x_0 y_0 + \bx \by, x_0 \by + y_0 \bx)
\end{equation}
(\cite[pp.~25, 31]{FK94}). 
This Jordan algebra embeds into the Clifford algebra
${\rm Cl}(\R^d)$ generated by anticommuting elements
$\be_1, \ldots, \be_d$ with $\be_j^2 = \1$, via
  \[ \iota \: \R^{1,d} \to \CAR(\C^d), \quad
    x = (x_0, \bx) \mapsto x_0 \1 + \sum_{ j = 1}^d x_j \be_j\]
  and the product \eqref{eq:jordanmatrix}. It is also
  called {\it the spin factor}.    
  \end{rem}

\begin{exs} \mlabel{ex:2.6} We describe the natural Euler elements in
  the non-exceptional hermitian Lie algebras from Table~1.

  \nin (a) For $\K = \R,\C, \H$, the Lie algebra 
  \begin{equation}
    \label{eq:urrcb}
\fu(\Omega, \K^{2r})  := \{ x \in \gl_{2r}(\K) \: x^* \Omega+ \Omega x = 0\}
    = \Big \{ \pmat{ a & b \\ c & - a^*} \:
    a \in \gl_r(\K), b,c \in \Herm_r(\K)\Big\} 
  \end{equation}
  contains the Euler element
$ h := \shalf \diag(\1_r, -\1_r)$ 
  and $e = \1$ is a natural identity element in the
  Jordan algebra $\Herm_r(\K)$. This covers three cases at once.
 For $\K = \C$,
  the Lie algebra $\fu(\Omega, \C^{2r}) \cong \fu_{r,r}(\C)
  = \R i\1 \oplus \su_{r,r}(\C)$ is not semisimple,
  for $\K = \R$ we obtain $\sp_{2r}(\R)$, and for
  $\K = \H$ we get $\so^*(4r)$. 

\nin (b) In $\so_{2,d}(\R)$ a natural Euler element is given by
\[ h := E_{1,d+2} + E_{d+2, 1}, \]
where $E_{k\ell} \in M_{d+2}(\R)$ are the matrix units.
\end{exs}

\begin{defn} \mlabel{def:moeb-act}
  The M\"obius group $\PSL_2(\R)$ acts naturally on
  $M \cong G/Q_h$ as the integral subgroup corresponding to the Lie subalgebra
  $\Delta_\fs$  generated by $\{h, e, \theta(e)\}$,
  where $e \in \jV \cong \g_1(h)$ is the Jordan unit
  (cf.\ \eqref{eq:hef}). 
  As $\PSL_2(\R)$ is generated by $\R e$-translations,
  dilations and the map $\theta_\jV(z) = - z^{-1}$.
  The corresponding pulled back action of $\GL_2(\R)_+$ takes
  the form of fractional linear maps 
  \begin{equation}
    \label{eq:fraclin}
    \pmat{a & b \\ c & d}.z = (a z + b e)(c z + d e)^{-1}.
  \end{equation}
  In particular, the action of $\SO_2(\R)$ leads to the ``conformal time
  translations''
  \[ \rho(t).z:=
    \pmat{ \cos(t/2) & \sin(t/2) \\ - \sin(t/2)
      & \cos(t/2)}.z
    = (\cos(t/2) \cdot z + \sin(t/2) \cdot e)
    (- \sin(t/2) \cdot z + \cos(t/2) \cdot e)^{-1}.\]
  We write $z_\fk = \frac{1}{2}(e + \theta(e))\in \Delta_\fs$
  for the element corresponding to $\frac{1}{2}\pmat{ 0 & 1 \\ -1 & 0}
  \in \fsl_2(\R)$,
  so that $\rho(t)= \exp(t z_\fk)$.
 For $t =  \pi$, we thus obtain the Cartan involution 
 \[ \rho(\pi).z = \exp(\pi z_\fk).z = -z^{-1} = \theta_\jV(z),\]
 and for $t = \pi/2$ the real Cayley transform 
 \begin{equation}
   \label{eq:cayley-real-sl2}
   c(z) := \rho(\pi/2).z =  \exp\Big(\frac{\pi}{2} z_\fk\Big).z
   = \pmat{1 & 1 \\ - 1& 1}.z= (e+z)(e-z)^{-1}.
 \end{equation}
  The Euler element $h = \frac{1}{2} \diag(1,-1) \in \fsl_2(\R)$ generates dilations
  \[ \exp(th).z = \pmat{e^{t/2} & 0 \\ 0 & e^{-t/2}}.z = e^t z,\]
  and the Euler element
  $k := \frac{1}{2} \pmat{0 & 1 \\ 1 & 0}$ generates a flow of the form
  \begin{equation}
    \label{eq:k-rot}
 \exp(tk).z = \pmat{\cosh(t/2) & \sinh(t/2) \\ \sinh(t/2) &
      \cosh(t/2)}.z =
    (\cosh(t/2) z + \sinh(t/2) e)(\sinh(t/2) z + \cosh(t/2) e)^{-1}, 
  \end{equation}
fixing~$\pm e$.

  The action \eqref{eq:fraclin} extends naturally to the action of $\PGL_2(\C)$
  on the conformal compactification $M_\C$ of the complex Jordan algebra
  $\jV_\C$. From this complex perspective,
  $\PSL_2(\R)$ preserves the upper tube domain
  $\cT_\jV := \jV +  i \jV_+$, generalizing the complex upper half plane,
  and the subgroup  $\PU_{1,1}(\C)$ preserves the open unit ball 
  $\cD \subeq \jV_\C$ and acts on its Shilov boundary~$\Sigma = \partial^\vee \cD$.
This action also provides the Cayley transform
  \begin{equation}
    \label{eq:cayley-C}
C(z) := \exp\Big(-\frac{\pi i}{2}k\Big).z
    = \pmat{      \cos(\pi/4) & -i \sin(\pi/4) \\ 
      -i\sin(\pi/4) &  \cos(\pi/4)}.z
    = \pmat{ 1 & -i \\ -i & 1}.z
    = \frac{z - i e}{-i z +  e}. 
  \end{equation}
It  satisfies $C(\cT_\jV) = \cD$ (\cite[Thm.~X.4.3]{FK94}), which implies that 
\begin{equation}
  \label{eq:GcintoGD}
 G^c := \exp\Big(-\frac{\pi i}{2}k\Big) G \exp\Big(\frac{\pi i}{2}k\Big)
 \subeq \Inn(\g_\C) 
\end{equation}
is a subgroup preserving $\cD$ and its Shilov boundary~$\Sigma$.
\end{defn}

\subsection{Orientability and connectedness of $\Inn(\g)^h$}
\mlabel{sec:2.2}

In this subsection we address the orientability of causal
flag manifolds. Theorem~\ref{thm:orient} below asserts that it 
is equivalent to the connectedness of the centralizer
$\Inn(\g)^h = \Ad(G)^h$. When this condition is satisfied
is listed explicitly in \cite[Thm.~7.8]{MNO23}.

We assume that $\g$ is simple and also that
$Z(G) = \{e\}$, so that $G \cong \Inn(\g)$.

\begin{lem} \mlabel{lem:2.10} A flag manifold $M = G/Q$ is orientable if and only if
  $\Ad_{\g/\fq}(Q) \subeq \GL(\g/\fq)_+$.
\end{lem}

\begin{prf} If $\Ad_{\g/\fq}(Q) \subeq \GL(\g/\fq)_+$,
  then $M = G/Q$ inherits a $\GL(\g/\fq)_+$-structure, hence
  is orientable. If, conversely, $M$ is orientable,
  then the connectedness of $G$ implies that it acts
  by orientation preserving diffeomorphisms.
  Hence the stabilizer group $Q$ of the base point $eQ$ with
  $T_{eQ}(G/Q) \cong \g/\fq$ maps into 
  $\GL(\g/\fq)_+$, which is equivalent to $\Ad_{\g/\fq}(Q) \subeq \GL(\g/\fq)_+$.
\end{prf}

In view of the preceding lemma, the connectedness of $G^h = \{ g \in G \:
\Ad(g)h = h\}$
(which is equivalent to the connectedness of $Q_h
\cong N \rtimes G^h$) implies that 
$M \cong G/Q_h$ is orientable. But it may happen that
$Q_h$ is not connected and still preserves the orientation on $\jV
\cong \g_1$ (see Example~\ref{ex:sl2n}). So the converse is more subtle.

\begin{thm} \mlabel{thm:orient}
  For an Euler element $h$ in a hermitian simple
  Lie algebra $\g$, the flag manifold $M = G/Q_h$ 
  is orientable if and only if  $\Inn(\g)^h$ is connected. 
The non-orientable flag manifolds correspond to 
  $\g = \so_{2,d}(\R)$ for $d \geq 3$ odd, and
to $\g = \sp_{2r}(\R)$ for $r$ even. 
\end{thm}

\begin{prf} In view of \cite[Thm.~7.8]{MNO23},
  for $\g$ hermitian,  the groups $\Inn(\g)^h$ are not connected for
  $\g = \so_{2,d}(\R)$ if $d \geq 3$ is odd, and for 
  $\g = \sp_{2r}(\R)$ if $r$ is even. In view of Lemma~\ref{lem:2.10}, we
  have to show that in  both cases $M$ is not orientable.

\nin (a) For $\jV = \R^{1,d-1}$, $\g = \so_{2,d}(\R)$
  and $M = (\bS^1 \times \bS^{d-1})/\{\pm \1\}$ with
  $d \geq 3$ odd, we know that $M$ is not orientable
  and the corresponding centerfree group is $G = \SO_{2,d}(\R)_e$.
We write $r_{k\ell}(t) \in \SO_{2,d}(\R)$ for the $t$-rotation 
in the plane spanned by $\be_k$ and $\be_\ell$ for $k <\ell$.
Then 
  \[ g_- = r_{12}(\pi) r_{3,d+2}(\pi)
    = \diag(-1,-1,-1,1,\ldots,1,-1) \in G^h \]
  represents the non-trivial component of $G^h$
  (cf.\ Example~\ref{ex:2.6}(b)).
    As $r_{3,d+2}(\pi) \in \SO_d(\R)_e$ for $d \geq 3$, it acts
  on $\jV \cong \g_1$ with determinant $1$. Moreover
  $r_{12}(\pi)$ acts by $-\id_\jV$ of determinant
  $(-1)^d = -1$. Therefore $\Ad_{\g_1}(Q_h)$ is not
  contained in $\GL(\g_1)_+$ if $d$ is odd. 

\nin (b) For $\jV = \Sym_r(\R)$ and $\g = \sp_{2r}(\R)$,  
  we consider the group $G^* := \Sp_{2r}(\R)$
  with maximal compact subgroup $K^* \cong \U_r(\C)$, $Z(G^*) = \{ \pm \1\}$, 
  and the Euler element $h = \frac{1}{2}\diag(\1_r, - \1_r)$. Then
  \[ (G^*)^h  = \{ \diag(g,g^{-\top}) \: g \in \GL_r(\R)\}
    \cong \GL_r(\R) \]
  is not connected, but if $r$ is even, then $-\1 \in \SL_r(\R)
  \subeq (G^*)^h_e$, so that $\Inn(\g)^h$ is not connected.
  Any $g \in \GL_r(\R)$ with $\det(g) < 0$ represents an element
  in the non-trivial component of $\Inn(\g)^h$.
  Consider $g_- := \diag(-1,1,\ldots, 1).$ 
  Then ${\det}_\jV(\Ad(g_-)) = (-1)^{r-1} = -1$ 
  and therefore $M$ is not orientable.
\end{prf}

The other examples of Euler elements in simple Lie algebras,
where $\Inn(\g)^h$ is not connected (\cite[Thm.~7.8]{MNO23}) 
do not correspond to euclidean Jordan algebras,
but they are also instructive, because they
show that the preceding characterization does {\bf not} extend to
arbitrary Euler elements in simple Lie algebras.
The corresponding flag manifolds $M = G/Q_h$ are
non-causal symmetric $R$-spaces. 

\begin{ex} \mlabel{ex:sl2n} Let $\g = \fsl_{2n}(\R)$ and
  $h = \frac{1}{2}\diag(\1_n,-\1_n)$ with
  $\jV \cong M_n(\R)$. We consider the connected Lie group
  $G^* = \SL_{2n}(\R)$ with $Z(G^*) = \{ \pm \1\}$ and find
  \[ (G^*)^h \cong {\rm S}(\GL_n(\R) \times \GL_n(\R)),\]
  which has $2$ connected components.
  If $n$ is even, then $-\1$ is contained in its identity component,
  so that $\Inn(\g)^h = \Ad(G^*)^h$ is not connected.
  The group $(G^*)^h$ acts on $\jV$ by
$(g_1, g_2).A = g_1 A g_2^{-1},$   so that
  \[  \det(\Ad_\jV(g_1, g_2))
    = \det(g_1)^n \det(g_2)^{-n}  = \det(g_1)^{2n} > 0 \]
  for $(g_1, g_2) \in (G^*)^h$. Therefore
  $M$ is orientable, but $\Inn(\g)^h$ is not connected.   
\end{ex}

\section{Causal Makarevi\v c spaces}
\mlabel{sec:3}

To any involutive automorphism $\alpha$ of the
euclidean Jordan algebra $\jV$, we associate 
two open submanifold $M^{(\pm \alpha)} \subeq M$
as orbits of reductive subgroups $G^{(\pm \alpha)} \subeq G$, 
so that both are causal symmetric spaces for~$G^{(\pm \alpha)}$.
The spaces $M^{(\pm \alpha)}$ are called
causal Makarevi\v c spaces 
and the subgroups $G^{(\pm \alpha)} \subeq G$ are 
the identity components of fixed point groups
of certain involutive automorphisms $\theta_{\pm \alpha}^G$ of~$G$.
The embedding into $M$ provides flat coordinates on
$M^{(\pm \alpha)}$ and can be used to study 
their causal structure and the properties of modular flows. 

In Subsection~\ref{subsec:def-causal-sym-space} we collect
some observations on the causal structure on the spaces
$M^{(\pm \alpha)}$ from the Euler element perspective, as 
developed in \cite{MNO23} and \cite{NO23a, NO23b}.
On a causal symmetric space $M = G/H$, we call a
flow $\alpha_t(m) = \exp(th).m$, generated by an Euler element $h \in \g$,
a modular flow (cf.\ \cite{MN24}).
In \cite{NO23a} positivity regions of modular flows have been studied in modular compactly causal symmetric spaces,
because the existence of an Euler element in $\g$ already
implies the existence of a modular structure
(Subsection~\ref{subsec:scLa}, \cite[Prop.~2.7]{NO23a}),
and this is needed for wedge regions
and positivity regions to be defined.
As this involves Euler elements in $\fh$, it is important
that \cite[Thm.~3.2]{NO23a} provides a classification 
of the $\Inn(\fh)$-orbits of such Euler elements, i.e., 
the different possibilities for modular structures.
For applications to reductive spaces,
we extend this result in Theorem~\ref{thm:7.4}, so that it applies
in particular to the reductive compactly causal symmetric
spaces $M^{(\alpha)} \subeq M$.

In Subsection~\ref{subsec:what-jordan} we explain how
Jordan involutions classify several structures, including
modular compactly causal symmetric spaces, 
i.e., those with a modular flow having a non-trivial
positivity region (cf.\ \cite{NO23a}).
Subsection~\ref{subsec:classif-be96} contains a table
with the complete classification of causal Makarevi\v c spaces
and the relevant data. Finally,
Section~\ref{sec:inv-jordan} provides a description of
the involutions on the Jordan algebras $\Herm_r(\K)$.

\subsection{Causal symmetric Lie algebras}
\mlabel{subsec:scLa}

To describe the basic properties of the causal symmetric spaces
$M^{(\pm \alpha)}$, we first
recall some terminology and observations concerning symmetric spaces and
symmetric Lie algebras:

\begin{itemize}
\item A  {\it symmetric Lie algebra}
is a pair $(\g,\tau)$, where $\g$ is a finite-dimensional real Lie algebra 
and $\tau$ is an involutive automorphism of~$\g$. 
We write 
\begin{equation}
  \label{eq:hq-decomp}
 \g = \fh \oplus \fq \quad \mbox{ with } \quad 
\fh = \g^\tau= \ker(\tau -\1) \quad \mbox{ and } \quad 
\fq = \g^{-\tau}= \ker(\tau +\1).
\end{equation}
\item A {\it causal symmetric Lie algebra}
is a triple $(\g,\tau,C)$, where $(\g,\tau)$ is a symmetric Lie algebra 
and $C \subeq \fq$ is a pointed generating 
closed convex cone, invariant under the group $\Inn_\g(\fh) := \la e^{\ad \fh}\ra
\subeq \Aut(\g)$. 
We call $(\g,\tau,C)$
\begin{itemize}
\item {\it compactly causal~(cc)} if 
$C$ is {\it elliptic} in the sense that, for any $x \in C^\circ$
(the interior of $C$), the 
operator $\ad x$ is semisimple with purely imaginary spectrum. 
\item {\it non-compactly causal (ncc)} if 
$C$ is {\it hyperbolic} in the sense that, for any $x \in C^\circ$, the 
operator $\ad x$ is diagonalizable. 
\end{itemize}
\item For a symmetric Lie algebra $(\g,\tau)$, the
  pair $(\g^c, \tau^c)$ with   
  $\g^c := \fh + i \fq$ and $\tau^c(x + iy)  = x- iy$
  is called the {\it c-dual symmetric Lie algebra}.
  Note that $(\g,\tau,C)$ is non-compactly causal
  if and only if $(\g^c,\tau^c,iC)$
  is compactly causal. 
\item   A~{\it modular causal symmetric Lie algebra} is a 
  quadruple $(\g,\tau, C, h)$, where
  $(\g,\tau,C)$ is a causal symmetric Lie algebra,
  $h \in \g^\tau$ is an Euler element,  
  and the involution $\tau_h$ satisfies $\tau_h(C) = - C$.
  This structure is motivated by the geometric significance
  of the corresponding {\it modular flow} 
  $\alpha_t(m) = \exp(th).m$ on $M$
  (see \cite{NO23a} for details). 
\end{itemize}

\begin{rem} \mlabel{rem:modular-duality}
  $(\g,\tau, C, h)$ is modular if and only if the
  $c$-dual quadruple $(\g^c, \tau^c, i C, h)$ is modular.   
\end{rem}

For the following we recall that a symmetric Lie algebra
$(\g,\tau)$ is called {\it effective} if $\fh = \g^\tau$ contains
no non-zero ideal of $\g$.

\begin{lem} \mlabel{lem:mod-struc}
{\rm (Modular structures on compactly causal symmetric Lie algebras)} 
Let $(\g,\tau,C)$ be an effective reductive compactly causal
symmetric Lie algebra with $C^\circ \cap [\g,\g] \not=\eset$.
If $\g$ contains a non-central Euler element, then
there exist an Euler element $h' \in \fq = \g^\tau$
and a cone $C' \subeq C$ such that
$(\g,\tau,C',h')$ is a modular causal symmetric Lie algebra.
\end{lem}

\begin{prf} (a) First we use the Extension Theorem
    \cite[Thm.~2.4]{NO23a} to find a pointed generating
    $\Inn(\g)$-invariant
  cone $C_\g$ in $\g$ with $-\tau(C_\g) =C_\g$ and $C = C_\g \cap \fq$.
  It follows in particular that $\g$ is quasihermitian, i.e.,
  its simple ideals are either compact or hermitian.
  We write
  $\g = \fz(\g) \oplus \g_h \oplus \fu$ with $\fu$ compact semisimple
  and $\g_h$ a sum of hermitian simple ideals.
  Projecting along the compact semisimple ideal $p_\fu \: \g \to \fz(\g) + \g_h$
  (the fixed point projection
  of the compact group $\Inn(\fu)$), it follows that
  \[ C_\g^\circ \cap (\fz(\g) + \g_h) = p_{\fu}(C_\g^\circ) \not=\eset\]
  (cf.\ Lemma~\ref{lem:coneint})   and likewise 
  \[ C_\g^\circ \cap \g_h = p_{\fu}(C_\g^\circ \cap [\g,\g]) \not=\eset.\]
  Here we use that our assumption implies that
  \[ \eset \not= C^\circ \cap [\g,\g]
    = C_\g^\circ \cap \fq \cap [\g,\g].\]
  
  \nin (b) Let $h_1 \in\g$ be an Euler element.
  Then the ideal $\g_1 \trile \g$ generated by $[h_1,\g]$ has trivial
  center and contains no compact ideal, hence only simple hermitian
  ones with an Euler elements, so that they are  of tube type.
  The $\tau$-invariant ideal $\g_2 := \g_1 +\tau(\g_1)$ also has
  only simple hermitian tube type ideals. We may thus replace
  $h_1$ by an Euler element $h_2 \in [\g,\g]$ generating the ideal $\g_2$.

\nin (c)  Let $\fj \trile \g_2$ be a minimal $\tau$-invariant ideal.
  Then either $\fj$ is simple or a sum of two simple ideals exchanged
  by $\tau$. In the latter case $\fj \cong \fb \oplus \fb$ with
  $\tau$ acting by $\tau(a,b) = (b,a)$ (\cite{NO23a}).
  Any generating Euler element
  in $\fj$ has non-zero components, and all these are conjugate under
  inner automorphisms (\cite[Prop.~3.11]{MN21}).
  So the projection of $h_2$ to $\fj$ is conjugate
  to an element of the form $(x,x) \in \fj^\tau$.
  If $\fj$ is simple, then $\fh = \g^\tau$ contains an Euler element 
  by \cite[Prop.~2.7]{NO23a}. Putting these results on minimal
  invariant ideals together, we see that $h_2$ is conjugate to an
  element of $\g^\tau$, i.e., $\g^\tau$ contains an Euler element $h_3$
  generating~$\g_2$. 

  \nin (d) The involution $\tau_3 := \tau_{h_3}$ commutes with $\tau$. 
  Next we observe that $\g^{-\tau_3} \subeq \g_2$ is contained in
  a sum of hermitian simple
  ideals. Therefore \cite[Prop.~2.7(d)]{NO23a} implies that
the cones $C_\g^{\rm min}$ and $C_\g^{\rm max}$ are $-\tau_{3}$-invariant and
\[ (C_\g^{\rm max})^{-\tau_{3}} = (C_\g^{\rm min})^{-\tau_{3}} = C_\g^{-\tau_{3}}.\]
As $\g_2$ intersects the interior of $C_\g$ and the cone 
$C_\g^{\rm min} \subeq \g_2$ is generating, it
follows with (Lemma~\ref{lem:coneint}) that
\[ \eset \not = (C \cap \g_2^{-\tau_3})^\circ
  =(C_\g \cap \g_2^{-\tau_3})^\circ
  =C_\g^\circ \cap \g_2^{-\tau_3}. \] 
Now 
  \[ C' := C \cap (-\tau_{3}(C))  \subeq \fq \]
  is an $\Inn(\fh)$-invariant pointed cone in $\fq$.
  As it contains $C_\g \cap \g_2^{-\tau_3} \cap \fq
  = C \cap \g_2^{-\tau_3}$, hence interior points of $C_\g$,
  it has non-trivial interior.
    Therefore $(\g,\tau, C', h_3)$ is modular.
\end{prf}
  
\begin{rem}
  The preceding lemma shows that, if $(\g,\tau, C)$ is 
  compactly causal and $\g$ is simple, then 
the existence of a modular structure 
is equivalent to $\g$ being hermitian of tube type
(cf.\ \cite[Prop.~3.11]{MN21}). 

 We refer to \cite[Thm.~5.4]{MNO23} for a discussion
    of this fact from the perspective of non-compactly causal
    symmetric spaces: For the $c$-dual non-compactly causal
    symmetric Lie algebra $(\g^c, \tau^c,i C)$,
    the modularity of $(\g,\tau)$ corresponds to the
    symmetry of the causal Euler element 
    (cf.\ Definition~\ref{def:reductive-ncc} below).
\end{rem}

    \begin{defn} \mlabel{def:reductive-ncc} (The non-compactly causal
      symmetric Lie algebra  associated to an Euler element) \ 
  Let $\g$ be a reductive Lie algebra and
  $h \in \g$ be a (non-central) Euler element. 
We choose a Cartan involution $\theta$ on $\g$ in such a way that 
$\theta\res_{\fz(\g)} = -\id_{\fz(\g)}$ and $\theta(h) = -h$.
Then $\tau := \theta \tau_h$ is an involutive automorphism on $\g$.
We write $\fh := \g^\tau$ and $\fq := \g^{-\tau}$
  for the $\tau$-eigenspaces in $\g$.     Then there exists in $\fq$ 
    a pointed generating $e^{\ad \fh}$-invariant
    cone $C$ containing the Euler element $h$ in its interior 
    (\cite{MNO23, Ol91}).  
We call $(\g, \tau, C)$ a {\it non-compactly causal
      symmetric Lie algebra with causal Euler element $h$.}
  \end{defn}

  The construction in Definition~\ref{def:reductive-ncc}
  is classical for the case where $\g$ is simple
  (\cite{Ol91}),   but the non-simple case has to be treated with some
  caution, as the following examples show. 
  \begin{exs} (a) If $\g = \fz(\g) \oplus [\g,\g]$ and
    $0 \not= h \in \fz(\g)$ is central (a degenerate Euler element), then
    $\tau = \theta$ and $\fq = \fp\supeq \fz(\g)$. Then $\Inn(\fh)$ is a compact group
    fixing $h$, so that there exists a pointed generating
    $\Inn(\fh)$-invariant cone $C \subeq \fq$ containing $h$ in its interior.
    It is easy to see that there is no minimal or maximal cone of this type.
    In fact, pick an affine hyperplane
    $E \subeq \fq$ containing $h + (\fp \cap [\g,\g])$, but not~$0$.
    Then $E$ is $\Inn(\fh)$-invariant, and any $\Inn(\fh)$-invariant
    compact convex neighborhood $U$ of $h$ in $E$ yields a pointed generating
    invariant cone  $C_U := [0,\infty[ \cdot U \subeq \fp$.     

    \nin (b) If $\g = \g_1 \oplus \g_2$ is semisimple, each $\g_j$ is simple
    and $h \in \g_1$, then $\tau\res_{\g_2}$ coincides with
    $\theta\res_{\g_2}$. So
    \[ (\g,\tau) \cong (\g_1, \tau_h \theta) \oplus (\g_2, \theta) \]
    is a direct sum of a non-compactly causal symmetric Lie algebra
    and a Riemannian one, hence non-compactly causal by
    the Cone Extension Lemma (\cite[Lemma~B.2]{MNO23}).
    This lemma applies  because $\tr(\ad x\res_{\fq_1}) = 0$
    for every $x \in \fh_1$.

\nin (c)    A typical example is the Lie algebra
    $\g = \so_{1,d}(\R) \oplus \so_{1,k}(\R)$ and the Euler element $h = (h_0,0)$,
    for which      $M := \dS^d \times \bHy{k}$ is a corresponding non-compactly
    causal symmetric space (of Lorentzian type). 
\end{exs}

The following theorem provides a very convenient
test for modularity of a causal  symmetric Lie algebra:

\begin{thm} \mlabel{thm:7.4}
  {\rm(Modularity test for causal symmetric Lie algebras)} 
For a reductive causal symmetric 
Lie algebra $(\g,\tau)$ with $\fz(\g) \subeq \fq$, the following
assertions hold: 
\begin{itemize}
\item[\rm(a)] If $(\g,\tau)$ is non-compactly causal, then it is modular
  in the sense that there exists a non-central Euler element
$h'$ and a pointed generating $\Inn_\g(\fh)$-invariant closed convex
cone $C \subeq \fq$ satisfying $-\tau_{h'}(C) = C$, 
if and only if the corresponding non-central causal
Euler element $h$ with  $\tau = \theta \tau_h$ is symmetric.    
\item[\rm(b)] A compactly causal symmetric Lie algebra 
$(\g,\tau)$ is modular 
  if and only if its dual, the non-compactly causal symmetric Lie algebra 
  $(\g^c,\tau^c)$, is modular. 
\end{itemize}
\end{thm}

  \begin{prf} (a) Consider the decomposition
    $\g = \fz(\g) \oplus \g_1 \oplus \cdots \oplus \g_\ell$ into simple ideals
    and $h = h_1 + \cdots + h_\ell$ accordingly. Then
    every $\g_j$ is $\tau$-invariant (\cite[proof of Thm.~4.2]{MNO23}) and
    $h$ is symmetric if and only if each $h_j$ is symmetric. 

The case of simple Lie algebras with $h_j \not=0$
is covered by \cite[Thm.~5.4]{MNO23}. It follows that,
$h$ is symmetric if and only if $h_j \not=0$ implies the
existence of modular structures
$(\g_j, \tau\res_{\g_j}, C_j, h_j')$.
On all other simple ideals, $\tau\res_{\g_j}$ is a Cartan involution.

Suppose that $h$ is symmetric and let $h' := \sum_{j = 1}^N h_j' \in \fh$.
This is a non-central Euler element in $\g$
and $-\tau_{h'}(D_1) = D_1$ for $D_1 := \sum_{j = 1}^N C_j$.
Now we use the Cone Extension Lemma (\cite[Lemma~B.2]{MNO23})
to find a pointed generating cone
$D_2 \subeq \fq$ that is invariant under $-\tau_{h'}$ 
and $\Inn_\g(\fh)$. This lemma applies  because $\tr(\ad x\res_{\fq_j}) = 0$
for every $x \in \fh_j$. We conclude that $(\g, \tau, D_2, h')$ is modular.

If, conversely, $(\g, \tau, D, h')$ is modular for some
Euler element $h' \in \fh$, then the same holds for all ideals
$\g_j$ with  $h_j \not=0$, and this implies with 
\cite[Thm.~5.4]{MNO23} that $h_j$ is symmetric. It follows
that $h$ is symmetric.

\nin (b) follows directly from the definition. 
  \end{prf}

  \begin{rem}  \cite[Thm.~3.2]{NO23a} provides more details on
    the classification of $\Inn(\fh)$-orbits of $\g$-Euler elements in $\fh$,
    resp.,     the different classes of modular structures
    of an irreducible compactly causal symmetric Lie algebra $(\g, \tau, C)$:
    \begin{itemize} 
    \item[\rm(G)] For group type spaces ($\tau = \tau_{\rm flip}$)
      with $\g \cong \fh \oplus \fh$ and $\g^\tau = \Delta_\fh$
      the Lie algebra $\fh$ is simple hermitian of tube type
      and there is only one orbit of Euler elements intersecting~$\Delta_\fh$. 
    \item[\rm(C)] For Cayley type ($\tau = \tau_h$)
      of real rank $r$ there are $r+1$ orbits. 
    \item[\rm(S)] For split type ($\rk_\R \fh = \rk_\R \g$),
      there are $2$ orbits, lying in the  same
      orbit of the non-connected group $\Inn(\g)^\tau$. 
    \item[\rm(NS)] For non-split type ($2\rk_\R \fh = \rk_\R \g$),
      there is only one orbit. 
    \end{itemize}    
  \end{rem}

\subsection{General properties of  causal Makarevi\v c spaces}
\mlabel{subsec:def-causal-sym-space}

From here on, we assume that $\g$ is a sum of hermitian simple ideals,
that $h \in \fp$ is an Euler element, and that 
$M = G/Q_h$ is the corresponding causal flag manifold
(Theorem~\ref{thm:causal-flagman}), where $Z(G) = \{e\}$.
Then there exists a
minimal invariant closed convex cone $C_\g \subeq \g$ for which 
\[ C_+ := C_\g \cap \g_1(h) = \oline{\jV_+} \]
is the closed positive cone in the Jordan algebra $\jV := \g_1(h)$,
the tangent space of the base point on~$M$
(cf.\ \eqref{eq:x*y} in Section~\ref{sec:2}).
We also write
\begin{equation}
  \label{eq:hatv}
  \hat\jV = M,
\end{equation}
to emphasize that $M$ is a {\it conformal completion} of $\jV$.

In $\fz(\fk) \cap C_\g$ we choose an element $z_\fk$ for which
$\ad(z_\fk)$ defines a complex structure on $\fp$, so that
\begin{equation}
  \label{eq:thetazk}
  \theta = e^{\pi \ad z_\fk}
\end{equation}
(cf.\ Definition~\ref{def:moeb-act}). 
In particular, $\theta$ is inner, so that
\begin{equation}
  \label{eq:theta-cone}
  \theta(C_\g) = C_\g,
\end{equation}
and thus
\begin{equation}
  \label{eq:c-}
  C_- := - C_\g \cap \g_{-1}(h) = -\theta(C_+).
\end{equation}

Let $\alpha \in \Aut(\jV)\subeq \Co(\jV)$ be an involution, so that
$\alpha$ is fixed by the canonical Cartan involution $\theta^G$
of $\Co(\jV)$, given by $\theta^G(g)= \theta_\jV g \theta_{\jV}
= \exp(\pi z_\fk) g \exp(-\pi z_\fk)$,
where
\begin{equation}
  \label{eq:thetav}
  \theta_\jV(v) = - v^{-1}
\end{equation}
(cf.\ Definition~\ref{def:moeb-act}). 
We write $\sigma_\alpha \in \Aut(\g)$ for the canonical extension
of $\alpha$ to $\g$ with $\sigma_\alpha\res_{\g_1(h)} = \alpha$ fixing $h$,
  and $\sigma_\alpha^G \in \Aut(G)$ for the corresponding
  automorphism of~$G$. 
We thus obtain two involutive automorphisms of $G$: 
\begin{equation}
  \label{eq:theta-alpha}
\theta_\alpha = \theta^G \sigma_\alpha^G = \sigma_\alpha^G \theta^G
  \quad \mbox{ and } \quad
  \theta_{-\alpha} := \tau_h^G \theta_\alpha = \theta_\alpha \tau_h^G,
\end{equation}
where $\tau_h^G$ is the involution of $G$ integrating~$\tau_h$.
We put
\begin{equation}
  \label{eq:malpha}
 G^{(\pm\alpha)} := (G^{\theta_{\pm\alpha}})_e \quad \mbox{ and } \quad 
 M^{(\pm\alpha)} := G^{(\pm\alpha)}.o,
\end{equation}
where we identify $\jV$ with the open subset of $M$ arising from the embedding
\[ \iota_M \: \jV \into M = G/Q_h,\quad x \mapsto \exp(x) Q_h, \]
and consider $o := \iota_M(0)$ as the base point in~$M$. 
Following \cite[\S XI.3]{Be00}, we call $M^{(\pm\alpha)}$ 
{\it causal Makarevi\v c spaces} and write
\[  \g^{(\pm\alpha)} := \L(G^{(\pm \alpha)}) = \g^{\theta_{\pm\alpha}}\]
for the Lie algebra of $G^{(\pm \alpha)}$.

\begin{rem} As the involution $\sigma_\alpha$ fixes $e \in \jV = \g_1(h)$
  and commutes with $\theta$, it also fixes $\theta(e)$, hence the
  Lie algebra $\Delta_\fs \cong \fsl_2(\R)$, generated by these
  elements. In particular, it fixes the element $z_\fk \in \fz(\fk)$
  (cf.\ Definition~\ref{def:moeb-act}). 
\end{rem}

\begin{rem} (a) The classification of the spaces $M^{(\pm \alpha)}$
  can be read from
  Makarevi\v c's paper \cite{Ma73}, which unfortunately contains no proofs.
  He studies open orbits in symmetric R-spaces (he calls them Nagano spaces).
  By \cite{Lo85}, these are precise those flag manifolds $G/Q_h$, where
  $h$ is an Euler element.   Here we are only interested in the causal ones,
  called {\it causal Makarevi\v c spaces} in~\cite{Be96}. 
  As the classification in Section~\ref{sec:2} shows,
  causal flag manifolds are precisely those symmetric
  R-spaces for which the corresponding group $G$   is hermitian.
  We refer to \cite{GKa98} for more related references.

\nin (b)  Compact hermitian symmetric spaces are the symmetric $R$-spaces,
  where $G$ is a complex simple group
  whose Lie algebra contains an Euler element. 
The causal ones can be considered as real forms of these spaces,
corresponding to simple hermitian Lie algebras of tube type,
i.e., to embeddings $M \into M_\C$.
On the level of Jordan algebras, this corresponds to the embedding
of a euclidean Jordan algebra $\jV$ in its complexification~$\jV_\C$. 
  
\nin (c) Makarevi\v c starts with elements 
  $\alpha\in \Str(\jV) \subeq \Co(\jV)$
  satisfying $\theta^G(\alpha) = \alpha^{-1}$ (but
  not necessarily $\alpha^2 = \id_\jV$), which 
  ensures that $\theta_{\pm \alpha}$ are involutions.
However, \cite[Lemma~XI.3.1]{Be00} reduces their classification
to the special case where $\alpha \in \Aut(\jV)$ is an
involution, which we take as our starting point.
We refer to Sections~\ref{sec:inv-jordan} (matrix case) and \ref{subsec:minkowski}  (Minkowski case) for
the classification of involutive automorphisms of  
euclidean Jordan algebras.

\nin (d)  As we shall see below, the classification of
the involutions $\alpha \in \Aut(\jV)$ can be linked
to the classification of irreducible (non-)compactly causal
symmetric spaces (Proposition~\ref{prop:cc-alpha}), hence can also be derived
from the their classification in terms of
Euler elements (\cite{MNO23}, \cite{NO23a}). 
\end{rem}

\begin{prop} \mlabel{prop:3.2} For an involution $\alpha \in \Aut(\jV)$,
  the following assertions hold:
  \begin{itemize}
  \item[\rm(a)] $M^{(\pm\alpha)}$ is an open subset of $M$ and
    an effective causal symmetric space. The corresponding causal symmetric Lie
    algebra is $(\g^{(\pm \alpha)}, \tau_h, C^{(\pm \alpha)})$, where
    \[ C^{(\pm \alpha)} :=
      \{ x \pm \theta_\alpha(x) \: x \in C_+ = \oline{\jV_+} \}.\]
  \item[\rm(b)] The symmetric spaces $M^{(\alpha)}$ and $M^{(-\alpha)}$
    are c-dual to each other. 
  \item[\rm(c)] $M^{(\alpha)}$ is compactly causal with
    $C^{(\alpha)} = C_\g \cap \fq^{(\alpha)}$,
    and  $M^{(-\alpha)}$ is non-compactly causal.
\end{itemize}
\end{prop}

Assertion (c) refines corresponding statements in \cite{Be96}
showing the existence of hyperbolic, resp., elliptic elements in
the cones $C^{(\pm \alpha)}$. 

\begin{prf}   (a) First we show that
  the  symmetric Lie algebras $(\g^{(\pm \alpha)}, \tau_h)$ is  effective.
If $x \in \fh^{(\pm \alpha)}$ is contained 
  in an ideal of $\g^{(\pm \alpha)}$, then $\exp(\R x)$ acts trivially
  on $M^{(\pm \alpha)}$, hence also on $M$, but then
  $x \in \fz(\g) = \{0\}$.
  
For the $h$-eigenspaces $\g_j = \g_j(h)$, we have $\theta_{\pm\alpha}(\g_j)
= \g_{-j}$, so that
\[ \g^{(\pm\alpha)} =\g^{\theta_{\pm\alpha}} =  (\g_0)^\alpha + 
  \{ x \pm \theta_\alpha(x) \: x \in \g_1 \}
  \quad \mbox{ with } \quad
  (\g_0)^\alpha = \{ x \in \g_0 \: \ad x \circ \alpha = \alpha
  \circ \ad x\ \text{ on } \ \g_1\},\] 
which immediately implies that the orbit $M^{(\alpha)} = G^{(\alpha)}.0
\subeq M$ is open. We also see that $(\g^{(\pm\alpha)}, \tau_h)$ is a symmetric
Lie algebra, where
 \[ \g^{(\pm\alpha), \tau_h}
   = \g^{(\pm\alpha)} \cap \g_0 = (\g_0)^\alpha, \]
 so that the orbit $M^{(\alpha)} = G^{(\alpha)}.0$
is a corresponding symmetric space.

\nin (b)   (\cite[Prop.~X.1.3]{Be00}) 
The bracket
\[ [x + \theta_\alpha(x), y + \theta_\alpha(y)]
  = [\theta(\alpha(x)), y] + [x, \theta(\alpha(y))] \]
changes sign when $\alpha$ is replaced by $-\alpha$. Therefore
$(\g^{(\alpha)}, \tau_h)$ is c-dual to 
$(\g^{(-\alpha)}, \tau_h)$.

\nin (c) (cf.~\cite[Prop.~XI.3.2]{Be00}) 
In $\fq^{(\pm\alpha)} = \g^{(\pm\alpha)} \cap \g^{-\tau_h}$, the positive cone
is
\begin{equation}
  \label{eq:calpha}
 C^{(\pm\alpha)} := \{ x \pm \theta_\alpha(x) \: x \in C_+\}.
\end{equation}
As $\alpha(C_+) = C_+$ and $\theta(C_+) = -C_-$ by \eqref{eq:c-},
we obtain
\begin{equation}
  \label{eq:calpha1}
  C^{(\pm\alpha)} = (C_+ \pm \theta(C_+)) \cap \fq^{(\pm \alpha)}
  = (C_+ \mp C_-) \cap \fq^{(\pm \alpha)}. 
\end{equation}
Recall from \cite[Lemma~3.2(ii)]{NOO21} that,
for any choice of $C_\g \subeq \g$ with
$C_\g \cap \g_1(h) = C_+$, we have 
\begin{equation}
  \label{eq:interior-inter}
C_\g \cap \g^{-\tau_h} = C_+ - C_- \quad \mbox{ and } \quad 
  C_\g^\circ \cap \g^{-\tau_h} = C_\g^\circ \cap (\g_1 \oplus \g_{-1})
  = C_+^\circ \oplus - C_-^\circ.
\end{equation}
Now $\fq^{(\alpha)} \cap C_\g  \subeq \g^{-\tau_h} \cap C_\g$ implies
that 
\begin{equation}
  \label{eq:cdag}
  C_\g \cap \fq^{(\alpha)}
  = C_\g^{\theta_\alpha, -\tau_h} =  (C_+ - C_-)^{\theta_\alpha}
  = \{ x + \theta_\alpha(x) \: x \in C_+\} = C^{(\alpha)}.
\end{equation}
Next we observe that $\alpha(e) = e$ implies that 
$e + \theta_\alpha(e) = e + \theta(e) \in C_\g^\circ$, so that
$C^{(\alpha)}$ contains interior points of $C_\g$. 
Thus Lemma~\ref{lem:coneint} implies 
$C^{(\alpha),\circ} = \fq^{(\alpha)} \cap C_\g^\circ$, so that the cone
$C^{(\alpha)}$ is elliptic (because $C_\g$ is elliptic), i.e., that 
$M^{(\alpha)}$ is compactly causal.

For the cone $C^{(-\alpha)} \subeq
\fq^{(-\alpha)}
= \g^{(-\alpha)} \cap \g^{-\tau_h}$, we first observe that the hyperbolic element 
\begin{equation}
  \label{eq:e+theta}
e + \theta_{-\alpha}(e) =   e - \theta_\alpha(e) = e - \theta(e)
\end{equation}
lies in the interior of the hyperbolic cone $C_+ + C_-$
(\cite[\S 7.2]{NO23b}). 
From the fact that 
\[ C^{(-\alpha)} = (\1 - \theta_\alpha)(C_+) = \fq^{(-\alpha)} \cap (C_+ + C_-) \]
contains the hyperbolic element $e - \theta(e) \in C_+^\circ + C_-^\circ$, 
we obtain with Lemma~\ref{lem:coneint} that
\[ C^{(-\alpha), \circ} = \fq^{(-\alpha)} \cap (C^\circ_+ + C_-^\circ).\]
As $C_+ + C_-$ is hyperbolic, the cone $C^{(-\alpha)}$ is hyperbolic.
\end{prf}

\begin{ex} For $\alpha = \id_\jV$, we have $\theta_\alpha = \theta$,
  so that $G^{(\alpha)} = K$, which acts transitively on $M$, so that
  $M^{(\alpha)} = M$, considered as a compactly causal symmetric space
  of $K$.
  
  Accordingly, $\theta_{-\alpha} = \theta \tau_h$ is the
  non-compactly causal involution of $\g$ associated
  to $h$ (\cite{MNO23}), 
  \[ \g^{(-\id_\jV)} = \fh_\fk + \fq_\fp \]
  and $M^{(-\id_\jV)}$ is the open unit ball
  $\cD_\jV \subeq \jV$, in its Harish--Chandra embedding
  (cf.~\cite{MNO24}).
\end{ex}

For $v \in \jV$, let $\tau_v(x) := x + v$ denote the translation
by $v$, write $\theta_\jV(x) = -x^{-1}$ for the Cartan involution, 
and
\[ \hat\tau_v(x) := (\theta_\jV \tau_v \theta_\jV)(x)
  = -(v - x^{-1})^{-1} = (x^{-1} - v)^{-1}.\]
For $x,v \in \jV$, we consider the {\it  Bergman operator} of
the Jordan algebra $\jV$, defined for $\hat\tau_v(x) \in \jV$ by
\[ B(x,v)
  := \dd (\hat\tau_v\tau_x)(0)^{-1}\]
(\cite[\S 1.2]{Ber98}).
We recall from \cite[Prop.~1.2.1]{Ber98}
that $B$ extends to a polynomial operator-valued function
on $\jV$.
According to \cite[Lemma~1.2.5]{Ber98}, we have
\begin{equation}
  \label{eq:invcrit}
 B(x,y) \in \GL(\jV)
  \qquad \Leftrightarrow \qquad
  \hat\tau_{-x}(-y) \in \jV 
  \qquad \Leftrightarrow \qquad
  \exp(-\theta(x)).(-y) \in  \jV.
\end{equation}

\begin{lem} \mlabel{lem:B-inv-on-D} If $\jV$ is a euclidean Jordan algebra
  and $\cD_\jV \subeq \jV$ its open unit ball with respect to the
  spectral norm, then, for $x, y \in \cD_\jV$, the operator
  $B(x,y)$ is invertible. 
\end{lem}

\begin{prf} 
\cite[Lemma~XII.1.9]{Ne99} implies that,
for $z,w$ in the unit ball in $\jV_\C$,
the bounded symmetric domain of the conformal group $\Co(\jV)$,
we have $\exp(w).z \in \jV_\C$. For $z,w \in \cD_\jV$,
it follows that $\exp(w).z \in \jV$.  Therefore $B(z,w)$ is
invertible by \eqref{eq:invcrit}.   
\end{prf}

\begin{prop} \mlabel{Dv-inc} We always have $\cD_\jV \subeq M^{(\pm \alpha)}$.
\end{prop}

\begin{prf}
If $\alpha \in \Aut(\jV)$ is an involution, then
$\pm \alpha(\cD_\jV) = \cD_\jV$ implies that, for all
$x \in \cD_\jV$, the operator $B(x,\pm \alpha x)$ is invertible. 
As \cite[Thm.~2.1.1(ii)]{Ber98} implies that 
the intersection $M^{(\pm\alpha)} \cap \jV$ is a union of
connected components of the set of all $x \in \jV$ with
$B(x,\pm\alpha x)$ is invertible,
Lemma~\ref{lem:B-inv-on-D} shows that $\cD_\jV \subeq M^{(\pm \alpha)}$.
\end{prf}

The following proposition provides a sufficient criterion for modularity
which is easily verified in concrete cases. For (c), we recall
the concept of a positivity region. 

\begin{defn} \mlabel{def:posreg}
  If $M$ is a causal $G$-manifold with cone field $(C_m)_{m \in M}$, and
  $h \in \g$, we define the corresponding {\it positivity region} by 
 \[ W_M^+(h) =
  \Big\{ m \in M \: \frac{d}{dt}\Big|_{t = 0} \exp(th).m \in C_m^\circ \Big\}.\] 
\end{defn}

\begin{prop} \mlabel{prop:aut-modul} {\rm(Automatic modularity)}
  If $h' \in \fg^{(\pm \alpha)}$ is an {\bf Euler element of $\g$},
  then the following assertions hold:
  \begin{itemize}
  \item[\rm(a)]  If $h' \in \fh^{(\pm \alpha)}$,
    then $(\g^{(\pm \alpha)}, \tau_h, C^{(\pm \alpha)}, h')$ 
    is a modular causal symmetric Lie algebra, 
    i.e., $\tau_{h'}(C^{(\pm \alpha)}) = -  C^{(\pm \alpha)}$.
  \item[\rm(b)] If $h' \in [\g^{(\alpha)}, \g^{(\alpha)}]$, then
    $h'$ is conjugate under $\Inn(\g^{(\alpha)})$ to some $h'' \in \fh^{(\alpha)}$, 
    and the quadruple $(\g^{(\pm\alpha)}, \tau_h, C^{(\pm\alpha)}, h'')$ is modular. 
   \item[\rm(c)]     $W^+_{M^{(\pm \alpha)}}(h') \not=\eset$. 
  \end{itemize}
  \end{prop}

\begin{prf} (a) We consider the linear isomorphism
$\psi \:  \fq^{(\pm \alpha)} \to \jV = \g_1(h),$ 
corresponding to the tangent map in $0$ for the orbit map
$G^{(\pm\alpha)} \to M, g \mapsto g.o$. It corresponds to the
linear  projection $\g \to \g_1(h)$, restricted to $\g^{(\pm\alpha)}$.
Since $h' \in \fh^{(\pm \alpha)} \subeq \g_0(h)$, the involution
  $\tau_{h'}$ commutes with $\ad h$, hence leaves $\g_1(h)$ invariant, 
  and $\psi$ is $\tau_{h'}$-equivariant.
 
  As $h'$ is an Euler element of $\g$, we have 
  $-\tau_{h'}(C_\g^{\rm max}) = C_\g^{\rm max}$ (\cite{MNO23}).
  Thus $\oline{\jV_+} = C_\g^{\rm max} \cap \g_1(h)$
  implies that $-\tau_{h'}(\jV_+) = \jV_+$,
  and finally the equivariance of $\psi$ entails that
\[ -\tau_{h'}(C^{(\pm \alpha)}) = -\tau_{h'}(\psi^{-1}(\oline{\jV_+}))
  = C^{(\pm \alpha)}.\]

\nin (b)   As $\g^{(\alpha)}$ is quasihermitian
(all simple ideals of $\g^{(\alpha)}$ are hermitian of tube type),
$h'$ is contained in
  the commutator algebra if and only if it is symmetric. 
  We claim that it is conjugate under $\Inn(\g^{(\alpha)})$ to
  an Euler element in $\fh^{(\alpha)}$. 
  As $(\g^{(\alpha)}, \tau_h)$
  is reductive, effective and compactly causal
  (Proposition~\ref{prop:3.2}), 
  this follows from part (c) in the proof of Lemma~\ref{lem:mod-struc}.

If $h' \in  \g^{(\alpha)}$ is not contained in the
    commutator algebra, then $h' = h'_1 + h'_2$ with $h'_1$ central
    and $h'_2$ an Euler element of $\g^{(\alpha)}$, contained in the
    commutator algebra. Applying the preceding argument to $h'_2$,
    we see that $h'_2$ is $\Inn(\g^{(\alpha)})$-conjugate to an element
    of $\fh^{(\alpha)}$. The center of $\g^{(\alpha)}$ is contained
    in $\fq^{(\alpha)}$, so that $h'$ is not conjugate to an element of
    $\fh^{(\alpha)}$.    
  
\nin (c) 
That the positivity regions are non-empty for modular
causal symmetric Lie algebras follows
from \cite[Thm.~6.5]{NO23a} for general compactly causal symmetric Lie
algebras and from \cite[Thm.~7.1]{NO23b} for 
reductive non-compactly causal ones. 
\end{prf}

\subsection{What involutions in $\Aut(\jV)$ classify}
\mlabel{subsec:what-jordan}

We have seen above that involutive automorphisms $\alpha$ of euclidean
Jordan algebras $\jV$ lead to open submanifolds
$M^{(\pm \alpha)} \subeq M$ which are symmetric spaces of the subgroups
$G^{(\pm \alpha)} \subeq G = \Co(\jV)_e$. To understand the
diversity of this family, it is crucial to observe that 
these automorphisms classify various types of geometric objects.
In this subsection we briefly describe the connections to
\begin{itemize}
\item simple real Jordan algebras $\jV^c = \jV^\alpha  + i \jV^{-\alpha}$
  with the same complexification as $\jV$ (\cite{NK19}), and 
\item irreducible modular compactly causal symmetric spaces
  (cf.~\cite{NO23a}).
\end{itemize}

\subsubsection{Simple real Jordan algebras}

Let $(\jV, \alpha)$ be a euclidean Jordan algebra and
$\alpha \in \Aut(\jV)$ an involution. The pair $(\jV, \alpha)$ is called
an {\it involutive Jordan algebra}. 
We further assume that
$(\jV,\alpha)$ is {\it indecomposable}, i.e., $\{0\}$ and $\jV$
are the only $\alpha$-invariant ideals. 
Indecomposable involutive Jordan algebras
are in one-to-one correspondence with simple real Jordan algebras by 
the assignment
\[ (\jV, \alpha) \mapsto \jV^c := \jV^\alpha + i \jV^{-\alpha}\subeq \jV_\C. \]
For a fixed $\jV$, we thus obtain
all simple real Jordan algebras with the same complexification~$\jV_\C$ 
(see Table 2 and \cite{Hel69}, \cite{FG96}, \cite{BH98}, \cite{NK19}).
The symmetric $R$-space $M^c$ associated to~$\jV^c$
  can be realized as a totally real submanifold 
  of the complex hermitian $R$-space $M_\C = G_\C/Q_h$
  ($Q_h \subeq G_\C$ the corresponding parabolic subgroup;
  cf.~Section~\ref{sec:2}), 
  which is a hermitian symmetric space of compact type.

  \begin{rem} \mlabel{rem:realtube}
    (a)     If $\jV$ is not simple and $(\jV,\alpha)$ is indecomposable,
    then $(\jV,\alpha)$ is equivalent to the involutive Jordan algebra 
  $(\jV^\alpha \oplus \jV^\alpha, \alpha)$ with
  $\alpha(v,w) = (w,v)$ (flip involution), and then
  \[ \jV^c = \{ (x + iy, x- iy) \: x, y \in \jV^\alpha \}
    \cong (\jV^\alpha)_\C.\]
This case corresponds to the situations, where $\jV^c$ is a complex
  simple Jordan algebra and the involution is antilinear. 
    
\nin (b) The antiholomorphic 
  involution $\oline\tau_\alpha(z) := - \oline{\alpha(z)}$
  on the tube domain  $\cT_\jV = \jV + i \jV_+$
  has the fixed point set
  \[ \jV^{-\alpha} \oplus i \jV_+^\alpha = \cT_\jV \cap i \jV^c, \] 
  which is a {\it real tube domain}. 
  \end{rem}

  \vspace{5mm} 
The following table lists the data related to the classification 
of simple real non-complex Jordan algebras. 
For further details we refer to \cite[\S 1.5]{BH98} and \cite{NK19};
see also \cite{Hel69} and \cite{FG96}. If $\alpha$ is
not of Pierce type (P), we distinguish 
split type (S) ($\rk \jV = \rk \jV^\alpha$) and non-split type (NS)
($\rk \jV = 2\rk \jV^\alpha$). We write (S1), (S2), \ldots
and (NS1), (NS2), \ldots etc. for the sake
of easier reference below.\\

\hspace{-12mm}
\begin{tabular}{||l|l|l|l|l|l||}\hline
 \text{type} & $\jV$   & $\jV_\C$ & $\alpha$\ on\ $\jV$ 
& $\jV^\alpha$ &  $\jV^c$ \\ 
\hline\hline 
(P)  &  $\Sym_{p+q}(\R)$ &  $\Sym_{p+q}(\C)$ &
 $I_{p,q} x I_{p,q}$ & $\Sym_p(\R) \oplus \Sym_q(\R)$ & $\Sym_{p,q}(\R)$  \\
(NS1)  & $\Sym_{2s}(\R)$ & $\Sym_{2s}(\C)$ & 
$\Omega x\Omega^{-1}, \Omega^2 = -\1$ & $\Herm_s(\C)$ & $\Sym_{2s}(\C) \cap M_s(\H)$  \\
 & & & & & $\cong \Aherm_s(\H)$\\
\hline 
  (P)  & $\Herm_{p+q}(\C)$
& $M_{p+q}(\C)$   &  $I_{p,q} x I_{p,q}$ & $\Herm_p(\C) \oplus \Herm_q(\C)$  &  $\Herm_{p,q}(\C)$  \\

(S1) &$\Herm_r(\C)$ &  $M_r(\C)$ & 
\phantom{\Big(}
$\oline x$ &  $\Sym_r(\R)$ &   $M_r(\R)$  \\
  (NS2)& $\Herm_{2s}(\C)$ & $M_{2s}(\C)$ & $\Omega x\Omega^{-1}, \Omega^2 = -\1$  &$\Herm_s(\H)$
               &  $M_s(\H)$  \\
\hline
  (P) & $\Herm_{p+q}(\H)$ & $\Skew_{2(p+q)}(\C)$
 &$I_{p,q} x I_{p,q}$ &  
$\Herm_p(\H) \oplus \Herm_q(\H)$  &$\Herm_{p,q}(\H)$   \\
(S2) & $\Herm_r(\H)$ & $\Skew_{2r}(\C)$ 
& $(I x_{ij} I^{-1})_{ij}$ &  $\Herm_r(\C)$ &  $\Skew_{2r}(\R)$\\
\hline
  (P)&   $\Herm_3(\bO)$  & $\Herm_3(\bO)_\C$
& $I_{1,2} x I_{1,2}$ & 
$\R \oplus \Herm_2(\bO)$  & $\Herm_{1,2}(\bO)$ \\ 
 (S3) &  $\Herm_3(\bO)$ & $\Herm_3(\bO)_\C$  &$\alpha_0$  & 
 $\Herm_3(\H)$& $\Herm_3(\bO_{\rm split})$  \\ 
\hline
(P) &  $\R^{1,d-1}$ & $\C^d$  & 
$\diag(1, 1, -\1_{d-2})$  & $\R^{1,1}$ & $\R^{d-1,1}$ \phantom{\Big(}\\
(S4)& $\R^{1,p+q-1}$ & $\C^{p,q}$  & 
$\diag(1, - \1_{p-1}, \1_q)$  & $\R^{1,q}$&  $\R^{p,q}, p,q > 1$  \\
(NS3) & $\R^{1,d-1}$ & $\C^d$ & 
$\diag(1, - \1_{d-1})$ & $\R$  &  $\R^{d,0}$ \\
\hline
\hline
\end{tabular} \\ 
\begin{center}{Table 2: Simple real non-complex unital Jordan algebras 
    $\jV^c$ with euclidean dual $\jV$}
\end{center}

\subsubsection{Irreducible modular compactly causal symmetric spaces} 
\mlabel{subsubsec:3.3.2}

Irreducible compactly causal modular symmetric Lie algebras 
$(\g, \tau, C,h)$ (cf.\ Section~\ref{subsec:scLa})  arise in two types (\cite{NO23a}):
\begin{itemize}
\item {\it Group type:} $(\g, \tau) \cong (\fh \oplus \fh,\tau_{\rm flip})$, 
$(\fh, C_\fh)$ is a hermitian simple Lie algebra, and
$C_\fh \subeq \fh$ a pointed generating invariant cone.
Then $h = (h_0, h_0)$ holds for an Euler element $h_0 \in \fh$,
so that the Lie algebra $\fh$ must be simple hermitian of tube type
(\cite[Prop.~3.11]{MN21}). 
\item {\it Simple type:} 
$\g$ is simple hermitian of tube type and $\tau$ an
involution satisfying $-\tau(C_\g) = C_\g$ for a pointed generating
invariant cone $C_\g \subeq \g$ satisfying $C_\g \cap \fq = C$.
If $\tau$ commutes with the Cartan involution $\theta$, this is
equivalent to  $\tau(z_\fk) = - z_\fk$.
\end{itemize}

The following lemma follows from \cite[Prop.~3.12]{Oeh22}.
\begin{lem} \mlabel{lem:b.1} Let $\g$ be simple hermitian of tube type, 
$h \in \g$ an Euler element, and $\jV := \g_1(h)$ the corresponding euclidean 
Jordan algebra. For every involutive automorphisms $\alpha \in \Aut(\jV)$,  
there exists a unique 
automorphism $\sigma_\alpha \in \Aut(\g)$ with $\sigma_\alpha\res_\jV = \alpha$,
and then $(\g, \tau_h \sigma_\alpha, C_\g^{-\tau_h \sigma_\alpha}, h)$ is modular
compactly causal.
Conversely,  every modular compactly causal Lie algebra  
of simple type is of this form.
\end{lem}

The following proposition connects the pairs $(\jV,\alpha)$
to irreducible compactly causal symmetric Lie algebras, 
which have been classified in terms of Euler elements
in \cite{MNO23} and \cite{NO23a}. 

\begin{prop} \mlabel{prop:cc-alpha}
  Let $\g$ be simple hermitian of tube type, $h \in \cE(\g)$
  and $\jV = \g_1(h)$.
  \begin{itemize}
\item[\rm(a)] The assignment
  \[ (\jV, \alpha) \mapsto
    (\g, \tau_h \sigma_\alpha, C_\g \cap \g^{-\tau_h \sigma_\alpha}, h) \]
  from indecomposable involutive Jordan algebras to modular 
  causal symmetric Lie algebras,
  yields all modular compactly causal symmetric Lie algebras
  of simple type.
\item[\rm(b)] Every involutive automorphisms $\tau \in \Aut(\g)$
  commuting with $\theta$ and satisfying $\tau(z_\fk) = - z_\fk$
  and $\tau(h) = h$ is of the form $\tau_h \sigma_\alpha$,
  hence in particular compactly causal.
\item[\rm(c)]   For every involution $\alpha \in \Aut(\jV)$, the involution
    $\theta_{-\alpha}$ defines a modular compactly causal symmetric
  Lie algebra $(\g,\theta_{-\alpha}, C, k)$,
  where $C := C_\g \cap \fq^{(-\alpha)}$, and 
  \begin{equation}
    \label{eq:kh}
    k := e^{-\frac{\pi}{2} \ad z_\fk} h  
  \end{equation}
  is the corresponding Euler element   fixed by $\theta_{-\alpha}
  = \tau_k \sigma_\alpha$.
\item[\rm(d)] The element $k\in \g^{(-\alpha)}$ from
  {\rm Definition~\ref{def:moeb-act}} is a causal
  Euler element for the non-compactly causal symmetric Lie algebra
  $(\g^{(-\alpha)}, \tau_h, C^{(-\alpha)})$. In particular,
  $\eset \not= W_{M^{(-\alpha)}}^+(k) \subeq W_M^+(k)$. 
  \end{itemize}
\end{prop}

\begin{prf} (a) Lemma~\ref{lem:b.1} implies that the involution 
$\tau$ of a compactly causal modular symmetric Lie algebra 
$(\g,\tau, C,h)$ of simple type is of the form $\tau_h \sigma_\alpha$. 

\nin (b) follows from   \cite[Prop.~3.12(a)]{Oeh22}. 

\nin (c)   For the Euler element $k \in \g$
  we derive from $\theta = e^{-\pi \ad z_\fk}$ that 
  \begin{align*}
    \theta_\alpha
    &= \theta \sigma_\alpha = (\theta \tau_h) (\tau_h \sigma_\alpha)
      = e^{-\frac{\pi}{2} \ad z_\fk} \tau_he^{\frac{\pi}{2} \ad z_\fk} 
      (\tau_h \sigma_\alpha)
      = \tau_k \tau_h \sigma_\alpha =  \tau_h \tau_k\sigma_\alpha.       
  \end{align*}
  Accordingly,
  \begin{equation}
    \label{eq:cc-invol}
    \theta_{-\alpha} 
    = \tau_h \theta_\alpha=   \tau_k\sigma_\alpha.
  \end{equation}
  In view of (b), this   is a modular compactly
  causal involution on $\g$ because $\sigma_\alpha(k) = k$
  implies $\theta_{-\alpha}(k) = k$ 
  and $\theta_{-\alpha}(C_\g) =~- C_\g$ follows from $\tau_k(C_\g) = - C_\g$.

\nin (d)
    Equation \eqref{eq:e+theta} in
    the proof of Proposition~\ref{prop:3.2} shows that the element
  $e - \theta(e) \in C^{(-\alpha)}$ is hyperbolic.
  Comparing with the subalgebra
  $\fs \cong \fsl_2(\R)^{\oplus r}$ corresponding to the strongly orthogonal
  roots (see \eqref{eq:delts-fs} below), it follows that the Euler element $k$
  from \eqref{eq:kh} is a causal Euler element in~$\g^{(-\alpha)}$.
  In particular, it is contained in the diagonal subalgebra 
  $\fsl_2(\R) \cong \Delta_\fs \subeq \fs$,
  spanned by $h$, $e$ and~$\theta(e)$. 
\end{prf}

The following observation refines 
Proposition~\ref{prop:cc-alpha}(d).

\begin{prop} \mlabel{prop:wm+k} The positivity region of the Euler element
  $k = e^{-\frac{\pi}{2} \ad z_\fk} h$ in $M$ is
    the open unit ball for the spectral norm in $\jV$: 
  \[ W_M^+(k) = \cD_\jV= (e - \jV_+) \cap(-e+ \jV_+).\]
\end{prop}

\begin{prf} We derive from $W_M^+(h) = \jV_+$ (\cite[Lemma~2.7]{MN25})
  and \eqref{eq:kh} 
  that $W_M^+(k) = \exp(-\frac{\pi}{2} z_\fk).\jV_+.$ 
From \eqref{eq:cayley-real-sl2} in Definition~\ref{def:moeb-act},
we know that
$\exp(\frac{\pi}{2}z_\fk)$ acts on $\jV \subeq M$ as the Cayley
transform $c(x) = \frac{e +x}{e-x}$ 
of order $4$ with $c^2(x) = - x^{-1} = \theta^M(x)$. Therefore 
$\exp(-\frac{\pi}{2}z_\fk)$ acts by 
$c^{-1}(x) = \frac{x-e}{x+e}.$ 
As $c^{-1}(\R_+) = (-1,1)$, this shows that 
$W_M^+(k) = c^{-1}(\jV_+) = \cD_\jV.$
\end{prf}

We conclude that
\begin{equation}
  \label{eq:wnotempty}
  \eset\not= W_{M^{(\pm\alpha)}}^+(k) = M^{(\pm\alpha)} \cap \cD_\jV \subeq \jV
\end{equation}
(cf.\ Proposition~\ref{prop:cc-alpha}(d)).

\begin{thm} \mlabel{thm:pos-dom-ncc}
  For a non-compactly causal space $M^{(-\alpha)} \subeq M$
  and the causal Euler elemenet $k$ from
  {\rm Definition~\ref{def:moeb-act}}, the  positivity region is
  \[ W_{M^{(-\alpha)}}^+(k) = \cD_\jV  = (e - \jV_+) \cap(-e+ \jV_+).\]
 The subgroup $G^k_e$ acts transitively on this domain. 
\end{thm}

\begin{prf} We have seen in Proposition~\ref{prop:cc-alpha}(d)
  that $k$ is a causal Euler element for $M^{(-\alpha)}$,
  and in  Proposition~\ref{prop:wm+k} that
  $W_M^+(k) = \cD_\jV$. Hence the assertion follows from
  Proposition~\ref{Dv-inc}.

  That the subgroup $G^k_e$ acts transitively on $\cD_\jV$
  follows from the fact that $k$ is conjugate to $h$
  (cf.~\eqref{eq:kh}) 
  and that $G^h_e$ acts transitively on the open cone~$W_M^+(h)
  = \jV_+$ (\cite[Lemma~2.7]{MN25}). 
\end{prf}
  
The preceding theorem shows that the positivity
region is a so-called real bounded symmetric domain.
For a detailed analysis of these domains, we refer to
\cite{OSt20}.  We also refer to Table 1 in
  \cite{GKO03}, which lists those crown domains $\Xi$ of $G/K$
  that are bounded symmetric domains of larger groups. 

  The homogeneity of wedge regions seems to be curiously
  related to the discussion of
  one-parameter groups of isometries  in the recent paper
  \cite{GHZ25}.

\begin{rem} \mlabel{rem:3.7} 
  (a) If $(\g, \tau, h)$ is irreducible compactly causal modular, then
  $\tau =\tau_h \sigma_\alpha$ (Lemma~\ref{lem:b.1}),
  $\tau$ fixes $h$   and $\g_1^\tau = \jV^{-\alpha}$, so that
  \[ \fq = \jV^\alpha \oplus \g_0^{-\alpha} \oplus \g_{-1}^\alpha.\]

  \nin (b) The involution $\theta_\alpha = \theta \sigma_\alpha$
  fixes $C_{\g^{(\alpha)}} = C_\g \cap \g^{(\alpha)}$ and $z_\fk$,
  whereas $\tau_h(C_\g) = - C_\g$.
  Therefore
  $(\g^{(\alpha)}, \tau_h, C_{\g^{(\alpha)}})$ is 
  compactly causal (cf.~Proposition~\ref{prop:3.2}).

\nin (c) For the compactly causal involution $\tau := \tau_h \sigma_\alpha$, 
the c-dual symmetric Lie algebra $(\g^c, \tau^c, h)$ is
  non-compactly causal modular (Proposition~\ref{prop:3.2}),
  so that $\g^c\cong \fco(\jV^c)$ is the conformal Lie algebra of
  the simple non-euclidean Jordan algebra $\jV^c$.

\nin (d)  We have the causal symmetric Lie algebras
  \[ (\g^{(\pm\alpha)}, \tau_h\res_{\g^{(\pm\alpha)}}, C^{(\pm\alpha)})
    \quad \mbox{ with } \quad
     \fh^{(\alpha)} = \fh^{(-\alpha)}
      = \g_0(h) \cap \g^{(\pm \alpha)}, \qquad 
     \fq^{(\pm\alpha)} = \g^{(\pm\alpha)} \cap \g^{-\tau_h}.\] 
    The Euler element $k$ is contained in $\g^{(-\alpha)}$ and
    causal for $M^{(-\alpha)}$ (Proposition~\ref{prop:cc-alpha}(d)),
    but it is not contained in $\fh^{(\alpha)}$.

    Both Euler elements $h$ and $k$ are contained in the diagonal subalgebra
    \begin{equation}
      \label{eq:delts-fs}
      \fsl_2(\R) \cong \Delta_\fs \subeq\fs \cong \fsl_2(\R)^{\oplus r}
    \end{equation}
    (see Definition~\ref{def:moeb-act}). 
    On $\Delta_\fs$ the involution $\theta_\alpha$ fixes $z_\fk$,
    hence is a Cartan
    involution and
    \[ \g^{(\alpha)} \cap \Delta_\fs  = \R z_\fk = \so_2(\R).\]
    The involution $\theta_{-\alpha} = \tau_h \theta_\alpha = \tau_k \sigma_\alpha$
    (cf.\ \eqref{eq:cc-invol}) satisfies 
    \[ \g^{(-\alpha)} \cap \Delta_\fs  = \R k = \so_{1,1}(\R).\] 
    In particular, $\g^{(\alpha)}$ does not always contain Euler elements.
    This phenomenon is discussed in detail in Section~\ref{sec:7}.
\end{rem}

\subsection{The classification of causal Makarevi\v c spaces}
\mlabel{subsec:classif-be96}

Proposition~\ref{prop:cc-alpha}(c) connects irreducible modular
   compactly causal symmetric Lie algebras, resp., quadruples 
   $(\g, \tau_h \sigma_\alpha, C, h)$ with
   the non-compactly causal symmetric spaces
   $M^{(-\alpha)}$ of $G^{(-\alpha)}$.  
  In particular, it exhibits the subalgebras $\g^{(-\alpha)}$ as the fixed
  point algebras  of  compactly causal involutions on~$\g$,
which have been classified in terms of Euler elements
in \cite{MNO23} and \cite{NO23a}. 

  We distinguish the following types (see also
  the introduction of  \cite{NO23a}): 
  \begin{itemize}
  \item[\rm(C)] {\bf Flip involutions} $\alpha$ correspond to group type spaces:
    $(\g \oplus \g)^{(-\alpha)} \cong \g$,
    and these correspond to Cayley type spaces
    $M^{(-\alpha)} \cong G/G^h$, embedded in the product space 
    $(\jV \oplus \jV)\,  \hat{}\, \cong \hat\jV \times \hat\jV$
    (see Subsection~\ref{subsec:flip-ct}). 
  \item[\rm(P)] {\bf Pierce involutions} $\alpha$
    correspond to satellites $M^{(-\alpha)} \subeq \jV^\times$
    (see Section~\ref{sec:pierce} and in particular
    \eqref{eq:pierce-inv-alpha} for details).
  \item[\rm(S)] Spaces of {\bf spit type}, i.e., 
    $\tau \not=\tau_h$ and 
    $\rk_\R \g^{(-\alpha)} = \rk_\R \g$:
\item[\rm(NS)] Space of {\bf non-split type}, i.e., 
  $\tau \not=\tau_h$, $\rk_\R \g^{(-\alpha)} =\frac{ \rk_\R \g}{2}$: 
\end{itemize}

Combining the classification from \cite{NO23a} with the material
from \cite{Be96} and \cite{BH98}, leads to the following table.
Here $\fco(\jV)$ denotes the conformal Lie algebra corresponding to
a  pair $(\jV,\alpha)$, where $\alpha$ may also be a flip involution
and $\jV$ not simple. For the correspondence between
$\jV$ and $\fco(\jV)$, we refer to Table~1.
For Pierce type we also write $(p,q) = (r-j,j)$ in the table. 
Although the symmetric spaces $M^{(\pm \alpha)}$ are determined
up to covering by the symmetric pair
$(\g^{(\pm \alpha)}, \fh^{(\pm \alpha)})$, we sometimes
list a global representative. \\
\hspace{-20mm}
\begin{tabular}{||l|l|l|l|l|l|l|l||}\hline
  Type & $\fco(\jV)$ & $\g^{(\alpha)}$  & $\g^{(-\alpha)}$
  & $\fh^{(\alpha)} = \fh^{(-\alpha)}$ & $M^{(\alpha)}$ & $M^{(-\alpha)}$   \\ 
\hline 
\hline 
$\Sym_r(\R)$\phantom{\Big (} & $\sp_{2r}(\R)$ && &&& \\
\hline   
 (C) &$\sp_{2r}(\R)^{\oplus 2}$ & $\sp_{2r}(\R)$ & $\sp_{2r}(\R)$  &  $\gl_r(\R)$ && \\
  (P) & $\sp_{2r}(\R)$ & $\fu_{p,q}(\C)$ & $\gl_r(\R)$  &  $\fo_{p,q}(\R)$ &
$\U_{p,q}(\C)/\OO_{p,q}(\R)$  &   $\Sym_r(\R)^\times_{(p,q)}$ \\
  (NS1) & $\sp_{4s}(\R)$& $\sp_{2s}(\R)^{\oplus 2}$ & $\sp_{2s}(\C)$
  &  $\sp_{2s}(\R)$ & $\Sp_{2s}(\R)$           & \\
  \hline
\hline 
$\Herm_r(\C)$\phantom{\Big (} & $\su_{r,r}(\C)$ && &&& \\
\hline   
  (C)  & $\su_{r,r}(\C)^{\oplus 2}$& $\su_{r,r}(\C)$ & $\su_{r,r}(\C)$
  &  $\fsl_r(\C)\oplus \R \1$&& \\
(P) & $\su_{r,r}(\C)$ &
        $\fu_{p,q}(\C)^{\oplus 2}/\R (i\1,i\1)$ & $\fsl_r(\C)+ \R \1$
  &  $\su_{p,q}(\C)$&
 $\U_{p,q}(\C)$ & $\Herm_r(\C)^\times_{(p,q)}$ \\
  (S1)  & $\su_{r,r}(\C)$ & $\so^*(2r)$ & $\so_{r,r}(\R)$  &  $\so_r(\C)$ &
$\SO^*(2r)/\SO_r(\C)$     &\\
  (NS2)  &$\su_{2r,2r}(\C)$ & $\sp_{4r}(\R)$  & $\fu_{r,r}(\bH)$  &  $\sp_{2r}(\C)$
 & $\Sp_{4r}(\R)/\Sp_{2r}(\C)$ &  \\
  \hline
\hline 
$\Herm_r(\H)$\phantom{\Big (} & $\so^*(4r)$ && &&& \\
  \hline
 (C) &$\so^*(4r)^{\oplus 2}$ & $\so^*(4r)$ & $\so^*(4r)$  &  $\gl_r(\H)$& & \\
(P) & $\so^*(4r)$& $\fu_{2p,2q}(\C)$ & $\gl_r(\H)$ 
                                     &  $\fu_{p,q}(\H)$
      & $\U_{2p,2q}(\C)/\U_{p,q}(\H)$ & $\Herm_r(\H)^\times_{(p,q)}$\\
(S2) &$\so^*(4r)$ & $\so^*(2r)^{\oplus 2}$ & $\so_{2r}(\C)$  &  $\so^*(2r)$ &
 $\SO^*(2r)$        &\\
  \hline
\hline 
$\Herm_3(\bO)$ \phantom{\Big (}& $\fe_{7(-25)}$ &&& & &\\
\hline   
  (C)  & $\phantom{\Big(}$$\fe_{7(-25)}^{\oplus 2} $& $\fe_{7(-25)}$  & $\fe_{7(-25)}$
  &  $\fe_{6(-26)} \oplus \R$ &&
 \\
  (S3) & $\fe_{7(-25)}$ & $\su_{6,2}(\C)$ & $\fsl_4(\H)$ & $\fu_{3,1}(\H)$ &
$\SU_{6,2}(\C)/\U_{3,1}(\H)$  & \\
  (P), $j=0$ & $\fe_{7(-25)}$ & $\fe_6^c \oplus \R$ &$\fe_{6(-26)} \oplus \R $ & $\fe_4^c$ &&
$\Herm_r(\bO)_+$  \\
(P), $j = 1$ & $\fe_{7(-25)}$ & $\fe_{6(-14)} \oplus \R$  & $\fe_{6(-26)} \oplus \R$ & $\ff_{4(-20)}$
&& $\Herm_r(\bO)^\times_{(2,1)}$  \\
  \hline   
\hline   
$\R^{1,d-1}$ \phantom{\Big (}& $\so_{2,d}(\R)$ && &&& \\
\hline   
(C)  & $\so_{2,d}(\R)^{\oplus 2}$ & $\so_{2,d}(\R)$ & $\so_{2,d}(\R)$  &  $\R \oplus \so_{1,d-1}(\R)$
 && \\
 (S4) &$\so_{2,d}(\R)$ &$\so_{2,p-1}(\R)$ & $\so_{1,p}(\R)$ 
                             &$\so_{1,p-1}(\R)$
 & $(\AdS^{p} \times \bS^{q})/\{\pm\1\}$ & $\dS^{p} \times \bHy{q}$ \\
 \quad $p+q=d$ &&  $\quad \times \so_{q+1}(\R)$ &$\quad\times \so_{1,q}(\R)$
  &  $\quad \times \so_q(\R)$
&  &  \\  
(NS3), $q=0$ &$\so_{2,d}(\R)$ &$\so_{2,d-1}(\R)$ & $\so_{1,d}(\R)$ 
                             &$\so_{1,d-1}(\R)$
 & $\AdS^{d}$ & $\dS^d$  \\
(P), $j=0$      &$\so_{2,d}(\R)$  & $\R \oplus\so_d(\R)$ & $\R \oplus  \so_{1,d-1}(\R)$
     &  $\so_{d-1}(\R)$ 
   & $(\bS^1 \times \bS^{d-1})/\{\pm\1\}$ & $\R \times \bHy{d-1}$ \\
(P), $j =1$  &$\so_{2,d}(\R)$  & $\R \oplus \so_{2,d-2}(\R)$
 & $\R \oplus \so_{1,d-1}(\R)$  &  $\so_{1,d-2}(\R)$
                                         &$(\AdS^{d-1} \times \bS^1)/\{\pm\1\}$
     &$\dS^{d-1} \times \R$ \\
\hline\hline 
\end{tabular}
\\ \begin{center}{Table 3: Causal Makarevi\v c spaces}\end{center}

\subsection{Involutions on Jordan matrix algebras} 
\mlabel{sec:inv-jordan}

In this subsection we describe the involutive automorphisms 
on the matrix Jordan algebras $\Herm_r(\K)$
following \cite[\S 1.5]{BH98}.
Our grouping according to Pierce type  (P),
split type (S) and non-split type~(NS) corresponds to the
classification of irreducible compactly causal symmetric
spaces (Table~1 in \cite{NO23a}), which is based on 
the correspondence from Proposition~\ref{prop:cc-alpha}.
Our types correspond to types I, II and III in \cite[p.~8]{FG96}. 
Involutions on  the Minkowski spaces $\R^{1,d-1}$ are
discussed in Subsection~\ref{subsec:minkowski}  below.\\

\nin {\bf Pierce involutions (P):}
On the Jordan algebras $\jV$ of the form
$\Herm_r(\K)$ with $\K = \R,\C,\bH, \bO$, we have the involutions
\begin{equation}
  \label{eq:alphaq}
  \alpha_q(z) := I_{p,q} z I_{p,q}, \quad p + q = r,
\end{equation}
with 
\[ \jV^\alpha = \Herm_p(\K) \oplus \Herm_q(\K) \quad\mbox{ and } \quad 
  \jV^{-\alpha} \cong  M_{p,q}(\K)\]
(cf.\ Example~\ref{ex:4.1}).
Note that $\alpha_0 = \alpha_r = \id_\jV$. 
Here $\jV^c$ is the space of hermitian matrices with respect to the 
hermitian form $\beta(z,w) = z^* I_{p,q} w$
of signature $(p,q)$. The corresponding spaces $M^{(-\alpha)}$
are given by $\GL_r(\K).I_{p,q} \cong \GL_r(\K)/\U_{p,q}(\K)$
and discussed in some detail in Proposition~\ref{prop:pierce-concrete} below.\\

\nin {\bf Split type involutions (S):}
These are the automorphisms of $\Herm_r(\K)$ induced by
automorphisms $\alpha_0$ of $\K$, applied entry-wise: 
\[ \alpha((x_{jk})) = (\alpha_0(x_{jk})). \] 
So the involutions (S1,2,3) below correspond to the equal rank
embeddings of euclidean Jordan algebras: 
\[ \Sym_r(\R) \into \Herm_r(\C), \quad 
  \Herm_r(\C) \into \Herm_r(\H), \quad 
  \Herm_3(\H) \into \Herm_3(\bO). \]

\nin {\bf(S1)} $\jV = \Herm_r(\C)$:
Here $\alpha(z) = \oline z$ with 
$\jV^\alpha = \Sym_r(\R)$, $\jV^{-\alpha} = i \Skew_r(\R)$
and $\jV^c = M_r(\R)$. \\ 

\nin {\bf (S2)} $\jV = \Herm_r(\H)$: 
Consider the automorphism $\alpha_0(z) = IzI^{-1}$
of $\H$ whose fixed point space is $\R\1 + \R I \cong \C$.
We then have $\jV^\alpha = \Herm_r(\C)$.\\ 

  \nin {\bf(S3)} $\jV = \Herm_3(\bO)$: 
Here $\alpha_0$ is the involution on $\bO$ with
$\bO^{\alpha_0} = \H$. Then 
$\jV^\alpha = \Herm_3(\H)$, $\jV^{-\alpha}_\C \cong  M_{6,2}(\C)$, 
$\jV^c \cong \Herm_3(\bO_{\rm split})$. \\

\nin {\bf Non-split type involutions (NS):}
Here $r = 2s$ and $\alpha$ is related to conjugation with
\[ \Omega := \pmat{ 0 & -\1_{2s} \\ \1_{2s} & 0}.\]
The involutions of non-split type (NS1,2) correspond to the
half-rank embeddings
\[ \Herm_s(\C) \into \Sym_{2s}(\R), \quad 
  \Herm_s(\H) \into \Herm_{2s}(\C).\] 

\nin {\bf(NS1)} $\jV  = \Sym_{2s}(\R)$:
Here $\alpha(x) = \Omega x\Omega^{-1}$, 
$\jV^\alpha = \Herm_s(\C)$, $\jV^{-\alpha} = \Sym_s(\C)$. 
For $x = \pmat{a & b \\ b^\top & d} \in \Sym_{2s}(\C) = \jV_\C$, 
the antilinear extension $\oline\alpha$ of $\alpha$ to $\jV_\C$ takes the form 
\[ \oline\alpha\pmat{a & b \\ c & d}
= \pmat{0 & -\1 \\ \1 & 0}\pmat{\oline a & \oline b \\ \oline c& \oline d}\pmat{0 & \1 \\ -\1 & 0}
= \pmat{\oline d & -\oline c \\ - \oline b & \oline a}.\] 
Therefore $x \in \jV^c$ is equivalent to  $\oline d = a$ and $c = - \oline b$, 
which is equivalent to $x \in M_s(\H)$. We thus obtain 
$ \jV^c = \Sym_{2s}(\C) \cap M_s(\H).$ 

\nin {\bf(NS2)} On $\jV = \Herm_{2s}(\C)$ with $\jV_\C = M_{2s}(\C)$,
we consider the involution
$\alpha(z) =  \Omega \oline z \Omega^{-1}$ with 
$\jV^\alpha = \Herm_s(\H)$, $\jV^{-\alpha} = \Aherm_s(\H)$, 
and $\jV^c = M_s(\H)$.

\section{Flip involutions and Cayley type spaces} 
\mlabel{subsec:flip-ct}

In this brief section we discuss flip involutions
and show that the corresponding open orbits
correspond to Cayley type spaces $G/G^h$, embedded in $M_d := M \times M$
as open $G$-orbits (cf.~\cite{Ko94}, \cite{BN04}).

If $\alpha$ is an involution on a euclidean Jordan algebra
and there exists no non-trivial $\alpha$-invariant ideal, then either
the Jordan algebra is simple or consists of two simple ideals flipped by
$\alpha$. In the latter case, there exists a simple Jordan algebra $\jV$
such that we have a situation of the type 
\begin{equation}
  \label{eq:flip-jordan}
\jV_d  = \jV \oplus \jV \quad \mbox{ with }  \quad 
  \alpha(v,w) = (w,v).  
\end{equation}
These involutions correspond to so-called Cayley type spaces
which are briefly discussed in this subsection. 
Let $G_d := G \times G$ be the double of $G = \Co(\jV)_e$. 

    On the group $G_d \cong \Co(\jV_d)_e$, we obtain 
$\sigma_\alpha(g_1, g_2) = (g_2,g_1)$, and hence the automorphism
    \[ \theta_\alpha(g_1, g_2) = (\theta^G(g_2), \theta^G(g_1),\]
which leads to 
    \[ G_d^{(\alpha)} = \Gamma(\theta^G) \cong G = \Co(\jV)_e.\]
    Further, $\g_d^{(\alpha), \tau_h} = \{ (x,\theta(x)) \: x \in \g_0(h) \}
      \cong \g_0(h),$     so that
    \[ M_d^{(\alpha)} \cong G/G^h \cong \cO_h \] 
    is a Cayley type space.
    
    Writing $\tau_{\rm ncc} = \theta \tau_h$ for the non-compactly causal
    involution defined by $h$ (cf.\ Definition~\ref{def:reductive-ncc}), we obtain
    from $\theta_{-\alpha} = \tau_h^{G_d} \theta_\alpha
    = \tau_{\rm ncc}^{G_d} \sigma_\alpha$ the relation 
    \[ \theta_{-\alpha}(g_1, g_2) = (\tau_{\rm ncc}^G(g_2), \tau_{\rm ncc}^G(g_1)).\]
This leads to 
\[ G_d^{(-\alpha)} \cong \Gamma(\tau_{\rm ncc}) \cong G
  \quad \mbox{ and  } \quad G_d^{(-\alpha), \tau_h^{G_d}} \cong G^h.\]
This shows that $M_d^{(-\alpha)} \cong \cO_h$, as a symmetric space.
Specifically, we have
$C_{\g_d,+} = C_+ \oplus C_+ \subeq \g_{d,1}(h)$, and therefore
\begin{align*}
 C^{(\pm \alpha)}
&  = \{(x,y) \pm (\theta(y), \theta(x)) \: x,y \in C_+ \} \\
&  = \{(x \pm \theta(y), y \pm \theta(x)) \: x,y \in C_+ \} 
  = \Gamma(\pm\theta\res_{C_+ \mp C_-}) \cong C_+  \mp C_-.
\end{align*}
The cone $C^{(\alpha)} \cong C_+ - C_-$ is elliptic and 
$C^{(-\alpha)} \cong C_+ + C_-$ is hyperbolic
(\cite[\S 7.2]{NO23b}). Thus, in the first case, $\cO_h$ inherits
a compactly causal structure
and in the second case a non-compactly causal structure
(cf.~also Proposition~\ref{prop:3.2}).

To see how $M_d^{(\pm \alpha)}$ is embedded in $M \times M$,
we write $\theta^M$ for the involution 
defined by $\exp(\pi z_\fk)$  on~$M$ and 
\[ \infty := \theta^M(o) = \exp(\pi z_\fk).o.\]
Then
\[ M_d^\top := G.(0,\infty) \subeq M_d \]
is the set of ``transversal pairs'', which is open and dense.
It intersects the open subset $\jV \times \jV$ in
\[ \jV_d^\top = \{ (x,y) \in \jV^2 \: x - y \in \jV^\times \} \] 
(see \cite[\S 3.1]{MN25} and \cite[\S 1.5]{BN04}).
As $M_d^{(\alpha)} \subeq M \times M$ is the
$G$-orbit of $(o,o)$ under the action $g.(x,y) := (g.x, \theta^G(g).y)$,
we obtain
\begin{equation}
  \label{eq:mdalpha}
  M_d^{(\alpha)} = (\id_M \times \theta^M)(M_d^\top).
\end{equation}
Writing $\tau_{\rm ncc}^M$ for the involution induced by
$\tau_{\rm ncc}^G$ on $M$, the relation 
\begin{equation}
  \label{eq:md-alpha}
  M_d^{(-\alpha)} = (\id_M \times \tau_{\rm ncc}^M)(M_d^\top),
\end{equation}
follows from the fact that this is the
$G$-orbit of $(o,o)$ under the action
$g.(x,y) := (g.x, \tau_{\rm ncc}^G(g).y)$ and $\tau_{\rm ncc}^M(o) = \infty$.
Both orbits are open and dense. 

\begin{prop} \mlabel{prop:ct-posreg} For the Euler element
  $h_d := (h,-h) \in \fh^{(\pm \alpha)}$,  the following assertions hold:
  \begin{itemize}
  \item[\rm(a)] $W_{M_d^{(-\alpha)}}^+(h_d) = W_M^+(h) \times W_M^+(-h)
    \cong \jV_+ \times (-\jV_+)$, and this domain is causally convex,
    hence globally hyperbolic. 
  \item[\rm(b)] $W_{M_d^{(\alpha)}}^+(h_d) \subeq \jV_+ \times (-\jV_+)$
    is the subset of all pairs $(x,y)$, for which
    $x + y^{-1} \in \jV^\times$.
  \item[\rm(c)] All Euler elements in $\g_d^{(\pm\alpha)} \cong \g$ are conjugate
    to $h_d$. 
  \end{itemize}
\end{prop}

\begin{prf} The Euler element $h_d$ defines on
  $(\g_d^{(\pm \alpha)}, \tau_{h_d}, C^{(\pm \alpha)}, h_d)$
  the structure of a modular causal symmetric Lie algebra,
  and the inclusion into $M_d$ is equivariant for the modular flow.
    First we use  \cite[Lemma~2.7]{MN25} to obtain 
  $W_M^+(\pm h) = \pm\jV_+$, which leads to
  \[  W_{M_d}^+(h_d) = \jV_+ \times (- \jV_+) \subeq M_d^\top.\]
  As $\tau^M_{\rm ncc}(x) = x^{-1} = -\theta^M(x)$ for $x \in \jV$
  leaves $-\jV_+$ invariant,
  we derive from $\jV_+ + \jV_+ \subeq \jV_+ \subeq \jV^\times$ that 
  \[ W_{M_d^{(-\alpha)}}^+(h_d) = \{ (x,y) \in \jV_+ \times (-\jV_+) \:
    x - \tau_{\rm ncc}^M(y) = x - y^{-1} \in \jV^\times\}
    = \jV_+ \times (-\jV_+).\]
  \nin (b) follows with a
  similar argument with $\theta_M$ instead of $\tau_{\rm ncc}^M$. 
  
  \nin (c) follows from \cite[Prop.~3.11]{MN21}
    because $\g$ is of tube type (see also  Theorem~\ref{thm:6.5}). 
\end{prf}

Note that $\theta^M$ is a causal involution on $M$ but 
$\tau_{\rm ncc}^M = \theta^M \tau^M_h$ is not
because $\tau_h^M$ flips the causal structure.

\begin{ex} \mlabel{ex:cayley-ads} It is instructive to visualize the preceding proposition
  for the simplest case where $\jV = \R$ and
  $M = \R_\infty \cong \bS^1$ is the one-point compactification.
  Here $\theta^M(x) = - x^{-1}$ and $\tau_{\rm ncc}^M(x) = x^{-1}$.
  Further $W_{M_d}^+(h_d) = \R_+ \times (-\R_+)$ and
  $M_d^\top = M_d \setminus \Delta_M$. Therefore
  \[  M_d^{(-\alpha)} = M_d \setminus \Gamma(\tau_{\rm ncc}^M)
    = \{(x,y) \in \R_+ \times (-\R_+) \: y \not= x^{-1}\}
    = \R_+ \times (-\R_+)  \]
  and 
  \[  M_d^{(\alpha)} = M_d \setminus \Gamma(\theta^M)
    = \{(x,y) \in \R_+ \times (-\R_+) \: y \not= -x^{-1}\} \]
  is the complement of a hyperbola in a quater plane
  (cf.\ Proposition~\ref{prop:wads-incl} and Example~\ref{ex:2d-ads}).   
\end{ex}

\section{Pierce involutions and satellites of the positive cone}
\mlabel{sec:pierce}

Beyond the flip involutions (cf.~Section~\ref{subsec:flip-ct}), 
the Pierce involutions can easily be understood in terms
of the Pierce decomposition of the euclidean Jordan algebra.
We use \cite{FK94}  as a reference for euclidean Jordan algebras. 
In this section $\jV$ is always a simple euclidean Jordan algebra of 
$\rank(\jV) = r$ and $(c_1, \ldots, c_r) \in \jV$ is a Jordan frame.
 We consider the \textit{Pierce decomposition}
\begin{equation}
  \label{eq:pierce}
  \jV = \bigoplus_{k = 1}^r \R c_k \oplus \bigoplus_{k <\ell} \jV_{k\ell}
  \quad \mbox{with} 
\quad \jV_{k\ell} := \big\{ v \in \jV \:  c_k v = c_\ell v = \shalf v\big\}
\end{equation}
(\cite[\S IV.1]{FK94}).
For $\g = \fco(\jV)$ and
$h \in \fa$, where $\fa \subeq \fp$ is maximal abelian,
we choose a maximal set 
$\{ \gamma_1, \ldots, \gamma_r\} \subeq \Sigma(\g,\fa)$
of strongly orthogonal roots
(cf.\ Appendix~\ref{app:strong-orth}).
We then have
\[ \R c_j = \g_{\gamma_j} \quad\mbox{ and  }\quad
  \jV_{k\ell} = \g_{\frac{1}{2}(\gamma_k + \gamma_\ell)}.\]

\begin{defn} \mlabel{def:pierce-invol} 
For $x_j := c_1 + \ldots + c_j, j = 0,\ldots, r$, the element
$\alpha_j := e^{2\pi i L(x_j)} \in \Str(\jV_\C)$ fixes each $c_j$ in~$\jV$,
and for $v \in \jV_{k\ell}$ we find 
\begin{equation}
  \label{eq:pierce-inv-alpha}
 \alpha_j(v) = 
  \begin{cases}
    v \quad\text{ for }\quad k,\ell  \leq j  \ \text{ or } \ k,\ell > j \\
    -v \quad \text{ for } \quad k \leq j < \ell.    
  \end{cases}
\end{equation}
This shows that $\alpha_j\in \Aut(\jV)$ is an involutive automorphism of~$\jV$,
called a {\it Pierce involution}.
Note that  $\alpha_0 = \alpha_r = \id_{\jV}$.
\end{defn}

Let $\fs \subeq \g$ be the subalgebra isomorphic to $\fsl_2(\R)^{\oplus r}$
generated by the elements $c_j \in \g_1(h)$ and $\theta(c_j) \in \g_{-1}(h)$.
We put 
\[ \fs_j
  := \R c_j  + \R [c_j, \theta(c_j)] + \R \theta(c_j)
  = \g_{\gamma_j} +\g_{-\gamma_j} + \R \gamma_j^\vee
  \cong \fsl_2(\R)\]
(cf.\ Appendix~\ref{app:strong-orth}). 
Then $h \in \fs$ can be written as a sum
\[ h = h_1 + \cdots + h_r = \frac{1}{2}(\gamma_1^\vee + \cdots + \gamma_r^\vee),\]
where $h_j = \shalf\gamma_j^\vee \in \fs_j$ is an Euler element
satisfying $[h_j,c_k] = \delta_{jk} c_j$. In $\Aut(\g)$, we then have
\begin{equation}
  \label{eq:sigma-alphaj}
  \sigma_{\alpha_j}
  = \exp(2\pi i \ad(h_1 + \cdots + h_{r-j})) 
  = \exp(2\pi i \ad(-h_{r-j+1} - \cdots - h_r)).
\end{equation}

\begin{exs} \mlabel{ex:4.1}
  (a) For $\jV = \Herm_r(\K)$ a natural Jordan frame consists
  of the matrix units $c_j = E_{jj}$ and $\jV_{k\ell} = \K (E_{k\ell} + E_{\ell k})$.
  In this case we find with \eqref{eq:alphaq} 
  \[ \alpha_j(x) = I_{r-j,j} x I_{r-j,j}.\] 

  \nin (b) For $\jV = \R^{1,d-1}$ we have $r = 2$ and  the two
  lightlike vectors 
  \[ c_1 = \be_0 + \be_1, \quad  c_2 = \be_0 - \be_1 \]
  form a Jordan frame. The only Pierce space is 
  $\jV_{12} = \Spann \{\be_2, \ldots, \be_{d-1}\}$
  and
  \[ \alpha_1 = \diag(1,1,-1,\ldots, -1) = I_{2,d-2}.\] 
\end{exs}

\begin{prop} \mlabel{prop:satj}
Let $\alpha_j \in \Aut(\jV)$, $j =0,\ldots, r$,
  be a Pierce involution and consider the Euler elements
\begin{equation}
  \label{eq:hj}
  h^j := h_1 + \cdots + h_{r-j} - h_{r-j+1} - \cdots - h_r
  \quad \mbox{ and }\quad
k^j := e^{-\frac{\pi}{2} \ad z_\fk} h^j.  
\end{equation} 
  Then the following assertions hold:
  \begin{itemize}
  \item[\rm(a)] $\theta_{-\alpha_j}$ is the non-compactly causal
    involution $\tau_{\rm ncc} = \theta \tau_{h^j}$
    on $\g$, specified by the Euler element $h^j$
    in the sense of {\rm Definition~\ref{def:reductive-ncc}}. 
  \item[\rm(b)] $\theta_{-\alpha_j} = \tau_{k^j}$ for the Euler element $k^j$.
    In particular,  the symmetric Lie algebra $(\g,\theta_{-\alpha_j})$ is of Cayley type and $k^j \in \fq^{(-\alpha_j)}$.  
  \end{itemize}
\end{prop}

\begin{prf}
  (a)   The subalgebra $\fs$ also contains the element $z_\fk \in \fz(\fk)
\cap C_\g$ for which
$e^{\pi \ad z_\fk} = \theta$ is the Cartan involution of~$\g$
(cf.~\eqref{eq:thetazk}). 
Therefore $\theta_{\alpha_j}$ is given by conjugation (in $G_\C$) with the
element 
\begin{equation}
  \label{eq:pierce-theta-alpha}
  \exp(\pi z_\fk) \exp(2\pi i(h_1 + \cdots + h_{r-j})),
\end{equation}
and $\theta_{-\alpha_j}$ by conjugation with 
\begin{align*}
  & \exp(\pi z_\fk)  \exp(-\pi i h) \exp(2\pi i(h_1 + \cdots + h_{r-j}))\\
&  = \exp(\pi z_\fk)  \exp(\pi i(h_1 + \cdots + h_{r-j}-  h_{r-j+1} - \cdots -  h_r))  = \exp(\pi z_\fk)\exp(\pi i h^j).
\end{align*}
Hence $\theta_{-\alpha_j}$ is
the non-compactly causal involution $\tau_{\rm ncc} = \theta \tau_{h^j}$
on $\g$, specified by the
Euler element $h^j$ of $\g$ (Definition~\ref{def:reductive-ncc}).

\nin (b) For the Euler element
\begin{equation}
  \label{eq:kj}
 k^j = e^{-\frac{\pi}{2} \ad z_\fk} h^j
\ {\buildrel \eqref{eq:90deg-rot}\over =}\   k_1 + \ldots + k_{r-j} - k_{r-j+1} - \cdots - k_r
\in \fs \cap \fp,
\end{equation}
we then derive from $\tau_{h^j}(z_\fk) = - z_\fk$ and (a) that 
\begin{equation}
  \label{eq:tau-rel}
\theta_{-\alpha_j} = \theta \tau_{h^j}
 = e^{-\pi\ad z_\fk} \tau_{h^j}
 = e^{-\frac{\pi}{2} \ad z_\fk} \tau_{h^j} e^{\frac{\pi}{2} \ad z_\fk}
  = \tau_{k^j}
\end{equation}
(cf.~Example~\ref{ex:sl2a}).
Finally, $k^j \in \fq^{(-\alpha_j)} = \g^{(-\alpha_j), -\tau_h}$ follows from
  \[  \tau_h(k^j)
    = e^{\frac{\pi}{2} \ad z_\fk} \tau_h(h^j)
    = e^{\frac{\pi}{2} \ad z_\fk} h^j 
    = e^{\pi\ad z_\fk} k^j  
    = \theta(k^j) = - k^j.\qedhere\]     
\end{prf}

The set $\jV^\times$ of invertible elements of $\jV$ 
has $r+1$ connected components that can be described 
as follows. We fix a spectral decomposition 
$x = \sum_{j = 1}^r x_j \tilde c_j$ of an element $x \in \jV$, 
where $(\tilde c_1,\ldots, \tilde c_r)$ 
is a Jordan frame (\cite[Thm.~III.1.1]{FK94}). 
This means that, under the automorphism group $\Aut(\jV)$, 
the element $x$ is conjugate to $\sum_{j = 1}^r x_j c_j$.
\begin{footnote}
{For $\jV = \Herm_r(\K)$, this corresponds to the conjugation 
of a hermitian matrix $x$ by $\U_r(\K)$ to a diagonal matrix, 
and for Minkowski space, it corresponds to conjugation 
of an element $x \in \R^{1,d-1}$ under $\OO_{d-1}(\R)$ to one 
with $x_2 = \cdots = x_{d-1} = 0$.}
\end{footnote}
The element $x$ is invertible if and only if all coefficients
$x_j \in \R$ are non-zero,
and if this is the case, we say that $x$ is of {\it signature} $(p,q)$ if
$p$ is the number of positive coefficients
$x_j$ and $q = r - p$ the number of negative ones.
We write
\begin{equation}
  \label{eq:vtimespq}
  \jV^\times_{(p,q)} \subeq \jV^\times, \quad p + q = r,
\end{equation}
for the open subset of elements of signature  $(p,q)$. 
These are orbits of the group $G^h_e \cong \Str(\jV)_e$.
For $0 < p < r$, they are called 
{\it satellites} of the open convex cone $\jV_+ = \jV^\times_{(r,0)}$. Note that
$-\jV_+ = \jV^\times_{(0,r)}$ and 
\begin{equation}
  \label{eq:cj}
  c^j = c_1 + \cdots + c_{r-j} - c_{r-j+1} - \cdots - c_r
  \in \jV^\times_{(r-j,j)}.
\end{equation}

The following proposition relates
the domains $\jV^\times_{(p,q)}$
to the open orbits $M^{(-\alpha)}$ for Pierce involutions.

\begin{prop} \mlabel{prop:pierce-concrete}
  In the subalgebra $\fs \cong \fsl_2(\R)^{\oplus r}$
  from \eqref{app:strong-orth}, 
  we  consider the elliptic element
    \[ z_\fk^j := z_{\fk,1} + \cdots + z_{\fk,r-j} - z_{\fk,r-j+1}
      - \cdots - z_{\fk,r}\]
    and put $d_j := \exp\big(\frac{\pi}{2} z^j_\fk\big)$. 
    Then the following assertions hold:
    \begin{itemize}
    \item[\rm(a)] $\theta_{-\alpha_j} = \Ad(d_j)^{-1} \tau_h \Ad(d_j)$. 
    \item[\rm(b)] $d_j.M^{(-\alpha_j)} = \jV^\times_{(r-j,j)}$
      and $d_j(0) = c^j$.
    \item[\rm(c)] The involution $\tau_h$ on $\g^{(-\alpha_j)}
      = \Ad(d_j)^{-1} \g_0(h)$ 
      corresponds to the non-compactly causal involution on $\g_0(h)$
      defined by the Euler element $h^j \in \g_0(h)$.
    \item[\rm(d)] The positivity region of $h^j$ in $M$ is contained 
      in $\jV^\times_{(r-j,j)}$ and given by the convex cone 
\begin{equation}
  \label{eq:w-def}
W_M^+(h^j) =  W_{\jV}^+(h^j)  
  = (\oline{\jV_+} \cap \jV_1(h^j))^\circ 
   -  (\oline{\jV_+} \cap \jV_{-1}(h^j))^\circ +  \jV_0(h^j)
\end{equation}
    \end{itemize}
\end{prop}

For $r = 2$ and $j = 1$ (Lorentz boost on Minkowski space),
  the wedge domain $W_\jV^+(h^1)$ in Minkowski space $\R^{1,d-1}$,
  is known as the {\it Rindler wedge}.
  For $j = 0,r$, we obtain the open cones
  $\jV_+$ and $- \jV_+$ for any Jordan algebra.

\begin{prf} (a), (b): The Cayley transform 
$c(x) = \exp\big(\frac{\pi}{2} z_\fk\big).x
    = \frac{e+x}{e-x}$ 
from \eqref{eq:cayley-real-sl2} satisfies $c(0) = e$ and $c^{-1}(0) = - e$. 
  Using the subalgebra $\fs \cong \fsl_2(\R)^{\oplus r}$,
  we apply $c$ in the first $r-j$ components and
  its inverse in the last $j$ components.
  Concretely, we obtain 
\begin{equation}
  \label{eq:mixcayley}
d_j.0 =  \exp\Big( \frac{\pi}{2} z^j_\fk\Big).0 = c^j. 
\end{equation}
  Therefore the $\Str(\jV)_e$-orbit $\jV^\times_{(r-j,j)}$ of $c^j$ corresponds to the
  orbit of $0$ under the connected subgroup with the Lie algebra 
  \[ \Ad(d_j)^{-1} \g_0(h) = e^{-\frac{\pi}{2} \ad z_\fk^j} \g_0(h)
    \ {\buildrel\eqref{eq:htok}\over =}\ \g_0(k^j) \
    {\buildrel\ref{prop:satj}(b)\over =}\ \g^{(-\alpha_j)}, \]
  where we use \eqref{eq:90deg-rot} in Appendix~\ref{app:a.3} to see that 
  \begin{equation}
    \label{eq:htok}
k^j =  e^{-\frac{\pi}{2} \ad z^j_\fk} h
 =  k_1 +  \cdots k_{r-j} -  k_{r - j+1} -  \cdots -  k_r.
  \end{equation}
  With Proposition~\ref{prop:satj}, this shows that
  $\g^{(-\alpha_j)} = \g_0(k^j) = \Ad(d_j)^{-1}\g_0(h)$,
  hence~(a) and (b).

\nin (c)  On $\g^{(-\alpha_j)}$ the involution is given by $\tau_h$, so that
  \begin{equation}
    \label{eq:star}
\Ad(d_j)h =  e^{\frac{\pi}{2} \ad z_\fk^j} h = -k_1 - \cdots - k_{r-j}
 + k_{r-j+1} +  \cdots + k_r = -k^j
  \end{equation}
  and $\tau_{-k^j} = \tau_{k^j} =  \theta \tau_{h^j}$ 
  imply that  the corresponding involution on $\g_0(h)$ is~$\tau_{k^j}$ 
  (cf.~\eqref{eq:tau-rel}).
  In view of
  Definition~\ref{def:reductive-ncc}, this proves~(c).

  \nin (d) The eigenspaces of $h^j$ on $\jV$ are given by 
  \begin{equation}
    \label{eq:h-v-eigen}
    \jV_1(h^j) = \bigoplus_{\ell = 1}^{r-j} \R c_\ell \oplus
    \bigoplus_{\ell < m \leq r-j}
  \jV_{\ell m},  \quad 
\jV_0(h^j) = \bigoplus_{\ell \leq r-j < m} \jV_{\ell m} , \quad 
\jV_{-1}(h^j) = \bigoplus_{\ell > r-j} \R c_\ell \oplus
\bigoplus_{r-j < \ell < m} \jV_{\ell m},
  \end{equation}
where $\jV_{\pm 1}(h^j)$ are Jordan subalgebras of $E$.  
By \cite[eq.~(B.4)]{NOO21}, the positivity region of $h^j$ in~$\jV$ is 
\begin{equation}
  \label{eq:w-def2}
 W_{\jV}^+(h^j) = \{ x \in \jV \: [h^j,x] \in \jV_+ \}
  = (\oline{\jV_+} \cap \jV_1(h^j))^\circ 
   -  (\oline{\jV_+} \cap \jV_{-1}(h^j))^\circ +  \jV_0(h^j)
\end{equation}
Further \cite[Cor.~B.5]{NOO21} shows that
$W_{\jV}^+(h^j)  \subeq M_j := \jV^\times_{(r-j,j)}$, so that
$W_{\jV}^+(h^j) = W_{M_j}^+(h^j)$. 
\end{prf}

The cones $W_{\jV}^+(h^j)$ are in particular causally convex subsets of $\jV$. 

\begin{lem} \mlabel{lem:inv1}
  $x \in \jV^\times$ if and only if $h \in [x,\g_{-1}(h)]$.
\end{lem}

\begin{prf} In view of the Spectral Theorem for euclidean Jordan Algebras
  (\cite{FK94}), we may assume that
  $x = \sum_{j = 1}^r x_j c_j$ with $x_j \in \R$.
  With respect to the strongly orthogonal
  roots $\gamma_1, \ldots, \gamma_r$ with $\g_{\gamma_j} = \R c_j$, we then have
  \[ [x,\g_{-1}(h)] \cap \fa
    = \Big[x, \sum_{j = 1}^r \g_{-\gamma_j}\Big]
    = \sum_{x_j \not=0} \R \gamma_j^\vee.\]
  Hence $h = \shalf \sum_j \gamma_j^\vee$ is contained in
  $[x,\g_{-1}(h)]$ 
  if and only if all $x_j$ are non-zero, i.e., if $x \in \jV^\times$.   
\end{prf}

\begin{lem} \mlabel{lem:fixedpoints}
$\jV_0(h^j) \cap \jV^\times \not=\eset$ if and only if $r= 2j$. 
\end{lem}

\begin{prf} We have
  \[ \jV_0(h^j)
    = \sum_{k \leq r-j < \ell} \jV_{k\ell} 
    = \sum_{k \leq r-j < \ell} \g_{\shalf(\gamma_k + \gamma_\ell)}.\]
An element $x = \sum_{k \leq r-j < \ell} x_{k\ell}$ with
  $x_{k\ell} \in \g_{\shalf(\gamma_k + \gamma_\ell)}$
  is invertible in $V_0(h^j)$ if and only if 
  Lemma~\ref{lem:inv1} implies that $h = [x,y]$ for some
  $y = \sum_{k,\ell} y_{k\ell} \in \g_{-1}(h)$. We then must have
$h = \sum_{k,\ell} [x_{k\ell}, y_{k\ell}]$ 
  because all other brackets lie in different $\fa$-eigenspaces.
  For any $\beta \in \Sigma(\g,\fa)$, we have
  $p_\fa([\g_\beta, \g_{-\beta}]) = \R \beta^\vee$
($\beta^\vee$ is the corresponding coroot), 
so that 
\[ h \in \sum_{k \leq r-j < \ell} \R \Big(\frac{\gamma_k + \gamma_\ell}{2}\Big)^\vee
  \subeq (h^j)^\bot.\]
This implies that $h \bot h^j$ with respect to the Cartan--Killing form.
As
\[ h = \frac{1}{2}(\gamma_1^\vee + \cdots + \gamma_r^\vee) \quad \mbox{ and } \quad
  h^j = \frac{1}{2}(\gamma_1^\vee + \cdots + \gamma_{r-j}^\vee -
  \gamma_{r-j+1}^\vee - \cdots - \gamma_r^\vee),\]
and the $\gamma_j^\vee$ are pairwise orthogonal of the same length,
it follows that $r = 2j$.

Assume, conversely, that $r = 2j$. Then $r$ is even, so that either
$\jV \cong  \R^{1,d}$ (for $r = 2$), or $\jV \cong \Herm_r(\K)$ for $r > 2$
and $\K \in \{\R,\C,\H\}$. 
For $r = 2$, the space $V_0(h^1)$ ($h^1$ a Lorentz boost) 
contains non-zero spacelike points (because $d > 1$ by simplicity of $\jV$),
and these are invertible elements of $\jV$.

For $r > 2$ even, $\jV \cong  \Herm_r(\K)$, 
and for $j = r/2$ we have 
\begin{equation}
  \label{eq:even-invert}
  \jV_0(h^j) = \Big\{ \pmat{ 0 & b \\ b^*  & 0}\: \ b \in M_j(\K)\Big\}.
\end{equation}
The invertible elements in this set are those for which $b$ is invertible. 
\end{prf}

The non-compactly causal symmetric spaces
\begin{equation}
  \label{eq:def-mj}
  \jV^\times_{(r-j,j)}\ {\buildrel\ref{prop:pierce-concrete}\over \cong}\  M^{(-\alpha_j)}
\end{equation}
carry natural flows generated by     the Euler elements $h^\ell$,
$\ell = 0,\ldots, r$, of $\g$ contained in $\g_0(h) \cong \g^{(-\alpha_j)}$.
The following  proposition clarifies the existence of positivity regions and
 fixed points: 

 \begin{prop} \mlabel{prop:fixed-points}
   {\rm(Modular flows on $\jV^\times_{(r-j,j)}$)} 
For the flow generated by the Euler element $h^\ell \in \g_0(h)$ on 
$\jV^\times_{(r-j,j)}$, the following assertions hold:
\begin{itemize}
\item[\rm(a)] The positivity region of $h^\ell$ 
is non-empty if and only if $\ell = j$. 
\item[\rm(b)] The symmetric space $\jV^\times_{(r-j,j)}$ is the 
 non-compactly causal symmetric space corresponding to the pair
 $(\g_0(h), h^j)$ {\rm(cf.\ Definition~\ref{def:reductive-ncc})}. 
\item[\rm(c)] It has fixed points if and only if $r = 2j$ and $\ell = j$.
\item[\rm(d)] For $r = 2j$, there exists  an Euler element
  $h' \in \fh^{(-\alpha_j)} \cap [\g^{(-\alpha_j)}, \g^{(-\alpha_j)}]$  for which 
  \[ (\g^{(-\alpha_j)}, \theta_{-\alpha_j}, C^{(-\alpha_j)}, h') \]
  is a modular causal symmetric Lie algebra, i.e., 
 $\theta_{-\alpha_j}(h') = h'$ and $\tau_{h'}(C^{(-\alpha_j)}) = - C^{(-\alpha_j)}$.
\end{itemize}
\end{prop}

Assertion (a) implies that
all Euler elements $h^\ell, \ell \not=j$, have empty
positivity regions in $\jV^\times_{(r-j,j)}$
(cf.\ Definition~\ref{def:posreg}).

\begin{prf} (a)  As the positivity region $W^\ell$ of $h^\ell$ 
  in $\jV$ is contained in  $\jV^\times_{(r-\ell,\ell)}$ by
  Proposition~\ref{prop:pierce-concrete}(d),
it intersects $\jV^\times_{(r-j,j)}$ non-trivially if and only if $j = \ell$.

\nin (b) follows from the fact that
the corresponding involution on $\g_0(h)$ is 
$\tau_{h^j} \theta$ (Proposition~\ref{prop:pierce-concrete}(c)),
and that $\fz(\g_0(h)) \subeq \g_0(h)^{-\tau_{h^j\theta}}$.

\nin (c) By Lemma~\ref{lem:fixedpoints}, 
the existence of fixed points in $\jV^\times$ 
is equivalent to $r = 2\ell$. It remains to verify that they are have
signature $(r-\ell,\ell) = (\ell,\ell)$. 
We thus assume that $r$ is even.
For $r = 2$, the flow of $h^1$ (the Lorentz boost) 
has only spacelike fixed points, and these are of signature $(1,1)$.

For $r > 2$ even, we have $\jV \cong  \Herm_r(\K)$ for
$\K \in \{\R,\C,\H\}$ and \eqref{eq:even-invert} implies that
the invertible elements in $\jV_0(h^\ell)$ are of the form
\[ \pmat{ 0 & b \\ b^*  & 0}, \quad b \in \GL_\ell(\K).\] 
It contains the elements
\begin{equation}
  \label{eq:b-diag}
  \pmat{ 0 & b \\ b & 0}, \quad b = \diag(b_1, \ldots, b_\ell), \ b_k\in \R^\times,
\end{equation}
and it is easy to see that all these elements are of signature $(\ell,\ell)$.
For $\K = \C,\H$, the open subset $\GL_r(\K)$ is connected, and for
$\K = \R$, it has $2$ connected components. They contain 
matrices as in \eqref{eq:b-diag}, so that they are completely contained
in $\jV^\times_{(\ell,\ell)}$.

\nin (d) Assume that $r = 2j$.
  Then the proof of Proposition~\ref{prop:4.9}
below implies that the Euler element $h^j \in \g_0(h)$ is symmetric,
hence in particular contained in the commutator algebra
$[\g_0(h), \g_0(h)]$, which is semisimple. This implies the
existence of a $\theta$-invariant
$\fsl_2$-subalgebra $\fs_0\subeq \g_0(h)$ containing
$h^j$ (\cite[Thm.~5.4]{MNO23}). The subalgebra $\fs_0$ is also invariant
under $\theta_{-\alpha_j} = \theta \tau_{h^j}$. The fixed point space for the
latter involution on $\fs_0$ is $\R h'$ for some Euler element
$h' \in \Inn(\fs_0)h \subeq \Inn(\g)h$.
Then $h' \in \g^{(-\alpha_j)}$ is an Euler element of $\g$.
By  construction, $h' \in \g_0(h)$, so that
$h'\in \fh^{(-\alpha_j)}$. Now (d) follows from
Proposition~\ref{prop:aut-modul}(b). 
\end{prf}

\begin{prop} \mlabel{prop:5.9} The positivity regions
  $W_{M^{(-\alpha_j)}}(h^j)$ are real tube domains on which the
  subgroup $(G^{h^j})_e$ acts transitively.   
\end{prop}

Note that the subgroup $(G^{h^j})_e$ is in general not contained in
$G^{(-\alpha_j)}$, hence does not act by automorphisms of the symmetric space 
$M^{(-\alpha_j)}$ on $W_{M^{(-\alpha_j)}}(h^j)$ (Example~\ref{ex:5.9}). 

\begin{prf} Since the group $G = \Co(\jV)_e$ acts transitively
  on the set $\cE(\g)$ of Euler elements in $\g$
  (\cite[Prop.~3.11]{MN21}), there exists a
  $g \in G$ with $\Ad(g) h = h^j$. Then
  \[ \Ad(g) \jV_+ = \Ad(g) W_M^+(h)  = W_M^+(h^j)
    \ {\buildrel\ref{prop:pierce-concrete}(d) \over=}\ W_{M^{(-\alpha_j)}}^+(h^j),\]
  and since $(G^h)_e$ acts transitively on $\jV_+ = W_M^+(h)$, the group
  $(G^{h_j})_e\cong \Str(\jV)_e$ acts transitively on $W_{M^{(-\alpha_j)}}^+(h^j)$.   
\end{prf}

\begin{ex} \mlabel{ex:5.9} For $\jV = \R^{1,d-1}$ and the Rindler
  wedge $W_R = \{ (x_0, \bx) \: x_1 > |x_0|\} = W_\jV^+(h^1),$ 
  the affine subgroup
  \[ \Spann \{\be_2, \ldots, \be_{d-1}\} \rtimes \R_+^\times\SO_{1,1}(\R)_e
    \cong \jV_0(h^1) \rtimes \exp(\R h + \R h^1)
  \subeq G^{h^1} \]
acts transitively on~$W_R$. It is not contained in $G^{(-\alpha_1)} 
\cong \Str(\jV)_e = \R_+ \SO_{1,d-1}(\R)_e$.
\end{ex}

\section{Modular flows on causal Makarevi\v c spaces}
\mlabel{sec:7}

The Cayley type spaces discussed in Section~\ref{subsec:flip-ct}
are always modular, but the modularity property is 
not always satisfied for the non-compactly causal spaces
$M^{(-\alpha_j)}$ of Pierce type, for which it is equivalent to $r = 2j$,
and this condition also characterizes the modularity
of the dual compactly causal space~$M^{(\alpha_j)}$.
The main result of this section describes those
compactly causal symmetric Lie algebras $(\g^{(\alpha)}, \tau_h, C^{(\alpha)})$
which are modular in the sense that $\fh^{(\alpha)} = \g^{(\alpha)} \cap \g_0(h)$
contains an Euler element $h'$ of $\g^{(\alpha)}$ satisfying 
$-\tau_{h'}(C^{(\alpha)}) = C^{(\alpha)}$. 
If, in addition, $h' \in \fh^{(\alpha)}$ is an Euler element of $\g$,
then the latter invariance condition is automatically
satisfied by Proposition~\ref{prop:aut-modul}.
In this context it is remarkable that 
Theorem~\ref{thm:6.5} asserts that the existence of   
non-central Euler elements in $\g^{(\alpha)}$ already implies
the existence of an Euler element of $\g$, contained in $\fh^{(\alpha)}$.


\begin{thm} \mlabel{thm:6.5} {\rm(Classification of the
    modular spaces $M^{(\alpha)}$)} 
  Whenever the Lie algebra $\g^{(\alpha)}$ contains a non-central Euler element,
  then $\fh^{(\alpha)} \cap [\g^{(\alpha)}, \g^{(\alpha)}]$
  contains an Euler element of $\g$.
  All Euler elements in $[\g^{(\alpha)}, \g^{(\alpha)}]$
    are conjugate under  $\Inn(\g^{(\alpha)})$, and for any such
    Euler element $h'$, the corresponding positivity region
    $W_{M^{(\alpha)}}^+(h')$ is non-empty. 
This happens for Cayley type and in the following cases: 
\begin{itemize}
\item[\rm(P):] For $r$ even and $j = r/2$,  
\item[\rm(S):] for {\rm(S1)} and {\rm(S2)} if $r$ is even, and
  for {\rm(S4)} and 
  $\g^{(\alpha)} \cong \so_{2,k}(\R) \oplus \so_{d-k}(\R)$
  if $1 \leq k < d-2$. 
\item[\rm(NS):] always for {\rm(NS1)}, {\rm(NS2)}, {\rm(NS3)}. 
\end{itemize}
\end{thm}

\begin{prf}   In view of Proposition~\ref{prop:aut-modul}(b),
  all Euler elements of $\g$, contained
  in $[\g^{(\alpha)}, \g^{(\alpha)}]$
  are $\Inn(\g^{(\alpha)})$-conjugate to some
  $h'$ for which 
  $(\g^{(\alpha)}, \tau_h, C^{(\alpha)},h')$
  is a modular compactly causal symmetric Lie algebra,
  and by Proposition~\ref{prop:aut-modul}(c),
  this implies that $W_{M^{(\alpha)}}^+(h')$ is non-empty. 
It therefore remains to check when
$\g^{(\alpha)}$ contains non-central Euler elements $h'$
(then its component in the commutator algebra is also Euler), 
and that all Euler elements of $\g$, contained in 
$[\g^{(\alpha)}, \g^{(\alpha)}]$ lie in the same $\Inn(\fg^{(\alpha)})$-orbit.\\

\nin {\bf Cayley type:}  For flip involutions we have
$\g_d^{(\alpha)} \cong \g$ for
$\g_d = \g^{\oplus 2}$ (see Section~\ref{subsec:flip-ct}).
If $\g$ contains an Euler element 
$h' = (\theta(h''),h'')$, then $h''$ is an Euler element of $\g$
and $h'$ an Euler element of~$\g_d$.
Then $\g$ is simple hermitian of tube type, so that
all  Euler elements in $\g$ are conjugate (\cite[Prop.~3.11]{MN21}).


Having dealt with the Cayley type involutions, we now
address the other three types (P), (S), and (NS), according to the type
of the corresponding simple euclidean Jordan algebra $\jV$. \\

\nin {\bf $\jV = \Sym_r(\R)$ and $\g = \sp_{2r}(\R)$:} 
\begin{itemize}
\item[(P)] In this case
  $\g^{(\alpha)} \cong \fu_{j,r-j}(\C) \cong \R \oplus \su_{j,r-j}(\C)$,
  and this Lie algebra contains a non-central Euler element if and only
  if $r = 2j$. Assuming that $r = 2s$ is even, it corresponds to the
  embedding
  \[ \g^{(\alpha)} \cong \fu_{s,s}(\C) \into \sp_{2r}(\R). \]
  Realizing $\fu_{s,s}(\C) \cong \fu(\Omega, \C^{2s})$
  as in \eqref{eq:urrcb} (Examples~\ref{ex:2.6}(a)),
\[ h =   \shalf \diag(\1_s, -\1_s) \in \fh^{(\alpha)}
  \cap [\g^{(\alpha)}, \g^{(\alpha)}] \] 
  is an Euler element, and since it has only 
  two eigenvalues, it is also Euler in~$\gl_{4s}(\R) =\gl_{2r}(\R)$,
  hence in particular in $\sp_{2r}(\R)$.
  All other Euler elements of the commutator algebra
  $\su_{s,s}(\C)$, which is hermitian of tube type,
  are conjugate by \cite[Prop.~3.11]{MN21}.
\item[(NS1)] Here $\g^{(\alpha)} \cong \sp_{2s}(\R)^{\oplus 2}  \subeq
  \g \cong \sp_{4s}(\R)$. In this case $\g^{(\alpha)}$ contains
  $3$ conjugacy classes of Euler elements, represented by
  $(h_0, h_0), (h_0, 0)$ and $(0, h_0)$, where
$h_0 =   \shalf \diag(\1_s, -\1_s) \in \sp_{2s}(\R).$ 
  The Euler elements $(h_0,0)$ and $(0, h_0)$ have $3$-eigenvalues
  on $\R^{4s}$, hence are not Euler in $\sp_{4s}(\R)$
  because this Lie algebra contains only one class of Euler elements
  (\cite[Prop.~3.11]{MN21}). 
  But $h := (h_0, h_0) \in \fh^{(\alpha)}$ is also Euler in $\sp_{4s}(\R)$,
  corresponding to $\shalf \diag(\1_{2s}, -\1_{2s})$.
\end{itemize}
\nin{\bf $\jV = \Herm_r(\C)$ and $\g = \su_{r,r}(\C)$:} 
\begin{itemize}
\item[(P)]   In this case
  $\g^{(\alpha)} \cong \R(i\1, -i\1) \oplus \su_{j,r-j}(\C)^{\oplus 2}$,
  and this Lie algebra contains a non-central Euler element if and only
  if $r = 2j$, which corresponds to the modularity
  of the compactly causal symmetric Lie algebra
  $(\g^{(\alpha)}, \theta_\alpha, C^{(\alpha)})$
  (Proposition~\ref{prop:aut-modul}). 
  Assuming that $r = 2s$ is even and $j = s$, it corresponds to the
  embedding
  \[ \R \oplus \su_{s,s}(\C)^{\oplus 2} \into \su_{r,r}(\C). \]
  Realizing $\fu_{s,s}(\C) \subeq  \fu(\Omega, \C^{2s})$
  as in \eqref{eq:urrcb} in Examples~\ref{ex:2.6}(a),
  we see that
$h_0 =   \shalf \diag(\1_s, -\1_s)$ 
is an Euler element in $\su_{s,s}(\C) = [\g^{(\alpha)}, \g^{(\alpha)}]$,
and since it has only 
  two eigenvalues, $h := h_0 \oplus h_0$
  is also Euler in $\su_{r,r}(\C)$ (there is only one conjugacy class
  by \cite[Prop.~3.11]{MN21}). The other Euler elements in
  the commutator algebra of $\g^{(\alpha)}$ are conjugate to
  $(h_0, 0)$ or $(0, h_0)$, which are not Euler in $\su_{r,r}(\C)$.
\item[(S1)] Here $\g^{(\alpha)} \cong \so^*(2r) = \so^*(4s)$, which contains
  an Euler element if and only if $r = 2s$ is even, i.e.,
  $\g^{(\alpha)}$ is of tube type.
  Then we have a natural embedding
  \[ \so^*(4s) \cong \fu(\Omega, \H^{2s})
    \into  \su_{2s,2s}(\C) \cong \su_{r,r}(\C) = [\fu(\Omega, \C^{2s}), \fu(\Omega, \C^{2s})], \]
  (Example~\ref{ex:2.6}(a)).
  This Lie algebra contains a single orbit of Euler elements,
  represented by
$h = \shalf\diag(\1_s,-\1_s) \in \gl_{2s}(\H),$ 
  and this also is an Euler element of $\su_{r,r}(\C)$.
\item[(NS2)] Here $\g^{(\alpha)} \cong \sp_{4r}(\R) \subeq
  \g \cong \su_{2r,2r}(\C)$. In this case $\g^{(\alpha)}$ contains a single
  conjugacy class of Euler elements (\cite[Prop.~3.11]{MN21}), represented by
$h =   \shalf \diag(\1_{2r}, -\1_{2r}),$ 
  and this is also an Euler element of~$\su_{2r,2r}(\C)$.
\end{itemize}
\nin {\bf $\jV = \Herm_r(\H)$ and $\g = \so^*(4r)$:} 
\begin{itemize}
\item[(P)] In this case
  $\g^{(\alpha)} \cong \fu_{2j,2r-2j}(\C)$,
  and this Lie algebra has a non-central Euler element if and only
  if $r = 2j$. Assuming that $r = 2s$ is even and $j = s$,
  it corresponds to the   embedding
  \[ \fu_{r,r}(\C) \cong \fu(\Omega, \C^{2r}) \into \fu(\Omega, \H^{2r}) \cong
    \so^*(4r) \]
(Examples~\ref{ex:2.6}(a)). 
So $h_0 =   \shalf \diag(\1_r, -\1_r)$ 
is an Euler element in the commutator algebra
$\su_{r,r}(\C)$, which is also Euler
in $\fu(\Omega, \H^{2r})$.
Uniqueness of the conjugacy class in
$\su_{r,r}(\C)$ follows from \cite[Prop.~3.11]{MN21}.
\item[(S2)] Here $\g^{(\alpha)} \cong \so^*(2r)^{\oplus 2}
  \subeq \so^*(4r)$. In this case $\g^{(\alpha)}$ contains 
  an Euler element if and only if $r = 2s$ is even
  (\cite[Prop.~3.11]{MN21}). 
  There exist three orbits of Euler elements, represented by
    $(h_0, h_0), (h_0, 0)$ and $(0, h_0)$, where
    $h_0 =   \shalf \diag(\1_s, -\1_s) \in \fu(\Omega, \H^{2s})
    \cong \so^*(4s).$ 
    The elements $(h_0,0)$ and $(0,h_0)$ are not Euler in $\so^*(4r)$,
    but $h' := (h_0, h_0) \in \fh^{(\alpha)}$ is also Euler in  $\so^*(4r)$. 
\end{itemize}
\nin {\bf $\jV = \Herm_3(\bO)$ and $\g = \fe_{7(-25)}$:} 
\begin{itemize}
\item[(P)] Here it suffices to consider the case $j = 1$, where 
$\g^{(\alpha)} \cong \R \oplus \fe_{6(-14)}$. Its commutator algebra 
is hermitian, but not of tube type, so that it contains no non-central 
Euler element. Here $r = 3$ is odd. 
\item[(S3)] In this case $\g^{(\alpha)} \cong \su_{2,6}(\C)$, which is 
  hermitian, not of tube type, so that it contains
  no Euler element (\cite[Prop.~3.11]{MN21}). 
\end{itemize}
\nin {\bf $\jV = \R^{1,d-1}$ and $\g = \so_{2,d}(\R)$:}
From \cite[\S 2.3]{MNO25} it follows that, for $d > 2$,
and a subalgebra of type $\so_{2,2}(\R) \cong \fsl_2(\R)^{\oplus 2}$,
the only conjugacy class of Euler elements in $\so_{2,2}(\R)$
that consists of Euler elements of $\so_{2,d}(\R)$
is the one consisting of Euler elements whose centralizer
contains no non-zero ideal. We shall use this several times 
below for the uniqueness of the conjugacy class in~$\g^{(\alpha)}$. 
\begin{itemize}
\item[(P)] Here $j= 1 = \frac{r}{2}$
  and $\g^{(\alpha)} \cong \R \oplus \so_{2,d-2}(\R)$.
  Its commutator algebras $\so_{2,d-2}(\R)$ contains 
  Euler elements of $\g$, and, by \cite[\S 2.3]{MNO25},
  they are all conjugate, even for $d = 4$.
\item[(S4)] Here $\g^{(\alpha)} \cong \so_{2,k}(\R) \oplus \so_{d-k}(\R)$ 
  with $0 \leq k < d-2$ (see \eqref{eq:mink-galpha} below),
  and this Lie algebra contains non-central Euler 
  elements which are also Euler in $\so_{2,d}(\R)$ if and only
  if $k \geq 1$.
  Even for $k = 2$ the Lie algebra $\so_{2,2}(\R)$ is not simple, but
  there is only one such class by \cite[\S 2.3]{MNO25}.
\item[(NS3)] Here $\g^{(\alpha)} \cong \so_{2,d-1}(\R)$ 
  contains Euler elements which are also Euler in $\so_{2,d}(\R)$. 
  Even for $d = 3$, there is only one such class by \cite[\S 2.3]{MNO25}.
\qedhere\end{itemize}
\end{prf}

\hspace{-22mm}
\begin{tabular}{||l|l|l|l|l|l|l|l||}\hline
  Type & $\fco(\jV)$ & $\g^{(\alpha)}$  & $\g^{(-\alpha)}$
  & $\fh^{(\alpha)} = \fh^{(-\alpha)}$ & $M^{(\alpha)}$ & $M^{(-\alpha)}$   \\ 
\hline 
\hline 
$\Sym_r(\R)$\phantom{\Big (} & $\sp_{2r}(\R)$ && &&& \\
\hline   
 (C) &$\sp_{2r}(\R)^{\oplus 2}$ & $\sp_{2r}(\R)$ & $\sp_{2r}(\R)$  &  $\gl_r(\R)$ && \\
  (P) & $\sp_{4s}(\R)$ & $\fu_{s,s}(\C)$ & $\gl_{2s}(\R)$  &  $\fo_{s,s}(\R)$ &
$\U_{s,s}(\C)/\OO_{s,s}(\R)$  &   $\Sym_{2s}(\R)^\times_{(s,s)}$ \\
  (NS1) & $\sp_{4s}(\R)$& $\sp_{2s}(\R)^{\oplus 2}$ & $\sp_{2s}(\C)$
  &  $\sp_{2s}(\R)$ & $\Sp_{2s}(\R)$           & \\
  \hline
\hline 
$\Herm_r(\C)$\phantom{\Big (} & $\su_{r,r}(\C)$ && &&& \\
\hline   
  (C)  & $\su_{r,r}(\C)^{\oplus 2}$& $\su_{r,r}(\C)$ & $\su_{r,r}(\C)$
  &  $\fsl_r(\C)\oplus \R$&& \\
(P) & $\su_{2s,2s}(\C)$ &
        $\fu_{s,s}(\C)^{\oplus 2}/\R (i\1,i\1)$ & $\fsl_{2s}(\C)+ \R \1$
  &  $\su_{s,s}(\C)$&
 $\U_{s,s}(\C)$ & $\Herm_{2s}(\C)^\times_{(s,s)}$ \\
  (S1)  & $\su_{2s,2s}(\C)$ & $\so^*(4s)$ & $\so_{2s,2s}(\R)$  &  $\so_{2s}(\C)$ &
$\SO^*(4s)/\SO_{2s}(\C)$     &\\
  (NS2)  &$\su_{2r,2r}(\C)$ & $\sp_{4r}(\R)$  & $\fu_{r,r}(\bH)$  &  $\sp_{2r}(\C)$
 & $\Sp_{4r}(\R)/\Sp_{2r}(\C)$ &  \\
  \hline
\hline 
$\Herm_r(\H)$\phantom{\Big (} & $\so^*(4r)$ && &&& \\
  \hline
 (C) &$\so^*(4r)^{\oplus 2}$ & $\so^*(4r)$ & $\so^*(4r)$  &  $\gl_r(\H)$& & \\
(P) & $\so^*(8s)$& $\fu_{2s,2s}(\C)$ & $\gl_{2s}(\H)$ 
            &  $\fu_{s,s}(\H)$  & $\U_{2s,2s}(\C)/\U_{s,s}(\H)$ & $\Herm_{2s}(\H)^\times_{(s,s)}$\\
(S2) &$\so^*(8s)$ & $\so^*(4s)^{\oplus 2}$ & $\so_{4s}(\C)$  &  $\so^*(4s)$ &
 $\SO^*(4s)$        &\\
  \hline
\hline 
$\Herm_3(\bO)$ \phantom{\Big (}& $\fe_{7(-25)}$ &&& & &\\
\hline   
  (C)  & $\fe_{7(-25)}^{\oplus 2} $& $\fe_{7(-25)}$  & $\fe_{7(-25)}$
  &  $\fe_{6(-26)} \oplus \R$ &&
 \\
  \hline   
\hline   
$\R^{1,d-1}, d \geq 3$ \phantom{\Big (}& $\so_{2,d}(\R)$ && &&& \\
\hline   
(C)  & $\so_{2,d}(\R)^{\oplus 2}$ & $\so_{2,d}(\R)$ & $\so_{2,d}(\R)$  &  $\R \oplus \so_{1,d-1}(\R)$
 && \\
 (S4), $p \geq 2$ &$\so_{2,d}(\R)$ &$\so_{2,p-1}(\R)$ & $\so_{1,p}(\R)$ 
                             &$\so_{1,p-1}(\R)$
 & $(\AdS^{p} \times \bS^{q})/\{\pm\1\}$ & $\dS^{p} \times \bHy{q}$ \\
 \quad $p+q=d$ &&  $\quad \times \so_{q+1}(\R)$ &$\quad\times \so_{1,q}(\R)$
  &  $\quad \times \so_q(\R)$
&  &  \\  
(NS3), $q=0$ &$\so_{2,d}(\R)$ &$\so_{2,d-1}(\R)$ & $\so_{1,d}(\R)$ 
                             &$\so_{1,d-1}(\R)$
 & $\AdS^{d}$ & $\dS^d$  \\
(P), $j =1$  &$\so_{2,d}(\R)$  & $\R \oplus \so_{2,d-2}(\R)$
 & $\R \oplus \so_{1,d-1}(\R)$  &  $\so_{1,d-2}(\R)$
                                         &$(\AdS^{d-1} \times \bS^1)/\{\pm\1\}$
     &$\dS^{d-1} \times \R$ \\
\hline\hline 
\end{tabular}\\[2mm] \begin{center}
  {Table 4: Modular causal Makarevi\v c spaces}\end{center}

\section{The Lorentzian case}
\mlabel{sec:lorentz} 

As Lorentzian causal manifolds are particularly relevant
for applications ins physics, we now discuss
the Lorentzian Makarevi\v c spaces in some detail.
It is well-known that the Lorentzian 
symmetric spaces   \[     \AdS^k \times \bS^{d - k}, \quad 
   \dS^k \times  \bHy{d - k},\quad k = 1, \ldots, d\]
  are conformally flat, i.e., locally conformal
  to $d$-dimensional Minkowski space  $\R^{1,d-1}$
  (cf.\ \cite[\S XI.5.3]{Be00}).
  Therefore it is no surprise that (up to coverings),
  these are precisely the Lorentzian Makarevi\v c  spaces. 
  For $k = 1$, we obtain in particular
  $\bS^1 \times \bS^{d-1}$, a two-fold covering
  of the conformal completion $M$ of $\R^{1,d-1}$,
  and $\R \times \bHy{d-1}$, considered  as causal symmetric spaces.

  We start this section by describing the conformal completion
  of $\jV = \R^{1,d-1}$ as the isotropic quadric $Q(\R^{2,d})$
  of isotropic lines in $\R^{2,d}$ (Subsection~\ref{subsec:7.1}).
  The simple classification of the involutive automorphisms of the
  Jordan algebra $\R^{1,d-1}$ is given in Subsection~\ref{subsec:minkowski},
  and this leads to simple descriptions of the corresponding
  spaces $M^{(\pm \alpha)} \subeq Q(\R^{2,d})$.
  In Subsection~\ref{subsec:7.3} we describe the boundary
  orbits of the $M^{(\pm \alpha)}$ in $Q(\R^{2,d})$, which is motivated in
  particular by potential applications to holographic aspects
  of AQFT. For the two most important examples, 
  de Sitter space and anti-de Sitter space, we describe
  the wedge regions in Subsection~\ref{subsec:desitter}
  and Subsection~\ref{subsec:antidesitter} and show that 
 the wedge regions in anti-de Sitter space
  are not globally hyperbolic (Proposition~\ref{prop:wads-incl}).

\subsection{The conformal completion of Minkowski space}
\mlabel{subsec:7.1}
For the conformal Lie algebra $\g = \so_{2,d}(\R)$
of $\jV = \R^{1,d-1}$, we consider $\tilde\jV := \R^{2,d}
= \R \oplus \jV \oplus \R$, endowed with the symmetric bilinear form
\[ \tilde\beta((t,v,s), (t',v',s')) := tt' + \beta(v,v') - ss'.\]
We write $\be_1, \ldots,\be_{d+2}$ for the canonical basis vectors
and consider the Euler element $h \in \so_{2,d}(\R)$, defined by 
\begin{equation}
  \label{eq:h-lorentz}
 h \be_1 = \be_{d+2}, \quad 
  h \be_{d+2} = \be_1 \quad  \mbox{ and } \quad
  h \be_j = 0, \ j = 2,\ldots, d-1.
\end{equation}
Its eigenspaces are
\[ \tilde\jV_1(h) = \R (\be_1 + \be_{d+2}), \quad 
  \tilde\jV_{-1}(h) = \R (\be_1 - \be_{d+2}) \quad  \mbox{ and } \quad
  \tilde\jV_0(h) = \Spann \{  \be_2, \ldots, \be_{d+1}\}.\] 
For $v = (v_0, \bv) \in \jV$, we write $v^\flat = (-v_0, \bv)$,
considered as a row matrix.
Then a short calculation yields
\[  \g_1(h) = \Bigg\{\pmat{
    0 &  v^\flat & 0 \\ 
    -v & 0 & v \\ 
    0 & v^\flat & 0} \: v \in \jV \Bigg\}, \]
where the $3 \times 3$-block structure corresponds to the
decomposition $\tilde\jV = \R \be_1 \oplus \jV \oplus \R \be_{d+2}$.
We write
\[ X_v := \pmat{    0 & - v^\flat & 0 \\     v & 0 & -v \\     0
    & -v^\flat & 0}  \in \g_1(h),
  \quad \mbox{ so that } \quad
 X_v^2 = \pmat{   -\beta(v,v) & 0 & \beta(v,v) \\
   0 & 0 & 0 \\   -\beta(v,v) & 0 & \beta(v,v)}\]
leads to 
\[ e^{X_v}(\be_1 - \be_{d+2})
  = (1 -\beta(v,v)) \be_1 + 2 v - (1+ \beta(v,v)) \be_{d+2}.\]

We realize the conformal completion $M$ of $\jV$ as the quadric
\begin{equation}
  \label{eq:quadtildev}
 Q := Q(\R^{2,d}) := \{ [\tilde v] \in \bP(\tilde\jV) \:
 \tilde\beta(\tilde v, \tilde v) = 0\}
\end{equation}
with the open dense embedding
\begin{equation}
  \label{eq:eta-flow}
\eta \: \jV \to Q, \quad  \eta(v) := \Big[ \frac{1 - \beta(v,v)}{2} : v : - \frac{1 + \beta(v,v)}{2}\Big] = \exp(X_v).[\be_1 - \be_{d+2}] \in \bP(\tilde\jV), 
\end{equation}
corresponding to the action of the translation group
$(\jV,+) \cong \g_1(h)$ on $Q$. Its range is 
\[  \eta(\jV) = \{  [t:v:s] \in \bP(\R \oplus \jV \oplus \R) \: t -s \not=0\}.\]
We also note that
\begin{equation}
  \label{eq:euler-flow}
  \exp(th).\eta(v) = \eta(e^t v) \quad \mbox{ for }  \quad
  v \in \jV, t \in \R
  \quad \mbox{ and } \quad
  \exp(\pi i h).\eta(v) = \eta(-v).
\end{equation}

  \begin{rem} The eigenspace $\tilde\jV_{-1}(h) = \R (\be_1- \be_{d+2})$
    is isotropic and invariant under the parabolic subgroup~$Q_h$.
    Maximality of the subgroup $Q_h \subeq G = \SO_{2,d}(\R)_e$
    implies that it coincides with the stabilizer
    of this eigenspace, so that we obtain an embedding
\[ M \cong G/Q_h \cong G.\tilde\jV_{-1}(h) =: Q(\R^{2,d}) \subeq \bP(\R^{2,d}).\]
  \end{rem}

  \begin{rem} Note that $Q(\R^{2,d})\setminus \eta(\R^{1,d-1})$
    is the set of all rays $[x_1:\ldots: x_{d+2}]$ with
    $x_1 = x_{d+2}$. For $x_1 = x_{d+2} = 0$, this specifies
    a projective quadric $Q(\R^{1,d-1}) \cong \bS^{d-2}$ of Minkowski space,
    and for $x_1 = x_{d+2} \not=0$, we may normalize
    to $x_1 = x_{d+2} = 1$, which leads to a manifold diffeomorphic to
    the light cone $\cN = \{ x \in \R^{1,d-1} \: \beta(x,x) = 0\}$.
  \end{rem}

\subsection{The Jordan involutions on Minkowski space} 
\mlabel{subsec:minkowski} 

  For the $d$-dimensional Minkowski space
  $\jV = \R^{1,d-1}$ with Jordan unit $\be_0$, we have 
\[ \Aut(\jV) \cong \OO_{1,d-1}(\R)^{\be_0} \cong \OO_{d-1}(\R). \] 
Involutions in this group 
are orthogonal reflections $\sigma$
in subspaces $F_\sigma \subeq \R^{d-1}$, and their  
conjugacy classes are determined by 
$\dim F_\sigma \in \{0,\ldots, d-1\}$. We write 
$r_j \in \Aut(\jV)$ for the involutive automorphism with
\[ r_j = I_{d-j,j} := \diag(\1_{d-j}, - \1_{j}) \in \OO_{1,d-1}(\R)^{\be_0}, \quad
  j = 0,\ldots, d-1.\]

Here $r = \rk \jV = 2$ and for the Jordan frame
\[ c_1 = \be_0 + \be_1, \quad  c_2 = \be_0 - \be_1, \]
there is only one Pierce space
$\jV_{1,2} = \Spann \{\be_2, \ldots, \be_{d-1}\}$. Therefore the only
non-trivial Pierce involution is
\[ \alpha_1 = r_{d-2} = \diag(1,1,-1,\ldots, -1).\]
The trivial ones are $\alpha_0 = \alpha_2 = \id_\jV$.

Deviating from our convention, we write in this section
$G := \SO_{2,d}(\R)_e$, which for $d$ even is a $2$-fold
covering of the conformal group~$\Co(\jV)_e$.
In  $G$, the involution
$\tau_h= - \id_\jV$ is implemented by
\[ \tau_h^G = \diag(-1,1,\ldots, 1,-1) \in \SO_{2,d}(\R)\]
(\cite[\S~17.4.2]{HN12}), and $\theta$ corresponds to
\[ \theta^G = \exp(\pi z_\fk) = \diag(-1,-1,1,\ldots, 1).\]
We may therefore implement $\theta_{-r_j} =\tau_h^G \theta^G \sigma_{r_j}$
by conjugation with 
\[\hat r_j = \diag(1,-1, \1_{d-j-1}, -\1_{j+1}) \in \OO_{2,d}(\R),\]
We then obtain 
\begin{equation}
  \label{eq:mink-gminusalpha}
\so_{2,d}(\R)^{(-r_j)}= \so_{2,d}(\R)^{\hat r_j}
  \cong
  \begin{cases}
\so_{1,d}(\R)  & \text{ for } j = d -1 \ \ \text{(NS3)} \\     
\so_{1,1}(\R) \oplus \so_{1,d-1}(\R)  & \text{ for } j = d -2  \  \ \text{(P)}\\ 
\so_{1,d-1-j}(\R) \oplus \so_{1,j+1}(\R)  & \text{ for } 0 \leq j  < d-2
\ \ \text{(S4)}.    
\end{cases}
\end{equation}

 For $p = j+1= 1,\ldots, d$, and $q := d-p$,
 we can thus identify the spaces $M^{(-r_j)}$ in the quadric 
  $Q(\R^{2,d})$ as follows. We  use 
  \eqref{eq:mink-gminusalpha} to obtain in 
  $G = \SO_{2,d}(\R)_e$ the subgroups 
\[ G^{(-r_j)} \cong \SO_{1,q}(\R)_e \times \SO_{1,p}(\R)_e, \]
which leads to the non-compactly causal symmetric spaces 
\[ M^{(-r_j)} := G^{(-r_j)}.\eta(0)
  = [\SO_{1,q}(\R)_e.\be_1 - \SO_{1,p}(\R)_e.\be_{d+2}]
  \cong   \bHy{q} \times \dS^p\]
(note the rearrangement of the indices, where $\be_2$ is exchanged with
$\be_{q+2}$). 
Here $\dS^p \subeq \R^{1,p}$ is \break $p$-dimensional {\it de Sitter space} 
and $\bHy{q} \subeq \R^{1,q}$ denotes $q$-dimensional hyperbolic space,
an irreducible Riemannian symmetric space.
\begin{footnote}{We use this notation, and not the more customary~$\bH$,
    to distinguish it from the skew-field
    $\bH$ of quaternions.}\end{footnote}

The involution $\theta_{r_j} = \theta^G \sigma_{r_j}$
can be implemented by conjugation with 
$\diag(-1,-1, \1_{d-j-1}, -\1_j, 1)$, which
is conjugate to 
\[\tilde r_j := \diag(-\1_{j+2}, \1_{d-j}). \]
We then obtain
\begin{equation}
  \label{eq:mink-galpha}
 \so_{2,d}(\R)^{(r_j)}= \so_{2,d}(\R)^{\tilde r_j}
  \cong
  \begin{cases}
\so_{2,d-1}(\R)  & \text{ for } j = d -1 \ \ \text{(NS3)} \\     
\so_{2,d-2}(\R) \oplus \so_2(\R) & \text{ for } j = d -2  \  \ \text{(P)}\\ 
 \so_{2,j}(\R) \oplus \so_{d-j}(\R)  & \text{ for } 0 \leq j  < d-2
\ \ \text{(S4)}.    
\end{cases}
\end{equation}
We obtain for $p = j+1 = 1,\ldots, d$ and $q = d- p$ 
the corresponding compactly causal symmetric spaces 
\begin{equation}
  \label{eq:mrj}
 M^{(r_j)} := G^{(r_j)}.\eta(0)
  = [\SO_{2,p-1}(\R)_e.\be_1 - \SO_{q+1}(\R)_e.\be_{d+2}]
  \cong (\AdS^p\times \bS^q)/\{\pm \1\},
\end{equation}
 where $\AdS^p \subeq \R^{2,p-1}$ is $p$-dimensional {\it anti-de Sitter space}. 
For $p = d = j+1$ we obtain with $\bS^0 \cong \{ \pm \1\}$ that
\begin{equation}
  \label{eq:mrj2}
  M^{(r_{d-1})} \cong \AdS^d.
\end{equation}

In all cases, the ``intersection with $\jV$'', is given by 
\begin{equation}
  \label{eq:malphaetav}
 M^{(\pm \alpha)} \cap \eta(\jV)
 = \{ [x] \in M^{(\pm \alpha)} \: x_1 - x_{d+2} \not=0\}.
\end{equation}
%

\subsection{Conformal boundaries}
\mlabel{subsec:7.3}

The compactness of $Q = Q(\R^{2,d})$ implies in particular that the
closure of $M^{(\pm \alpha)}$ is compact and invariant under
the group $G^{(\pm \alpha)}$. We now describe the
boundary orbits of this compactification, having in particular
applications to holographic constructions in AQFT in mind.

\begin{prop} \mlabel{prop:boundaries}
  The boundaries of $M^{(\pm \alpha)}$ in $Q(\R^{2,d})$ are as
  follows:
  \begin{itemize}
  \item[\rm(CC)] The open subset
    $M^{(\alpha)} \cong (\AdS^p\times \bS^q)/\{\pm\1\},
    q = 0,\ldots, d-1,$ $p = d - q$, 
    is dense in $Q(\R^{2,d})$ and its boundary
    is a $(p-1)$-dimensional quadric $Q(\R^{2,p-1})$.
\item[\rm(NCC)] The open subset
  $M^{(-\alpha)} \cong \dS^p\times \bHy{q}$, $q = 0,\ldots, d-1$, $p = d - p$, 
  is not dense.   Its boundary is a union 
  \[ \bS^{p-1}\ \dot\cup \ (\bS^{p-1} \times \R^\times \times \bS^{q-1}),\]
  where we put $\bS^{-1} = \eset$. 
  The other open orbit of $G^{(-\alpha)}$ is 
$\bHy{p} \times \dS^{q}$. 
  \end{itemize}
\end{prop}

\begin{prf}
  (CC)  The image of $\AdS^p\times \bS^q$ in
  $Q(\R^{2,d})$ consists of pairs $[\bv_1: \bv_2]$ with
$\bv_1 \in \R^{2,p-1}$ and $\bv_2 \in \R^{q+1}$ satisfying
$\tilde\beta(\bv_1, \bv_1) = -\tilde\beta(\bv_2, \bv_2) = 1$. 
As this is equivalent to $\bv_2 \not=0$, 
this is an open dense submanifold. 
Its complement, consisting of projective pairs with $\bv_2 = 0$,
is the $(p-1)$-dimensional quadric $Q(\R^{2,p-1})$.

\nin (NCC) The image of $\dS^p\times \bHy{q}$ 
consisting (after rearrangement of the coordinates
$\be_2 \leftrightarrow \be_{p+2}$)
of pairs $[\bv_1: \bv_2]$ with
$\bv_1 \in \R^{1,p}$ and 
$\bv_2 \in \R^{1,q}$ satisfying
$\tilde\beta(\bv_1, \bv_1) = -\tilde\beta(\bv_2, \bv_2) = -1$.
If $\tilde\beta(\bv_2, \bv_2) > 0$, we may normalize
to $1$ and obtain elements of $\dS^p\times \bHy{q}$. 

For $\tilde\beta(\bv_2, \bv_2) = 0$ and $\bv_2 \not=0$,
we obtain the boundary manifold
\[ (\R^\times \bS^{p-1} \times \R^\times \bS^{q-1})/\R^\times 
  \cong \bS^{p-1} \times \R^\times\times \bS^{q-1}, \]
and for $\bv_2= 0$, we get $\bS^{p-1}$.
Since there are also elements with 
$\tilde\beta(\bv_2, \bv_2) < 0$, the
embedding of $\dS^p \times \bHy{q}$ is not dense.
The interior of its complement is a causal embedding of 
$\bHy{p} \times \dS^q$, on which the
group $\SO_{1,p}(\R)_e \times \SO_{1,q}(\R)_e$ acts by isometries. 
\end{prf}

The following observation is the geometric key to the
celebrated $\AdS/\CFT$-correspondence in QFT. 

\begin{cor} The boundary of $\AdS^d$ in $Q(\R^{2,d})$
  is the $(d-1)$-dimensional quadric $Q(\R^{2,d-1})$, i.e.,
  compactified $(d-1)$-dimensional Minkowski space, 
  and the boundary of $\dS^d$ is the sphere $\bS^{d-1}$.
\end{cor}

\subsection{De Sitter space} 
\mlabel{subsec:desitter}

We consider the $d$-dimensional de Sitter space
$M^{(-r_{d-1})} \cong \dS^d$ (type (NS3) in \eqref{eq:mink-gminusalpha}).
Applying \eqref{eq:conformal-flat} in Appendix~\ref{subsec:antidesitter}
to $(\jV,-\beta)$,
for which  $\dS^d$ is the unit sphere, we obtain a diffeomorphism
\begin{equation}
  \label{eq:stereo-dS}
 \dS^d \setminus \{ x_d = 1\} \to \jV \setminus S(\jV,\beta)
 = \{  v \in \jV \: \beta(v,v) \not=1\}.
\end{equation}
In particular, 
$\eta^{-1}(M^{(-r_{d-1})}) = \eta^{-1}(\dS^d)$ is the complement of 
$S(\jV, \beta)$, hence has 
$3$ connected components. 
The inverse of this map is given,
for  $\beta(v,v)\not=1$, by 
\eqref{eq:eta-flow}
\[ v \mapsto (y_0, \by),\quad
  \eta(v) = \Big[ 1 : \frac{2v}{1 - \beta(v,v)} : \frac{\beta(v,v)+1}
  {\beta(v,v)-1}\Big] = [1: y_0: \by]
  \quad \mbox{ with }\quad
  (y_0,\by) \in \jV^\vee\]
(notation from Appendix~\ref{subsec:stereo}). 
In these coordinates the map
\eqref{eq:stereo-dS} corresponds to the
stereographic projection with center~$\be_d \in \dS^d$. 
 Our base point in $Q$ is $[\be_1 - \be_{d+2}]$, corresponding to
  $-\be_d  \in \dS^d \subeq \R^{1,d}$.

 Note that
  \[ \eta(S(V,\beta)) = \{ [0:v:1] \: \beta(v,v) = 1 \}
    \cong \bHy{d-1}.\]

  The embedding of $\dS^d$ in $Q$ 
  is also compatible with the canonical causal orientations: 
  In $\dS^d \subeq \jV^\vee$,
  the curve $\gamma(t) = -\be_{d+2} + t \be_2$ is causal,
  and the corresponding curve in $\jV$ has a tangent vector in $0$
  that is a positive multiple of $\be_0 \in \R^{1,d-1}$.

  \begin{prop} \mlabel{prop:wds-incl}
  We consider the Euler element $h' \in \so_{1,d}(\R) \subeq \so_{2,d}(\R)$
  specified by 
  \begin{equation}
    \label{eq:h'dS}
 h' \be_2 = -\be_{d+2}, \quad 
  h' \be_{d+2} = -\be_2 \quad  \mbox{ and } \quad
  h' \be_j = 0\quad \mbox{ for } \quad  \ 2 < j < d+2. 
\end{equation}
Its positivity region $W_{\dS^d}^+(h')$ has the following properties:
  \begin{itemize}
  \item[\rm(a)] $W_{\dS^d}^+(h')$ 
    is contained in $\eta(\jV) \cong \jV$, where it corresponds to
    the double cone $\cD_{\be_0, -\be_0}$. In particular, it is causally convex
    and therefore globally hyperbolic. 
  \item[\rm(b)] The group 
    $(\SO_{2,d}(\R)^{h'})_e \cong \SO_{1,1}(\R)_e \times \SO_{1,d-1}(\R)_e$
    acts transitively on $W_{\dS^d}^+(h')$ and turns it into a
    non-compactly causal Lorentzian symmetric space
    diffeomorphic to $\R \times \bHy{d-1}$.
  \end{itemize}
\end{prop}

\begin{prf} (a) As $h'(-\be_{d+2}) = \be_2$ is positive timelike, $h'$
  is a causal Euler element (cf.\ Definition~\ref{def:reductive-ncc}).
    Its positivity region is
\[ W_{\dS^d}^+(h') = \{ y \in \dS^d \subeq \R^{1,d}\: -y_d > |y_0|\}
\subeq \eta(\jV). \]
The inequality $-y_d > |y_0|$ corresponds to
  $-x_{d+2} > |x_2|$ in the $\R^{2,d}$-coordinates.
  It implies in particular that 
$x_{d+2} < 0$, so that $x_{d+2} \not=1$, i.e., 
$W_{\dS^d}^+(h') \subeq \eta(\jV).$

We consider the curve
  \[ \zeta \: (-1,1) \to \dS^d \subeq Q(\R^{2,d}), \quad
    \zeta(t) = \eta(t \be_0)
    = \Big[\frac{1-t^2}{2} : t \be_2 : -\frac{1 + t^2}{2}\Big] 
    = \Big[1 : \frac{2t}{1-t^2} \be_2 : -\frac{1 + t^2}{1-t^2}\Big].\]
  Then $\zeta(t)$ corresponds to the curve
  \[ (-1,1) \to \dS^d, \quad
    t \mapsto  \frac{2t}{1-t^2} \be_2 +\frac{t^2+1}{t^2-1}\be_{d+2}\in \dS^d.\]
  For $t := \tanh(s/2)$, we have
$\zeta(t) = [1:\sinh(s) \be_2 : - \cosh (s)] = e^{sh'}\be_{d+2},$ 
  so that $\zeta((-1,1))$ is the causal geodesic $e^{\R h'}.(-\be_{d+2})
  \subeq \dS^d$.

  Its causally convex hull in $\dS^d$ is the positivity region
  $W_{\dS^d}^+(h')$ (\cite{MNO24}),
  and the causally convex hull of $(-1,1)\be_0$ in
  Minkowski space $\jV$ is the double cone $\cD_{\jV} =
  \cD_{\be_0, -\be_0}$. This implies   that 
  \[ \eta(\cD_{\be_0, -\be_0}) = W_{\dS^d}^+(h')\]
  (cf.\ Theorem~\ref{thm:pos-dom-ncc}). 

  \nin (b) It follows that the positivity region of the Euler element $h'$ 
  is diffeomorphic to $\R \times \cD_{\R^{d-1}}$,
  where $\cD_{\R^{d-1}} \subeq \R^{d-1}$ is the open unit ball,
  which is also a model of hyperbolic space $\bHy{d-1}$.
  This is a homogeneous space of the group 
  \[  (\SO_{2,d}(\R)^{h'})_e \cong \SO_{1,1}(\R)_e \times \SO_{1,d-1}(\R)_e,\]
  which acts isometrically with respect to the natural
  non-compactly causal symmetric space structure on
  $\R \times \bHy{d-1}$.   
\end{prf}

\subsection{Anti-de Sitter space} 
\mlabel{subsec:antidesitter}

We consider the $d$-dimensional anti-de Sitter space
\[ M^{(r_{d-1})} \cong \AdS^d \cong \SO_{2,d-1}(\R)_e.[1:0:-1]
\subeq Q(\R^{2,d})\]
(type (NS3) in \eqref{eq:mink-galpha}).
In view of \eqref{eq:conformal-flat} in Appendix~\ref{subsec:stereo}, 
$\eta^{-1}(M^{(r_{d-1})})$ is the complement of
$S(\jV, -\beta) \cong \dS^{d-1}$, hence has
$2$ connected components for $d > 2$.
For $\beta(v,v)\not=-1$, we obtain from
\eqref{eq:eta-flow}
\[ \eta(v) = \Big[ \frac{1-\beta(v,v)}{1+ \beta(v,v)}
  : \frac{2v}{1 + \beta(v,v)} : -1\Big] = [x_1:x_2: \cdots : x_{d+1}: -1] \]
with $x_1^2 + x_2^2 - x_3^2 - \cdots - x_{d+1}^2= 1$,
i.e., $\beta^\wedge(x,x) =1$   
(cf.~Appendix~\ref{subsec:stereo} for  $\beta^\wedge$).
  In these coordinates, the rational map
  \[ \AdS^d \setminus \{ x_1 = -1\}  \to \eta(\jV) \cong \jV\]
  corresponds to the
  stereographic projection discussed in Appendix~\ref{subsec:stereo},
  applied to the bilinear form $\beta^\wedge$ on $\jV^\wedge \cong \R^{2,d-1}$,  and
  \begin{equation}
    \label{eq:eta-inv-ads}
 \eta^{-1}(\AdS^d) = \jV \setminus \dS^{d-1}
 = \{ v\in \jV \: \beta(v,v) \not= -1\}.
  \end{equation}

  Our base point is $[\be_1 - \be_{d+2}]$, corresponding to
  $\be_1  \in \AdS^d \subeq \R^{2,d-1}$. 
  This embedding of $\AdS^d$
  is also compatible with the canonical causal orientations,
  determined by $[\be_2] \in T_{\be_1}(\AdS^d)$ being positive
  (\cite[\S 11]{NO23a}). 
  In $\AdS^d\subeq \jV^{\wedge} \cong \R^{2,d-1}$,
  the curve $\gamma(t) = \be_1 +  t \be_2$ is causal,
  and the corresponding curve in $\jV$ has a tangent vector in $0$
  that is a positive multiple of the basis vector $\be_2
  \in \Spann \{  \be_2, \ldots, \be_{d+1}\} \cong \jV$.

\begin{prop} \mlabel{prop:wads-incl} 
  We consider the Euler element $h'' \in \so_{1,d-1}(\R) \subeq\so_{2,d}(\R)$, 
  specified by
  \begin{equation}
    \label{eq:h''euler}
 h'' \be_2 = \be_{d+1}, \quad 
  h'' \be_{d+1} = \be_2 \quad  \mbox{ and } \quad
  h'' \be_j = 0, \ j \not=2, d+1.
  \end{equation}
Its positivity region $W_{\AdS^d}^+(h'')$ has the following properties:
  \begin{itemize}
  \item[\rm(a)] It  is contained
    in the open dense subset $\eta(\jV) \subeq M$
    and $W^+_{\AdS^d}(h'') \cong W^+_\jV(h'') \setminus \dS^{d-1}$. 
  \item[\rm(b)] It has two connected components. 
  \item[\rm(c)] As causal manifolds, its connected components
    are not globally hyperbolic. 
  \end{itemize}
\end{prop}

\begin{prf} (a) First we observe that 
\[ W_{\jV}^+(h'') \cap \eta^{-1}(\AdS^d) 
  = \{ (y_1, y_2, \cdots, y_{d+1}) \in \jV \setminus \dS^{d-1}
  \: y_{d+1} > |y_2| \} \]
is not convex. 

Next we observe that the positivity region of $h''$ in $\AdS^d$
is contained in its timelike domain, specified by the inequality 
$x_{d+1}^2 > x_2^2.$
If $x \in \AdS^d \setminus \eta(\jV)$, then $x_1 = -1$, so that
\[ 1
  = x_1^2 +x_2^2 - x_3^2 - \cdots - x_{d+1}^2 
  = 1 +x_2^2 - x_3^2 - \cdots - x_{d+1}^2 \]
leads to $x_2^2 =  x_3^2 + \cdots +  x_{d+1}^2,$ 
and hence to $x_{d+1}^2 \leq x_2^2$. 
Therefore $x$ cannot be timelike for $h''$.

\nin (b) By \cite[Lemma~11.2]{NO23a}, the positivity region of $h''$
in $\AdS^d$ has two connected components. We we also  see this by
applying the inverse of $\eta$ and Proposition~\ref{prop:wads-incl}: 
\[ W_{\AdS^d}^+(h'') \cong
  W_\jV(h'') \setminus \dS^d
  = \{ y \in \R^{1,d-1} \: y_{d-1} > |y_0|\} \setminus \dS^{d-1}.\]
That this set has exactly two components follows from the fact that
every ray in the Rindler wedge $W_{\jV}^+(h'')$ intersects
$\dS^{d-1}$ in exactly one point.

\nin (c) As both connected components are isomorphic
(cf.\ \cite[Lemma~11.2]{NO23a}), it suffices to show that one of them
is not globally hyperbolic.  The component
  \[ \cO := \{ y \in \R^{1,d-1} \: y_{d-1} > |y_0|, -1 < \beta(y,y) < 0\} \]
  is easily seen to be not globally hyperbolic.
  For $0 < \eps <  1$, it contains the two elements
  $\be_1 \pm \eps \be_0$ and the  $\cO$-causal interval 
  $\cD_{\be_1 + \eps \be_0, \be_1 - \eps \be_0}$
  between these two points clusters in the boundary point~$\be_1$,
  hence does not have compact closure in $\cO$.   
\end{prf}

\begin{ex} \mlabel{ex:2d-ads} We take a closer look at the $2$-dimensional case.
We consider the lightlike basis 
  \[ \bc_1 := \frac{\be_0 + \be_1}{2} \quad \mbox{ and }  \quad 
    \bc_2 := \frac{\be_0 - \be_1}{2} \]
  of $\R^{1,1}$ and write $(y_1, y_2)$ for the corresponding coordinates;
  called {\it lightray coordinates}.   Then
  \[ x_0 = \frac{y_1 + y_2}{2}, \quad x_1 = \frac{y_1 - y_2}{2}
  \quad \mbox{ and } \quad  y_0 = x_0 + x_1, \quad y_1 = x_0 - x_1.\]
  In lightray coordinates, we consider the conformal
  transformation $(y_0, y_1) \mapsto (y_0', y_1') = (y_0, - y_1^{-1}).$ 
  Then
  \[ x_0' + x_1' = x_0 + x_1 \quad \mbox{ and }  \quad
    x_0' - x_1' = \frac{1}{x_1 - x_0} \]
  leads to
  \[ x_0'
    = \frac{1}{2}\Big(x_0 + x_1 + \frac{1}{x_1-x_0}\Big)
= \frac{x_1^2 - x_0^2 + 1}{2(x_1 - x_0)}
= \frac{1 - \beta(x,x)}{2(x_1 - x_0)} 
= \frac{ \beta(x,x)-1}{2(x_0 - x_1)} \]
and
\[ x_1'
  = \frac{-1 - \beta(x,x)}{2(x_1 - x_0)} 
= \frac{\beta(x,x) + 1}{2(x_0-x_1)}. \]
The  conformal maps $\psi(x) = (x_0', x_1')$ maps $\dS^1$ to the $x_0$-axis:
\[ \psi(x)
  = \Big(\frac{1}{x_1 - x_0},0\Big)
  = (x_0 + x_1, 0),\]
and in particular
$\psi(\sinh t, \pm\cosh t) = (\pm e^{\pm t},0).$ 
In lightray coordinates it is clear that the Rindler wedge $W_R$ satisfies 
$\psi(W_R) = \jV_+,$ so that
\[ \psi(W_R \setminus \dS^1) = \jV_+ \setminus \R \be_0.\]
In this picture the two connected components of 
$W_{\AdS^2}^+(h'')$, the left and right half of
$\jV_+ \setminus \R \be_0$, are not globally hyperbolic. 
\end{ex}

\begin{ex} (Globally hyperbolic domains in $\AdS^d$; the elliptic type)
\mlabel{ex:ads-ball}  
  The double cone
  $\cD_{\be_0, - \be_0} \subeq~\jV$ is contained in the intersection
  of $\jV \subeq Q(\R^{2,d})$ with $\AdS^d$
  (Propositions~\ref{Dv-inc}), 
  which corresponds to
  $\jV \setminus \dS^{d-1}$. It is a globally hyperbolic domain in $\jV$,
  so that its image under $\eta$ in $\AdS^d$ also has this property.

  We have
  \[ \eta(\pm \be_0)\ {\buildrel\eqref{eq:eta-flow}\over =}\ [0: \pm \be_0: -1] = [\pm\be_2 - \be_d]
    \sim \pm \be_2 \in \AdS^d \subeq \R^{2,d-1}\]
  (notation from Appendix~\ref{subsec:stereo}). 
  Further, for $|t| < 1$, 
  \[ \eta(t\be_0)
    = \Big[ \frac{1 - t^2}{1 + t^2} : \frac{2t \be_0}{1 + t^2} :    -1\Big]
    = \Big[ \frac{1 - t^2}{1 + t^2} \be_1 + \frac{2t}{1 + t^2} \be_2
    - \be_{d+2}\Big]
    \sim  \frac{1 - t^2}{1 + t^2} \be_1 + \frac{2t}{1 + t^2} \be_2
    \in \AdS^d \]
parametrizes the half circle
$\{ x_1 \be_1 + x_2 \be_2 \: x_1 > 0, x_1^2 + x_2^2 = 1\}
    \subeq \AdS^d$ 
  connecting $\eta(-\be_0) = -\be_2$ by a causal curve with
  $\eta(\be_0) = \be_2$. Therefore
  $\eta(\cD_{\be_0, -\be_0}) \subeq \AdS^d$ is the open domain
  of all points on causal curves from $-\be_2$ to $\be_2$.

We recall from \cite[\S 11]{NO23a} 
the exponential function of the symmetric space $\AdS^d$:
\begin{equation}
  \label{eq:expfunct-ads}
  \Exp_p(y) = \cos(\sqrt{\beta^\wedge(y,y)}) p +
  \frac{\sin(\sqrt{\beta^\wedge(y,y)})}{\sqrt{\beta^\wedge(y,y)}} y.
\end{equation}
Timelike vectors $y \in T_p(\AdS^d)$ 
generate closed geodesics and $\beta^\wedge(y,y) = \pi^2$ implies
$\Exp_p(y) = -p.$ In particular,
all timelike geodesics through $-\be_2$ also pass   through $\be_2$. 
As all points in $\eta(\cD_{\be_0, - \be_0})$ 
lie on geodesics, we obtain 
  \[ \eta(\cD_{\be_0, -\be_0})
    = \Exp_{-\be_2}(\{y = (y_1, 0, \by) \: y_1 > 0,
 0 <  \beta^\wedge(y,y)  < \pi^2\}).\]
\end{ex}

\subsection{An application to the general case}
\mlabel{subsec:7.6}

We have seen in Proposition~\ref{prop:wads-incl}  that
the positivity regions $W_{\AdS^d}^+(h'')$ in anti-de Sitter space
is not globally hyperbolic. We now show that
the same is true for all modular compactly causal
spaces~$M^{(\alpha)}$. We start with a remark that will be needed
in the proof.

\begin{rem} (a) \mlabel{rem:proper} Let $M = G/H$ be a semisimple symmetric space,
  on which the connected Lie group~$G$ acts faithfully.
  In particular, the corresponding symmetric Lie algebra
  $(\g, \tau)$ is effective. We choose a Cartan involution
  $\theta$ commuting with $\tau$, so that we obtain polar decompositions
  \[ G = K \exp(\fp) \supeq H = H_K \exp(\fh_\fp).\]
  Now \cite[Thm.~IV.3.5]{Lo69} implies that the map
$K \times \fq_\fp \to M, (k,x) \mapsto k \Exp_{eH}(x)$ 
  factors through a diffeomorphism
  $K \times_{H_K} \fq_\fp \cong M$, thus defining on $M$ the structure of a
  vector bundle over $K/H_K$. 

  For any $\tau$ and $\theta$-invariant connected subgroup $S\subeq G$, we 
  have a symmetric subspace $M_S \cong S/H_S$, where $H_S = H \cap S$.
  Moreover $S = K_S \exp(\fp \cap \fs)$, and we obtain a
  diffeomorphism
  \[ K_S \times_{H_K \cap S} (\fq_{\fp} \cap \fs) \cong M_S.\]
  As the decompositions of $M$ and $M_S$ are compatible,
  the inclusion $M_S \into M$ is a proper map
  if the inclusion $K_S/(H_K \cap S) \into K/H_K$ is proper.

  \nin (b) If $f \: M_1 \to M_2$
  is a proper morphism of causal homogeneous spaces
  and $M_2$ is globally hyperbolic, then $M_1$ is also globally hyperbolic. 
\end{rem}

\begin{thm} \mlabel{thm:cc-emb} 
 Let  $M_G = G/H$ be a compactly causal symmetric  space
  corresponding to the irreducible modular compactly causal
  symmetric Lie algebra $(\g, \tau, C, h')$. 
  Then the connected components of $W_M^+(h')$ are not globally
  hyperbolic.
\end{thm}

\begin{prf} The dual non-compactly causal symmetric 
  Lie algebra $(\g^c, \tau^c, i C, h')$
  contains a causal Euler element $h'' \in \fp^c \cap i C^\circ$,
  which implies that $\tau_{h''}(h') = - h'$
  (recall that $\tau_{h''} =\theta^c \tau^c$). 
  Now the proof of \cite[Thm.~5.4]{MNO23} implies
  that $h'$ and $h''$ generate a 
  subalgebra $\fs^c \subeq \g^c$ isomorphic to~$\fsl_2(\R)$,
  which is invariant under $\tau^c$ and $\theta$. 
      So the inclusion $(\fs^c, \tau^c) \into (\g^c, \tau^c)$
    corresponds to an embedding of $2$-dimensional non-compactly causal
    symmetric spaces.     On the $c$-dual side,  we thus obtain an embedding 
    $(\fs,\tau) \into (\g,\tau)$ of compactly 
    causal symmetric Lie algebras, where $\fs \cong \fsl_2(\R)$
    is generated by the elliptic element $ih''$
    and the Euler element $h'$. Anti-de Sitter space 
    $\AdS^2$ is a global symmetric space corresponding to $(\fs,\tau)$.

        Proposition~\ref{prop:wads-incl} now implies that the 
    connected components of the positivity region
    $W_{M_S}^+(h') \subeq M_S$ are not globally hyperbolic.
    Here we use that the failure of global hyperbolicity
    of positivity regions is inherited by all covering spaces of $\AdS^2$
    because all connected components of $W_{M_S}^+(h')$ are diffeomorphic to
    $\R_+^2$, hence simply connected and therefore lift to coverings.
    Moreover, $K_S = \exp(\R i h'')$ is a closed subgroup of $K \subeq G$,
    so that we obtain with Remark~\ref{rem:proper} a proper map
    $M_S \to M$ on the level of causal symmetric spaces.
    As $W_{M_S}^+(h')$ is not globally hyperbolic, it contains a causal
    interval with non-compact closure in $M_S$, and then the corresponding
    interval in $M_G$ cannot have compact closure. 
\end{prf}

\begin{thm} \mlabel{thm:cc-emb2} 
  If the compactly causal symmetric space
  $M^{(\alpha)}$ is modular, then its positivity
  region is not globally hyperbolic.
\end{thm}

\begin{prf} If $M^{(\alpha)}$ is irreducible, then the assertion
  follows from Theorem~\ref{thm:cc-emb}.
  According to Table 4, this leaves us with the
  group case  $M^{(\alpha)} \cong \U_{s,s}(\C)$
  and the Lorentzian case~$\jV = \R^{1,d-1}$.

  In the latter case $M^{(\alpha)} \cong (\AdS^p \times \bS^q)/\{\pm \1\}$
  with $p\geq 2$ and $p + q = d \geq 3$ (cf.~\eqref{eq:mrj}). Then 
  the Euler element in $\fh^{(\alpha)}$ leaves a subspace isomorphic to
  $\AdS^2$ invariant. Hence the failure of global hyperbolicity
  for the components of $W_{\AdS^2}^+(h'')$ implies the same for 
  $W_{M^{(\alpha)}}^+(h'')$.

  For the group case $\U_{s,s}(\C)$, we have a causal embedding of
  \[ \U_{1,1}(\C) = (\bS^1 \times \SU_{1,1}(\C))/\{\pm \1\}
    \cong (\bS^1 \times \AdS^3)/\{\pm \1\} \supeq \AdS^3.\]
  Therefore the assertion follows from the failure of global
  hyperbolicity of positivity regions in~$\AdS^3$. 
\end{prf}

\section{Cayley transforms and group type spaces}
\mlabel{sec:8}

In this section we discuss some specific properties
of the group type spaces
\[ M^{(\alpha)} \cong \Sp_{2r}(\Omega), \U_{p,q}(\C), \SO^*(2r).\] 
They permit a unified approach via skew-hermitian forms over
$\K = \R,\C,\H$. The groups are simple for $\K = \R,\H$, but for
$\K = \C$ their center is a circle group.
We also  refer to \cite{KOe96} and \cite{KOe97} for a
uniform treatment of the $3$ series of group type spaces from a different
perspective (see also \cite[\S 1.6.1]{Be00}). 

\subsection{Unitary groups of skew-hermitian forms}

Let $\K\in \{\R,\C,\H\}$. We consider $\K^r$ is a $\K$-right module,
so that a skew-hermitian form on $\K^r$ is a map of the form
\[ \beta(z,w) = z^* B w, \quad \mbox{ where } \quad B \in \Aherm_r(\K).\]
Changing the basis by $g \in \GL_r(\K)$, leads to a change of
the representing matrix to $g^*Bg$. In this sense, the normal
forms are $B = \Omega_{2s}$ for $\K = \R$ ($r = 2s$ needs to be even),
$B = i I_{p,q}$, $p + q = r$, for $\K = \C$, and
$B = i \1_r$ for $\K = \H$ (\cite[p.~434]{Br85}).
This specifies the corresponding unitary Lie algebras
\[ \fu(B,\K^r) := \{ x \in \gl_r(\K) \: x^* B + B x = 0\},\]
and, concretely,
\[ \sp_{2s}(\R)  = \sp(\Omega_{2s},\R^{2s}), \quad 
 \fu_{p,q}(\C) = \fu(i I_{p,q}, \C^r), \quad  
 \so^*(2r) = \fu(i \1_r, \H^r).\]
The corresponding groups are $\Sp_{2s}(\R), \U_{p,q}(\C)$ and $\SO^*(2r)$.
These are precisely the groups among the
causal   Makarevi\v c spaces (Table 3). In this section we take a closer look
at these spaces and the corresponding geometric structures.

\subsection{Maximal isotropic Gra\ss{}mannians} 

On $E = \K^{2r}$ we consider the non-degenerate skew-hermitian form
\[ \beta(z,w) := z^* \Omega w, \quad \mbox{ where } \quad
  \Omega = \Omega_{2r}= \pmat{ 0 & \1_r \\ -\1_r & 0}.\]
Forms that can be represented in this way are precisely
those of Witt index $r$ (\cite{Br85}).
We write 
\[ G := \U(\Omega, E) 
  = \{ g \in \GL_{2r}(\K) \:
  (\forall z,w) \ \beta(gz,gw) = \beta(z,w) \} 
  = \{ g \in \GL_{2r}(\K) \: g^* \Omega g = \Omega \} \]
for the corresponding isometry group and
$M \subeq \Gr_r(\K^{2r})$ for the space of maximal isotropic
subspaces $L \subeq E$. Note that $G$ acts transitively
on $M$ by Witt's Theorem (\cite{Br85}). 

Euler elements $h \in \fu(\Omega, E)$ 
correspond to orthogonal decompositions
\begin{equation}
  \label{eq:decomp}
  E = E_+ \oplus E_-
\end{equation} 
into isotropic subspaces
$E_\pm = \ker(h \mp \shalf \1)$. Then the parabolic subgroup
$Q_h \subeq G$ (cf.\ Section~\ref{sec:2}) is the stabilizer
of $E_-$, so that $M \cong G/Q_h$ is the associated causal flag manifold.

\nin {\bf The tangent space of $M$:}
The decomposition \eqref{eq:decomp}
identifies the tangent space $T_{E_-}(M)$ naturally with
the set of all $\K$-linear maps
$\phi \: E_- \to E_+$ whose graph
\begin{equation}
  \label{eq:gammaphi}
  \Gamma(\phi) = \{ v + \phi(v) \: v \in E_-\}
\end{equation}
is isotropic. This means that
$\beta(v,\phi(w)) + \beta(\phi(v), w) = 0$ for $v, w \in E_-,$ 
hence  that
\[ \gamma_\phi(v,w) := \beta(v,\phi(w))
  = - \oline{\beta(\phi(w), v)}
  = \oline{\beta(w,\phi(v))} \]
defines a hermitian form on $E_-$. We thus obtain a natural
identification
\[ T_{E_-}(M) \cong \jV := \Herm_r(\K),\]
an embedding $\phi \mapsto \Gamma(\phi)$ into~$M$, 
and also a causal structure by the convex cone
of positive semi-definite forms on $E_-$.\\

\nin {\bf Fractional linear transformations:}
Any matrix $g = \pmat{ a & b \\ c & d} \in \U(\Omega, \K^{2r})$
(with the block structure according to $E = E_+ \oplus E_-$ as in \eqref{eq:decomp})
satisfies 
\begin{equation}
  \label{eq:fraclin2}
  g.\Gamma(\phi) = \Gamma((a \phi + b)(c \phi + d)^{-1}),
\end{equation}
so that $\U(\Omega, \K^{2r})$ acts naturally by fractional linear transformations 
on $\Herm_r(\K)$.\\

\nin {\bf Jordan automorphisms:}
On the Jordan algebra $\jV = \Herm_r(\K)$, we consider an 
automorphism $\alpha(x) = J x J^{-1}$, where  $J \in \U_r(\K)$
satisfies $J^* = -J$ and is given by: 
\begin{equation}
  \label{eq:thejs}
 J =
  \begin{cases}
    \Omega_{2s} \text{\ for\ }  \K= \R, r= 2s, \\
      i I_{p,q} \text{\ for\ } \K = \C, r = p + q, \\ 
      i \1_r \text{\ for\ }  \K = \H.
  \end{cases}
\end{equation}
Then $\sigma_\alpha(x) := J_d x J_d^{-1}$, for $J_d := J \oplus J$, is
the corresponding automorphism of $\g= \fu(\Omega,\K^{2r})$ and
thus $\theta_\alpha(g) = J_d (g^*)^{-1} J_d^{-1}$. The fixed point
group $G^{(\alpha)} = G^{\theta_\alpha}$ then consists of those matrices
commuting with the matrix
\[  J_d \Omega= \pmat{ 0 & J \\ - J & 0},\]
which is an involution. Its eigenspaces are 
\[ \ker( J_d \Omega\mp \1) = \Gamma(\pm J) = \{ (v, \pm J v) \: v \in \K^r \}, \]
and the restriction of $\beta$ to these spaces is represented
by the matrix $\pm 2 J$. We conclude that,
for $(F,\beta_F)$ with $\beta_F(v,w) = v^* 2 J w$, we have
\begin{equation}
  \label{eq:ebeta-deco}
  (E,\beta) \cong (F,\beta_F) \oplus (F,-\beta_F).
\end{equation}
For $g \in \GL(F)$, the graph $\Gamma(g)$ is isotropic
if and only if $g \in \U(F,\beta_F)$. We thus obtain an embedding
of the corresponding groups 
\begin{equation}
  \label{eq:unitary-emb}
 \Gamma\:  \U(F,\beta_F) \to M, \quad
 g \mapsto \Gamma(g).
\end{equation}
Then \eqref{eq:fraclin2} still holds with respect to the
matrix block structure defined by \eqref{eq:ebeta-deco}, i.e.,
$g = \pmat{a & b \\ c & d} \in \U_{2r}(\Omega, \K^{2r})$
acts naturally by fractional linear maps on the unitary group~$\U(F,\beta_F)$.

The decomposition \eqref{eq:ebeta-deco} of $E$ also leads to an inclusion
$\U(F,\beta_F)^2 \into \U(E,\beta), (g_1, g_2) \mapsto g_1 \oplus g_2.$ 
As $(g_1, g_2)\Gamma(\phi)  = \Gamma(g_2 \phi g_1^{-1}),$ 
the map $\Gamma$ in \eqref{eq:unitary-emb} is equivariant for the natural action of
$\U(F,\beta_F)^2$ on $\U(F,\beta_F)$ by
$(g_1, g_2).g = g_2 g g_1^{-1}$.

This leads to the desired realizations of three families of groups
as a space $M^{(\alpha)}$: 
\begin{itemize}
\item[\rm($\R$)] $\Sp_{2s}(\R)$ for $\jV = \Sym_{2s}(\R)$. 
\item[\rm($\C$)] $\U_{p,q}(\C)$ for $\jV = \Herm_r(\C)$, for
  $r = p + q$. For $q = 0$, we have a diffeomorphism
  $\U_r(\C) \to M$. 
\item[\rm($\H$)]  $\SO^*(2r)$ for $\jV = \Herm_r(\K)$.
\end{itemize}

To determine the corresponding biinvariant causal structures on
these groups, we observe that, for $A \in \End(F)$, the sesquilinear form
$\beta_F^A(z,w) := \beta_F(z,Aw)$ 
  on $F$ is hermitian if and only if $A^* = - A$ with respect to $\beta_F$.
This defines a natural bijection $\fu(F,\beta_F) \cong \Herm(F)$
  and also a natural specification of the invariant cone in $\fu(F,\beta_F)$ by
  \[  C := \{ A \in \fu(F,\beta_F) \: \beta_F^A \geq 0\}.\]
  This cone determines the biinvariant causal structure on $\U(F,\beta_F)$
  for which the embedding is a causal map. 
  Concretely, we have
  \begin{itemize}
  \item[($\R$)] $C = \{ A\in \sp_{2s}(\R) \: \Omega_{2s} A \geq 0 \}$,
    corresponding to the positive semidefinite forms on $\R^{2s}$.
  \item[($\C$)] $C = \{ A\in \fu_{p,q}(\C)  \: i I_{p,q} A \geq 0 \}$
    is a minimal invariant cone in $\fu_{p,q}(\C)$
     (\cite[\S 8.5]{HN93}).
   \item[($\H$)] $C = \{ A\in \so^*(2r) = \fu(i\1_r,\H^r) \: i A \geq 0 \}$.
  \end{itemize}

  \begin{rem} Conformal embeddings into flag manifolds defined by
    Euler elements exist also for other groups, 
but these spaces are not causal. In particular $\SO_{2,d}(\R)$ embeds
into non-Lorentzian quadrics (see \cite{GKa98} for more details).
\end{rem}

\begin{rem} (Lorentzian spaces of group type) 

  \nin (P) The group $\U_n(\C) \cong (\T \times \SU_n(\C))/C_n$
  carries a biinvariant Lorentzian structure, 
  but the causal structure on $\U_n(\C)$,
  as the conformal completion of $\Herm_n(\C)$,
  corresponds to the positive cone $\Herm_n(\C)_+$
  which is only Lorentzian for $n = 2$.
  So we only obtain natural Lorentzian structures on the groups
  $\U_{p,q}(\C)$ for $p + q = 2$.

\nin (NS1) The group $\Sp_{2s}(\R)$ is Lorentzian of dimension $3$ for $s = 1$. 
Note that $\SL_2(\R)\cong \Sp_2(\R) \cong \SU_{1,1}(\C) \cong \AdS^3$.

\nin (S2) The group $\SO^*(2r)$ is Lorentzian of dimension $6$ for $r = 1$.

\end{rem}

\begin{rem}
For Pierce involutions $\alpha$ on $\Herm_r(\K)$,
the group $G^{(-\alpha)}$ commutes with some Euler element,
hence preserves a pair of complementary isotropic subspaces,
which leads to
\[  G^{(-\alpha)} \cong \{ (g_1, g_2) \in \GL(E_-) \times \GL(E_+) \:
  g_2 = g_1^{-*} \} \cong \GL(E_-) \cong \GL_r(\K) \]
and $M^{(-\alpha)} \cong \Herm_r(\K)^\times_{(p,q)}.$ 
\end{rem}

\begin{rem}
For $(E,\beta) \cong (F,\beta_F) \oplus (F,-\beta_F)$, the subspaces 
  \[ E_\mp := \Gamma(\pm \1) = \{ (v,\pm v) \:  v \in F \} \cong F \]
  are maximal isotropic with $E = E_+ \oplus E_-$, which can be used
  to specify the Euler element $h = \shalf \diag(\1,-\1)$ in $\fu(\Omega, \K^{2r})$
\end{rem}

\subsection{Relations to Cayley transforms}

In this subsection we discuss charts of the group type spaces
coming from Cayley transforms.
The well-known Cayley transform $\Herm_r(\C) \to \U_r(\C)$ 
arises in our context by mapping $\Gamma(z) \in M$ for
$z \in \Herm_r(\C)$ to the corresponding element of $\U_r(\C)$. From
\[ (zv,v)
  = \Big(\frac{zv - i v}{2}, i \frac{zv - iv}{2}\Big)
  + \Big(\frac{zv + i v}{2}, -i \frac{zv + iv}{2}\Big),\]
it follows that $\Gamma(z)$ corresponds to the unitary matrix
\begin{equation}
  \label{eq:tildecz}
  \tilde C(z) = (z + i \1)(z - i \1)^{-1}.
\end{equation}

More generally, we write 
\[ (zv,v)
  = \Big(\frac{zv - J v}{2}, J \frac{zv - Jv}{2}\Big)
  + \Big(\frac{zv + J v}{2}, -J \frac{zv + Jv}{2}\Big),\]
and find the Cayley transform $\Herm_r(\K) \to \U(J,\K^r)$ 
\begin{equation}
  \label{eq:tildecz2}
  \tilde C(z) := (z + J)(z-J)^{-1}.
\end{equation}

\section{Perspectives and open problems} 
\mlabel{sec:perspec}

\subsection{Open problems of general nature} 

\begin{prob} Find an explicit description of the positivity
  regions $W_M^+(h_\lambda)$ for the Euler elements
  $h_\lambda = h^j + \lambda h$ in $\g^{(-\alpha_j)}$
  for the Pierce type involutions
  $\alpha_j$ (cf.\ Proposition~\ref{prop:4.9}).
  For the modular cases \cite[Prop.~6.1]{NO23a} may be useful.
\end{prob}

\begin{prob} Let $M^{(\alpha)}$ be a modular compactly causal symmetric space
  (cf.~Theorem~\ref{thm:6.5}).
  \begin{itemize}
  \item   Describe the wedge regions $W := W_{M^{(\pm \alpha)}}^+(h')$  more concretely (cf.~\cite{NO23a, NO23b}).
  \item  Determine the fundamental groups $\pi_1(W)$ of the
    wedge regions. 
    They are never globally hyperbolic by Theorem~\ref{thm:cc-emb2}. 
\item Does the identity component of $G_W:= \{ g \in G \: g.W = W\}$ 
  commute with the modular flow, i.e., $G_{W,e} \subeq G^{h'}$?
  Note that the Euler element $h'$ is also one in $\g$.

  For flip involutions,  
  $W_{M_d^{(\alpha)}}^+(h_d)$ is dense in $W_{M_d}^+(h_d)$
  (Proposition~\ref{prop:ct-posreg}(b)), so that 
  $G_{W,e}$ preserves $W_{M_d}^+(h_d) = \jV_+ \times (-\jV_+)$,
  hence  is contained in $G^h_e \times G^h_e$
  and thus commutes with~$h_d$.

\item When does $G_W$ act transitively on $W$?
  For the non-compactly causal spaces, this follows from
  Theorem~\ref{thm:pos-dom-ncc}, but for the compactly causal spaces
  it is false in general.   For $\AdS^2 \cong M_d^{(\alpha)}$ with
  $\jV = \R$
  and $\jV_d \cong \R^{1,1}$,   Example~\ref{ex:cayley-ads} implies that
  $G_{W,e} \subeq G_{W_{M_d}^+(h_d),,e} \cong \R_+^2$ is a proper closed subgroup,
  hence equal to $\exp(\R h_d)$. Therefore it cannot act transitively
  on the $2$-dimensional domain~$W$. 
\end{itemize}
\end{prob}

\begin{prob} Characterize those spaces $M^{(-\alpha)}$ which are
  dense in $M$.
  \begin{itemize}
  \item For Cayley type (C), this is always the case. 
  \item For Pierce type (P), this is never the case. 
  \item In the Lorentzian case $M^{(-\alpha)}$ is not dense
    (Proposition~\ref{prop:boundaries}).
  \end{itemize}
\end{prob}

Here is a positive result for compactly causal spaces: 
\begin{thm} The compactly causal subspaces $M^{(\alpha)}$
  are dense in $M$.   
\end{thm}

\begin{prf} This follows from \cite[Thm.~3.3.6]{Ber98}.
It can also be derived from Theorems 5.1 and 5.9 in \cite{Bet03},
which unfortunately uses the defective tables in \cite{Ma73}.
As a consequence, Betten misses the compactly causal space
of type (S3), corresponding to $\g^{(\alpha)} = \su_{6,2}(\C)$
and $\fh^{(\alpha)} \cong = \fu_{3,1}(\bH)$.
\end{prf}

\subsection{Causal convexity  of wedge regions}
\mlabel{sec:9}

\begin{defn} \mlabel{def:caus-convex}
  We call an open subset $\cO \subeq \jV$ of the euclidean Jordan algebra $\jV$ 
  {\it causally convex} if, for $a,b \in \cO$, the corresponding
double cone
\[ \cD_{a,b} := (a - \jV_+) \cap (b + \jV_+) \]
is contained in $\cO$.
\end{defn}

\begin{prob} \mlabel{prob:9.6} Are connected open subsets 
  $\cO \subeq \jV$ which are globally hyperbolic also causally convex?

  Causal convexity implies global hyperbolicity:
  For $a,b \in \cO$, pick $\eps > 0$ with
  $a + \eps \be_0$, $b - \eps \be_0 \in \cO$. Then the closure
  of the double cone   between $a$ and $b$ is contained in
$\cD_{a + \eps \be_0, b - \eps \be_0} \subeq \cO,$ 
  so that $\cO$ is globally hyperbolic.

Olaf M\"uller has an example 
of an open connected subset of $\R^{1,2}$ which is globally
hyperbolic but not causally convex (\cite{Mu24}).
It is a union of double
cones, aligned along a spiral with timelike axis. 
\end{prob}

\begin{ex}
  We have seen in Proposition~\ref{prop:ct-posreg} that, 
  for flip involutions (Cayley type), the wedge region in
  $M^{(-\alpha)}$ has a causally convex realizations as $\jV_+ \times -\jV_+
  \subeq\jV \times \jV$. 
\end{ex}

\begin{prob}
  Show that the open domains $\jV^\times_{(p,q)}$, the connected
  components of $\jV^\times$, are causally convex in $\jV$.
  We expect that, for $a,b$ in this domain with $b \in a - \jV_+$, 
  the pair $a,b$ is contained in a wedge region $W \subeq \jV^\times_{(p,q)}$,
  and this in turn contains the corresponding double cone~$\cD_{a,b}$.

  For $r = 2$, the critical case is $\jV^\times_{(1,1)}$, which is the 
  set of spacelike elements. If $a,b$ are spacelike and $\cD_{a,b} \not=\eset$,
  then $(a - \jV_+) \cap \oline{\jV_+} = \eset$ and 
  $(b + \jV_+) \cap -\oline{\jV_+} = \eset$ imply
  $\cD_{a,b} \subeq \jV^\times_{(1,1)}$. 
\end{prob}

\begin{prob} In the modular case the Euler element
  $h' \in \fh^{(\pm \alpha)}$ is also one in $\g$
  (Theorem~\ref{thm:6.5}), hence conjugate to $h$.
  So some conjugate of the wedge region $W = W_{M^{(\pm \alpha)}}(h')$
  is contained in $\jV$.
  When is there a $G$-transform $g.W$ which is causally convex in $\jV$?
  This is false for anti-de Sitter space by Proposition~\ref{prop:wads-incl}
  because this requires $W$ to be globally hyperbolic
  (cf.~Problem~\ref{prob:9.6}).

  By Theorem~\ref{thm:cc-emb2}, this is never the case for
  the compactly causal spaces~$M^{(\alpha)}$.
  For flip involutions and Pierce involutions,
  Propositions~\ref{prop:ct-posreg} and \ref{prop:pierce-concrete} 
  shows that it is true for $M^{(-\alpha)}$. 
  \end{prob}

\begin{prob}
Irreducible non-compactly causal
symmetric spaces $M = G/H$ possess maximal and minimal causal
structures. For the maximal one,
\cite[Lemma~5.2, Thm.~5.7]{MNO24} 
implies the causal convexity of the wedge regions
in the globally hyperbolic space $M$,
hence that wedge regions are globally hyperbolic.
What about the other invariant causal structures?
Are their wedge regions also globally hyperbolic?
Note that, for a given Euler element, the wedge region
depends on the causal structure.

The discussion of the group type spaces $M^{(\alpha)}$ in Section~\ref{sec:8}
shows that the causal structure need not be maximal, which is inherited
by the non-compactly causal duals $M^{(-\alpha)}$. Therefore
the arguments in \cite{MNO24} do not apply to these spaces. 
\end{prob}

\begin{prob} Let $G$ be a simple hermitian Lie group,
  $C_\g := C_\g^{\rm max}$ a maximal invariant cone in~$\g$,
  and $z_\fk \in \fz(\fk) \cap (C_\g^{\rm max})^\circ$ the unique element for 
  which $\ad z_\fk$ defines a complex structure on $\fp$.
  Then $C_\fk := C_\g^{\rm max} \cap \fk$ is an
  $\Ad(K)$-invariant cone satisfying
  $C_\g = \Ad(G) C_\fk$, $z_\fk \in C_\fk^\circ$, and
  \[ \cD_\fk := C_\fk^\circ \cap (2\pi z_\fk - C_\fk^\circ) \]
  is a causal interval in $\fk$.  Further 
  \[ \cD_\g :=  \Ad(G)\cD_\fk   \subeq C_\g^\circ \]
is an open elliptic domain in $\g$
mapped by the exponential function diffeomorphically
to a domain $\cD_G \subeq G$. In \cite[Thm.~6.1]{HN24} it is shown that
this domain is globally hyperbolic for the biinvariant causal
structure on $G$.

This domain lies in a compactly causal symmetric space~$M$.
Its structure resembles the domains we find
by considering for the dual non-compactly causal
symmetric space $M^c$  the open subset
\[ \Xi \cap M =\Exp_{e}(\Omega_{\fq}),\]
where $\Omega_{\fq} \subeq C_\fq$ is a suitable open subset
and $\Xi$ is a crown domain for $G^c$ (cf.\ \cite{NO23b}).
Are the regions $\Xi \cap M$ in the causal homogeneous space $M$
globally hyperbolic? 
\end{prob}

\subsection{Non-reductive spaces}

In this paper we only studied open orbits
of symmetric subgroups of the conformal group~$G= \Co(\jV)_e$.
Using \cite{BB10}, it should also be possible
to develop a similar theory for open orbits of 
non-reductive subgroups.  An important example arises from the group
$G_1 \rtimes G_0 = \exp(\g_1)G^h_e$ acting on~$\jV$ by affine maps.
Specifically, the Poincar\'e group $\R^{1,d-1} \rtimes \SO_{1,d-1}(\R)_e$
acts with open orbits on Minkowski space and its compactification
$Q(\R^{2,d})$. 

\subsection{Holography}

The open orbits $M^{(\pm \alpha)} \subeq M$ have natural
compactifications obtained by their closure in~$M$.
Can these embeddings be used to develop
``holographic'' techniques on the level of nets of real subspaces
on $M^{(\pm \alpha)}$ and its boundary orbits?
We refer to \cite{Re00} for the relations between
  wedge regions in anti-de Sitter space $\AdS^d$ and regions in Minkowski space
  ``at infinity''  by light rays.
  For more on holographic aspects, we refer to \cite{dB01}, \cite{St01}
  (which deals with de Sitter space),
 \cite{BCW25}, and \cite{Wo24}.

Describe the boundary orbits under $G^{(\pm \alpha)}$ which are causal.
Here the closed orbits should be of particular interest
because they should be causal flag manifolds of $G^{(\pm \alpha)}$.
For the embedding $\AdS^d \into Q(\R^{2,d})$, 
the boundary  is a copy of $Q(\R^{2,d-1})$, compactified Minkowski space
  of codimension~$1$ (Proposition~\ref{prop:boundaries}).

\subsection{More connections with Physics?} 
  
\cite{KK23} and \cite{KKR25} deal with super versions
of the embedding of $\AdS^d$ in $Q(\R^{2,d})$.
In this context the existence of such embeddings
characterizes some ``superconformal flatness'';
see also  \cite{BC18}, \cite{FL15}, \cite{Gu00, Gu01}. 

We also mention \cite{Pa04}, which develops a nice picture
relating causal structures on manifolds and positive energy conditions
for representations. It touches on topics such as horizons,
supergravity and string theory, from a group theoretic perspective.
In the non-super context, some of this is under development in~\cite{MN25}.

\appendix

\section{Essential tools} 

\subsection{Some facts on convex cones} 
  
\begin{lem} \label{lem:coneint} {\rm(\cite[Lemma~B.1]{MNO23})}
Let $E$ be a finite dimensional real vector space, 
$C \subeq E$ a closed convex cone and $E_1 \subeq E$ a linear subspace. 
If the interior $C^\circ$ of $C$ intersects $E_1$, then 
$C^\circ \cap E_1$ coincides with the relative interior $C_1^\circ$ 
of the cone $C_1 := C \cap E_1$ in~$E_1$.
\end{lem}

\begin{lem} \mlabel{lem:ext-eigenvalue}
  Let $V$ be a finite-dimensional real vector space,
  $A \in \End(V)$ diagonalizable, and let $C \subeq V$ be a
  closed convex cone invariant under $e^{\R A}$.
 Let  $\lambda_{\rm min}$ and  $\lambda_{\rm max}$ be the minimal/maximal
 eigenvalues of $A$. For an eigenvalue $\lambda$ of $A$ we write
 $V_\lambda(A)$ for the corresponding eigenspace and
 $p_\lambda \:  V \to V_\lambda(A)$ for the projection along all other
 eigenspaces. Then
 \begin{equation}
   \label{eq:plambda}
 p_{\lambda_{\rm min}}(C) = C \cap V_{\lambda_{\rm min}}(A)
   \quad \mbox{ and }  \quad 
   p_{\lambda_{\rm max}}(C) = C \cap V_{\lambda_{\rm max}}(A).
 \end{equation}

 If $A$ has only two eigenvalues, it follows that
$C = p_{\lambda_{\rm min}}(C) \oplus p_{\lambda_{\rm max}}(C).$ 
\end{lem}

\begin{prf} Since we can replace $A$ by $-A$, it suffices to verify the
  second assertion in \eqref{eq:plambda}.
  So let $v \in C$ and write it as a sum
  $v = \sum_\lambda v_\lambda$ of $A$-eigenvectors. Then
  \[ v_{\lambda_{\rm max}} = \lim_{t \to \infty} e^{-t \lambda_{\rm max}}
    e^{tA} v \in C \]
  implies that $p_{\lambda_{\rm max}}(C) \subeq C \cap V_{\lambda_{\rm max}}(A),$
  and the other inclusion is trivial.   
\end{prf}

\subsection{$\fsl_2$-data} 
\mlabel{app:a.3}

\begin{example} \mlabel{ex:sl2a}
  In $\fsl_2(\R)$, we fix the following notation.
  We consider the Cartan involution $\theta(x) = - x^\top$
  with $\fk = \so_2(\R)$. 
  For the Euler elements
  \[  h := \shalf\diag(1,-1) \quad \mbox{ and } \quad
    k:= \shalf\pmat{0 & 1 \\ 1 & 0}\]
  we then have
  \begin{equation}
    \label{eq:kande}
    k = \frac{1}{2}(e - \theta(e))\quad \mbox{ and }\quad
    h = \frac{1}{2}[e,\theta(e)]] 
    \quad \mbox{ for } \quad
    e = \pmat{0 & 1 \\0 & 0}.
  \end{equation}
Further, 
\[ z_\fk := \frac{1}{2} \pmat{ 0 & 1 \\ -1 & 0} = [h, k]
  \quad \mbox{ satisfies  } \quad
  [z_\fk, h] = \frac{1}{2} \pmat{ 0 & -1 \\ -1 & 0} =- k, \]
    so that we have 
    \begin{equation}
      \label{eq:90deg-rot}
e^{-\frac{\pi}{2} \ad z_\fk}h= - [z_\fk,h]  = \frac{1}{2} \pmat{ 0 & 1 \\ 1 & 0}= k.
    \end{equation}
\end{example}

\subsection{Strongly orthogonal roots}
\mlabel{app:strong-orth}

Let $\g$ be a simple real Lie algebra and 
$h \in \g$ an Euler element. 
Let $\theta$ be a Cartan involution with $\theta(h) = -h$ and 
$\fa \subeq \fp := \g^{-\theta}$ be a maximal abelian subspace  
containing~$h$. In 
\[ \Sigma_1 := \{ \alpha \in \Sigma(\g,\fa) \: \alpha(h) = 1\} \]
we pick a  maximal set $\{ \gamma_1, \ldots, \gamma_r\}$ of long roots
which are {\it strongly orthogonal}, i.e., 
neither $\alpha + \beta$ nor $\alpha - \beta$ is a root.
We refer to \cite[\S 3.2]{MNO23} and \cite[Thm.~2.6]{MNO25}
for more on systems of strongly orthogonal roots.

Then there exist elements $c_j \in \g_{\gamma_j}$ such that 
the Lie subalgebras 
\begin{equation}
  \label{eq:def-fsj}
 \fs_j := \Spann_\R \{c_j, \theta(e_j), [c_j, \theta(c_j)]\}
  = \g_{\gamma_j} + \g_{-\gamma_j} + \R \gamma_j^\vee, \quad 
  j =1,\ldots, r,
\end{equation}
are isomorphic to $\fsl_2(\R)$, so that
\begin{equation}
  \label{eq:def-fs}
  \fs := \sum_{j = 1}^r \fs_j\cong \fsl_2(\R)^{\oplus r}.
\end{equation}


\begin{prop} \mlabel{prop:4.9}
  If $\g$ is hermitian, then the Euler elements 
  \begin{equation}
  \label{eq:hj-decomp}
  h^j := \frac{1}{2} (\gamma_1^\vee + \cdots + \gamma_{r-j}^\vee - \gamma^\vee_{r-j+1}  - \cdots - \gamma_r^\vee), \quad j = 0,\ldots, r, 
\end{equation}
represent the $\Inn(\g_0(h))$-conjugacy classes in $\cE(\g) \cap \g_0(h)$.
The $\Inn(\g_0(h))$-conjugacy classes of Euler elements of $\g_0(h)$
are represented by the elements of the form
\[ \lambda h + h^j, \quad \lambda \in \R, j = 1,\ldots, r-1.\]
\end{prop}

Note that $h  = h^0 = - h^r$.

\begin{prf}   The first assertion follows from~\cite[Thm.~3.2]{NO23a}.
Now we turn to Euler elements of the subalgebra $\g_0(h)$.
An element $x = \frac{1}{2}\sum_{j = 1}^r x_j \gamma_j^\vee \in \fa$
is an Euler element in the subalgebra $\fh = \g_0(h)$ 
with the root system
\[ \Sigma(\fh,\fa) = \big \{ \shalf(\gamma_j - \gamma_k) \: j \not=k\big\}
  \cong A_{r-1} \]
if and only if
\begin{equation}
  \label{eq:halfcond}
  (\forall j \not=k) \quad \shalf(x_j - x_k) \in \{  0, \pm 1\}.
\end{equation}
Subtracting a real multiple of the central element $h \in \g_0(h)$
and acting with the Weyl group, we may assume that
$1 = x_1 \geq x_2 \geq \cdots \geq x_r$. 
Then \eqref{eq:halfcond} implies that there exists a $j$ with 
$x_1 = \cdots = x_j = 1$ and $x_{j+1} = \cdots = x_r = -1,$ 
i.e., $x= h^j$. This proves the assertion.
\end{prf}

\subsection{The general stereographic projection}
\mlabel{subsec:stereo} 

For a finite-dimensional linear space $\jV$, endowed with a
non-degenerate symmetric bilinear form $\beta$, we consider the extensions
\begin{equation}
  \label{eq:tildev-deco}
 \jV^\wedge := \R \oplus \jV
  \subeq \tilde \jV := \R \oplus \jV \oplus \R 
  \supeq \jV^\vee := \jV \oplus \R
\end{equation}
with
\[   \tilde\beta((t,v,s), (t',v',s')) = tt' + \beta(v,v) - ss' \quad
  \quad \mbox{ and }\quad
  \beta^\wedge :=  \tilde\beta\res_{\jV^\wedge \times \jV^\wedge}.\]

We write
\[ S(\jV^\wedge, \beta^\wedge) := \{ v \in \jV^\wedge \: \beta^\wedge(v,v) = 1\} \]
for the unit sphere in $\jV^\wedge$.
Then $-\be_1 = (-1,0,0)\in S(\jV^\wedge, \beta^\wedge)$ and a
line $-\be_1 + \R(v + \be_1)$, $v \in \jV$, is tangent to
$S(\jV^\wedge, \beta^\wedge)$ if and only if $\beta(v,v) = 0$.
All other lines of this form intersect $S(\jV^\wedge, \beta^\wedge)$ in 
exactly one other point $\eta_1(v)$, provided $\beta(v,v) \not=-1$.
This leads to a map
\begin{equation}
  \label{eq:def-eta1}
  \eta_1 \: \jV \setminus S(\jV,-\beta) \to S(\jV^\wedge, \beta^\wedge),\quad
  \eta_1(v) := \frac{1- \beta(v,v)}{1 + \beta(v,v)} \be_1 +
  \frac{2}{1 + \beta(v,v)} v,
\end{equation}
whose inverse is given by
\[ \eta_1^{-1} \: S(\jV^\wedge, \beta^\wedge) \setminus \{ w_1 \not= -1\} \to \jV,
  \quad
  w \mapsto \frac{\bw}{1 + w_1}\quad \mbox{
    for } \quad w = (w_1, \bw) \in \jV^\wedge.\]
This map is the {\it stereographic projection}
of $S(\jV^\wedge, \beta^\wedge)$ to $\jV$.
If
\[ \eta \: \jV \to 
  Q(\tilde \jV) := \{ [\tilde v] \in \bP(\tilde \jV) \:
  \tilde\beta(\tilde v,\tilde v) = 0\} \] 
is the natural embedding
\[ \eta \: \jV \to Q, \quad  \eta(v) := \Big[ \frac{1 - \beta(v,v)}{2} : v : - \frac{1 + \beta(v,v)}{2}\Big] 
  \in \bP(\tilde\jV) \] 
(cf.~also \eqref{eq:eta-flow}), then 
$\eta \circ \eta_1^{-1}$ extends to a map 
$S(\jV^\wedge, \beta^\wedge) \to Q(\tilde \jV), \hat v \mapsto [\hat v:-1]$ 
that induces a diffeomorphism
\begin{equation}
  \label{eq:conformal-flat}
\eta_1^{-1} \:  S(\jV^\wedge, \beta^\wedge) \setminus \{w_1 \not=-1\} \to
 \jV \setminus S(\jV,-\beta).
\end{equation}

\end{document}